\documentclass[upref]{amsart}
\usepackage{latexsym,amsmath,amstext,mathrsfs,amssymb,amsthm,amscd,amsopn,verbatim}

\usepackage{times}
\numberwithin{equation}{section}

\newcommand{\abs}[1]{\lvert#1\rvert}
\DeclareMathOperator{\obsmu}{\obs_{\boldsymbol{\mu}}}
\DeclareMathOperator{\supercon}{\sfC_{\boldsymbol{\lambda}}}
\DeclareMathOperator{\sP}{P}

\newcommand{\bde}{{\boldsymbol{\partial}}}
\newcommand{\rd}{\overleftarrow{\partial}} 
\newcommand{\ld}{\overrightarrow{\partial}} 
\newcommand{\da}{\partial \phi^{\alpha}}
\newcommand{\sqbr}[2]{{\left\langle{\,{#1}\,,\,{#2}\,}\right\rangle}_{\mathrm{Hodge}}}
\newcommand{\cirbr}[2]{{\left({\,{#1}\,,\,{#2}\,}\right)}_{\mathrm{Hodge}}}
\newcommand{\lmf}[1]{\frac{\overrightarrow{\delta} #1}{\delta \phi^\alpha}}
\newcommand{\rmf}[1]{\frac{#1 \overleftarrow{\delta}}{\delta \phi^\alpha}}
\newcommand{\lma}[1]{\frac{\overrightarrow{\delta} #1}{\delta \phi_\alpha^*}}
\newcommand{\rma}[1]{\frac{#1 \overleftarrow{\delta}}{\delta \phi_\alpha^*}}
\newcommand{\lmalpha}[1]{\frac{\overrightarrow{\delta} #1}{\delta \varphi_\alpha}}
\newcommand{\rmalpha}[1]{\frac{#1 \overleftarrow{\delta}}{\delta \varphi_\alpha}}

\DeclareMathOperator{\orient}{or}
\newcommand{\ii}{{\mathrm{i}}}
\DeclareMathOperator{\ev}{ev}
\newcommand{\ser}{\mathcal{S}_{\sfA,\sfB}}
\newcommand{\dad}{\partial \phi^{+}_{\alpha}} 
\newcommand{\der}{\partial \varphi_{\alpha}} 
\newcommand{\lA}[1]{\frac{\ld #1}{\partial \sfa}} 
\newcommand{\lB}[1]{\frac{\ld #1}{\partial \sfB}} 
\newcommand{\rA}[1]{\frac{#1 \rd}{\partial \sfa}} 
\newcommand{\rB}[1]{\frac{#1 \rd}{\partial \sfB}}\newcommand{\rfunc}[2]{\frac{#1 \rd}{\partial #2}}
\newcommand{\lfunc}[2]{\frac{\ld #1}{\partial #2}}
\newcommand{\Abs}[1]{\Big\lvert#1\Big\rvert}
\newcommand{\dd}{\mathrm{d}} 
\newcommand{\ddtilde}{{\tilde \dd}}
 
\DeclareMathOperator{\tr}{Tr}
\DeclareMathOperator{\obs}{\sfO}
\DeclareMathOperator{\Ug}{\mathcal{U}(\Lg)}

\DeclareMathOperator{\gh}{gh} 
\DeclareMathOperator{\bvec}{\boldsymbol{v}} 
\DeclareMathOperator{\dg}{deg} 
\DeclareMathOperator{\ddt}{\frac{d}{d\mathnormal{t}}
\Big\vert_{\mathnormal{t}=0}}
\newcommand{\holo}[3]{H(#1)\arrowvert_{#2}^{#3}}
\DeclareMathOperator{\sfhol}{{\mathsf{Hol}}}   
\newcommand{\sfholo}[1]{\sfhol(#1)}

\newcommand{\dA}{\dd_\sfA} 
 
\DeclareMathOperator{\ad}{ad}
\DeclareMathOperator{\sad}{\mathsf{ad}}
\DeclareMathOperator{\Ad}{Ad}
\DeclareMathOperator{\rk}{rk}
\newcommand{\lo}{\mathrm{L}}
\newcommand{\pser}{\mathcal{S}_{\sfA, \sfB}}
\newcommand{\pserr}[1]{\mathcal{S}_{\sfA, \sfB}^{#1}}

\newcommand{\Lg}{\mathfrak{g}}
\newcommand{\push}[2]{\pi_{1*} \dbraket{#1}{\pi_2^{*}#2}}
\newcommand{\pushw}[2]{\pi_{1*} \braket{#1}{\pi_2^{*}#2}}
\DeclareMathOperator{\Aut}{Aut}
\DeclareMathOperator{\End}{End}
\DeclareMathOperator{\Hom}{Hom}
\newcommand{\en}[1]{\End(#1)}
 
\newcommand{\sbv}[2]{\left(\!\left({\,{#1}\,;\,{#2}\,}\right)\!\right)}
\newcommand{\lb}[2]{[\![#1\,;#2]\!]}
\newcommand{\covar}[2]{\dd_{#1}#2}
 
\newcommand{\lfuna}[1]{\frac{\ld #1}{\da}} 
\newcommand{\lfunad}[1]{\frac{\ld #1}{\dad}} 
\newcommand{\rfuna}[1]{\frac{#1 \rd}{\da}} 
\newcommand{\rfunad}[1]{\frac{#1 \rd}{\dad}}
\newcommand{\integ}[2]{\int_M \dbraket{#1}{#2}}
\newcommand{\integw}[2]{\int_M \braket{#1}{#2}} 
\newcommand{\intrho}{\int [D\rho]\ }

\newtheorem{Thm}{Theorem}[section]
\newtheorem{Prop}[Thm]{Proposition}
\newtheorem{Lem}[Thm]{Lemma}
\newtheorem{Cor}[Thm]{Corollary}

\theoremstyle{remark}
\newtheorem{Rem}[Thm]{Remark}
\newtheorem*{Ack}{Acknowledgment}

\theoremstyle{definition}

\newtheorem{Def}[Thm]{Definition}

\newtheorem{Ass}{Assumption}
\newtheorem*{Reg*}{Regularization procedure}

\newcommand{\fullref}[1]{\ref{#1} on page~\pageref{#1}}

\newcommand{\braket}[2]{\left\langle{\,{#1}\,,\,{#2}\,}\right\rangle}
\newcommand{\dbraket}[2]{\left\langle\!\left\langle
{\,{#1}\ ;\,{#2}\,}\right\rangle\!\right\rangle}
\newcommand{\Lie}[2]{{\left[{\,{#1}\,,\,{#2}\,}\right]}}

\newcommand{\BV}[2]{\left({\,{#1}\,,\,{#2}\,}\right)}

\newcommand{\vevo}[1]{{\left\langle\;{#1}\;\right\rangle}_0}

\newcommand{\bbR}{{\mathbb{R}}}

\newcommand{\bbZ}{{\mathbb{Z}}}

\newcommand{\de}{\partial}

\newcommand{\nablat}{{\widetilde\nabla}}
\newcommand{\SBV}{{S}}
\newcommand{\Scl}{{S^{\mathrm{cl}}}}

\newcommand{\calS}{\mathcal{S}}

\newcommand{\calE}{\mathcal{E}}

\newcommand{\sfA}{{\mathsf{A}}}
\newcommand{\sfB}{{\mathsf{B}}}
\newcommand{\sfF}{{\mathsf{F}}}
\newcommand{\sfa}{{\mathsf{a}}}
\newcommand{\sfh}{{\mathsf{h}}}

\newcommand{\sfO}{{\mathsf{O}}}
\newcommand{\sfS}{{\mathsf{S}}}
\newcommand{\sfC}{{\mathsf{C}}}
\newcommand{\sfs}{{\mathsf{s}}}

\newcommand{\frg}{{\mathfrak{g}}}

\DeclareMathOperator{\bbr}{\mathbb{R}}

\DeclareMathOperator{\Imb}{Imb}

\newcommand{\dBRST}{{\delta_{BRST}}}
\newcommand{\dBV}{{\delta_{BV}}}
\newcommand{\DBV}{{\Delta_{BV}}}
\newcommand{\sDBV}{\boldsymbol{\Delta}}
\newcommand{\OBV}{{\Omega_{BV}}}
\newcommand{\sOBV}{{\boldsymbol{\Omega}}}
\newcommand{\sfdelta}{{\boldsymbol{\delta}}}
\newcommand{\sfDelta}{{\boldsymbol{\Delta}}}

\newcommand{\ImbfM}{\Imb_f(S^1,M)}

\DeclareMathOperator{\Pfaff}{Pfaff}

\newcommand\qq{\rm}
\newcommand\cmp[1]{{\qq Commun.\ Math.\ Phys.\ \bf #1}}
\newcommand\jmp[1]{{\qq J.\ Math.\ Phys.\ \bf #1}}
\newcommand\pl[1]{{\qq Phys.\ Lett.\ \bf #1}}
\newcommand\np[1]{{\qq Nucl.\ Phys.\ \bf #1}}
\newcommand\mpl[1]{{\qq Mod.\ Phys.\ Lett.\ \bf #1}}

\newcommand\lmp[1]{{\qq Lett.\ Math.\ Phys.\ \bf #1}}

\newcommand\ijmp[1]{{\qq Int.\ J. Mod.\ Phys.\ \bf #1}}
\newcommand\cqg[1]{{\qq Class.\ Quant.\ Grav.\ \bf #1}}
\newcommand\prept[1]{{\qq Phys.\ Rept.\ \bf #1}}

\newcommand\anp[1]{{\qq Ann.\ Phys.\ \bf #1}}
\newcommand\anm[1]{{\qq Ann.\ Math.\ \bf #1}}

\newcommand\jdg[1]{{\qq J.\ Diff.\ Geom.\ \bf #1}}
\newcommand\phs[1]{{\qq Physica Scripta \bf #1}}
\newcommand\topol[1]{{\qq Topology \bf #1}}

\begin{document} 

\title[Higher-dimensional $BF$ 
theories in the Batalin--Vilkovisky formalism]{Higher-dimensional $BF$ 
theories in the Batalin--Vilkovisky formalism: The BV action and generalized 
Wilson loops}
%\date{}

\author[A.~S.~Cattaneo]{Alberto~S.~Cattaneo}
\author[C.~A.~Rossi]{Carlo~A.~Rossi}
\address{Mathematisches Institut --- Universit\"at Z\"urich--Irchel ---  
Winterthurerstrasse 190 --- CH-8057 Z\"urich --- Switzerland}  
\email{asc@math.unizh.ch; crossi@math.ethz.ch}
\thanks{A.~S.~C. acknowledges partial support of SNF Grant
No.~2100-055536.98/1}

\begin{abstract}
This paper analyzes in details the 
Batalin--Vilkovisky
quantization procedure for $BF$ theories on an $n$-dimensional manifold 
and describes a suitable superformalism to deal with the master equation
and the search of observables.
In particular,
generalized Wilson loops for $BF$-theories with additional
polynomial $B$-interactions are discussed in any dimensions.
The paper also contains the explicit proofs to the Theorems
stated in \cite{CR}.
\end{abstract}

\maketitle

\tableofcontents

\section{Introduction}
Topological $BF$ theories \cite{S,H,BlT}
are in three dimensions just another way
of writing Chern--Simons theory \cite{W} (at least at the perturbative level
and disregarding anomalies) \cite{CCFM}. 
In particular, they produce $3$-manifold and
knot invariants. The $2$-dimensional version, at least in its 
{\em canonical}\/ version (see Remark~\ref{rem-cBF}),  is a particular
case of the Poisson sigma model \cite{Ik,SS} and
describes the deformation
quantization of the dual of the corresponding Lie algebra \cite{CF}.
Higher-dimensional generalizations did not have interesting topological
interpretations---apart from the partition functions which describe torsion
invariants---due to the lack of interesting observables.
These were recently introduced in \cite{CR} 
thanks to a superformalism
that simplifies a lot the combinatorics of the associated Batalin--Vilkovisky
(BV) cohomology. The meaning of these observables is that their
expectations values are cohomology classes on the space of framed
imbeddings of $S^1$ (as those described in \cite{CCL}).

The superformalism for $BF$ theories in the BV framework
is not entirely new as it was proposed in
\cite{Wall,Ikem}, but unfortunately the sign rules were not spelt out.
So the first aim of this paper, after reviewing $BF$ 
theories in Section~\ref{sec-BFtheories} and the BV 
formalism in Section~\ref{sec-BV}, is to give a complete description
in Section~\ref{sec-BVsuper} of the superformalism and of its properties,
including explicit sign rules. This leads to a straightforward
proof of the master equation for $BF$ theories in Section~\ref{sec-BVaction},
both for ordinary $BF$ theories and for their ``canonical'' version
(see subsection~\ref{sec-dual}). If the resulting BV action is interpreted
as a supersymmetric sigma model, it falls in the general framework
discussed in \cite{AKSZ}, see subsection~\ref{sec-sm}. However, the
superformalism used in this paper is better suited for
the introduction (see Sections~\ref{sec-gWl} and \ref{sec-lobs})
of generalized Wilson loops, which are constructed
in terms of Chen's iterated integrals \cite{Chen} involving the superfields.
We refer to \cite{CR} for an overview of the above topics and for
a discussion of the main Theorems, which we prove here.
In Section~\ref{sec-gen} we finally discuss the generalization of the above
construction to the case of nontrivial bundles.

Now we briefly comment on the observables defined in this paper
and in \cite{CR}. First,
they are defined on loops (observables related to higher-dimensional
submanifolds are of course of great interest and will be discussed
in a forthcoming paper; we refer to \cite{CM,CCFM,CCR} for previous attempts
in this direction). Second, the quantum BV formalism requires considering
the so-called BV Laplacian (see subsection~\ref{sec-Lapl}) and this forces one
to restrict to imbeddings (more precisely, to framed imbeddings).
Third, one needs restrictions on the Lie algebra underlying the definition
of the $BF$ theory; a semisimple Lie algebra (or more generally a Lie algebra
as in Assumption~\fullref{Ass3}) will do only in the odd-dimensional case
and for a specific observable whose expectation value involves Feynman
diagrams that, apart from an obvious dependence of the propagators on the 
dimension, are exactly the same as in the computation of knot invariants
from Chern--Simons theory (see \cite{CCL} and references therein); 
the main
characteristic of these diagrams is that they involve completely
antisymmetric trivalent vertices which satisfy a diagrammatic version of
the Jacobi identity (see~\cite{BN2}): 
namely, each vertex represents a binary operation
with the same properties of a Lie bracket.
More general observables (and, in particular, any observable in the 
even-dimensional case) requires stronger restrictions on the involved
Lie algebra; in particular, our construction works
when it comes from an associative algebra
(see Assumption~\fullref{trodd}). The Feynman diagrams for the expectation
values of these generalized Wilson loops contain $(k+1)$-valent vertices
(for any $k\ge2$) which can be interpreted as $k$-ary operations satisfying
certain conditions. The way $BF$ theories generate vertices (i.e., $k+1$-ary
operations)
is through one single binary operation (i.e., an associative product)
plus a trace;
but in principle there might exist more general algebraic structures generating
graph cocycles
(in the sense of \cite{K} and made precise in \cite{CCL}) which yield
cohomology classes of imbedded circles.

Observe that, both in the odd and the even dimensional cases,
one can define observables involving only $3$-valent vertices
(in the odd-dimensional case, they are the first observables pointed out
above which do not require the extra assumption on the Lie algebra).
The theory itself formally looks like a $3$-dimensional
$BF$ theory with cosmological term (so essentially a Chern--Simons theory)
``rigidly'' transported to higher dimensions and this might provide
an interpretation of the Chern--Simons theory for strings proposed in
\cite{W-1}. 

\begin{Ack}
We thank P. Cotta-Ramusino for stimulating discussions and
R.~Longoni for carefully reading the manuscript and for useful comments.
\end{Ack}

\section{Preliminaries}
In this Section we introduce the main objects that we will use throughout the paper.
We begin by considering a principal $G$-bundle $P\to M$, where $M$
is a connected orientable manifold of dimension $m\ge2$. 
We will denote by $\mathfrak{g}$ the Lie algebra of $G$.
We consider the associated bundle $P\times_{G} V$ for a $G$-module $V$.
In particular, we will be interested in 
$\ad P:=P\times_{G} \Lg$ and
$\ad^* P:= P\times_{G} \Lg^*$.

We denote by $\Omega^{*}(N;V)$ the space of $V$-valued forms on 
a manifold $N$.  By $\Omega_{\text{bas}}^*(P; V)$ we denote the invariant, horizontal forms on $P$ taking values in $V$; then
\[
\Omega_{\text{bas}}^*(P; V)\cong \Omega^*(M, P\times_{G} V),
\]
where $\Omega^*(M, W)=\Gamma(M,\bigwedge^*T^*M\otimes W)$ for
a vector bundle $W\to M$.
For $P$ trivial, one also has
\[
\Omega^*(M, P\times_{G} V)\cong\Omega^*(M; V).
\]

The gauge group $\mathcal{G}$ of $P$ is defined as the set of all equivariant automorphisms of $P$; it can be identified with the set $\Gamma(M, \Ad P)$ of all sections of the bundle $\Ad P:=P\times_{G} G$, where $G$ acts on itself
by conjugation.
For $P$ trivial, it can be identified with the group $\Gamma(M,G)$ of maps 
from $M$ to $G$.

Another important ingredient that we need is the Universal Enveloping Algebra
(UEA) of a Lie algebra.
We denote by $\Ug$ the UEA of $\Lg$. We recall that $\Ug$ is an associative algebra.
Further, we denote by $\iota:\Lg \to \Ug$ the canonical inclusion of $\Lg$ into $\Ug$ (which is a Lie algebra morphism).

Throughout the paper we will be confronted with the problem of integrating along fibers forms with values in some vector space $V$; for the main properties of the push-forward of real forms and for its generalizations to the case of forms with values in some algebra, we refer to Appendix~\ref{app-A}.

We end this section with some simplifying assumptions that we consider
throughout the paper unless explicitly stated otherwise.
\begin{Ass}\label{Ass1}
The manifold $M$ is compact and there is a flat connection $A_0$ on $P$, such that all the cohomology groups $H^*_{\dd_{A_0}}(M,\ad P)$ are trivial.
\end{Ass}

\begin{Ass}\label{Ass2}
The principal bundle $P$ is trivial.
\end{Ass}

\begin{Ass}\label{Ass3}
The Lie algebra $\frg$ is endowed with a symmetric, $\Ad$-invariant, nondegenerate bilinear form $\braket{}{}$ (e.g., if $\frg$ is semisimple, we may take the Killing form).
In the following, we will extend this form to $\Omega^*(M,\ad P)$ in the usual way.
\end{Ass}

\begin{subsection}{A brief discussion of Assumption \ref{Ass1}}
Assumption \ref{Ass1} is very strong; we want to briefly discuss it in view of future applications (definitions of loop observables).
Let us suppose for a moment that we consider a Lie group $G$, whose Lie algebra $\Lg$ satisfies the third assumption, and a compact, oriented manifold $M$ of dimension $m$. In particular, the $0$-th and the $m$-th De Rham cohomology groups are nontrivial. Let us suppose additionally that the Lie algebra $\Lg$ possesses some invariant elements under the adjoint action of $G$ (i.e., the $0$-th cohomology group $\mathcal{H}^0(G,\Lg)$ is nontrivial. Then it can be shown that the first of the above assumptions cannot hold true. More generally, if the $0$-th cohomology group $\mathcal{H}^0(G,\Lg)$ is nontrivial, and the manifold $M$ is compact and oriented, then there exists no flat connection $A_0$ on $P$ such that $H^*_{\dd_{A_0}}\left(M,\ad P\right)$ is trivial. This implies that, e.g., Assumption \ref{Ass1} is {\em not} compatible with the case of a compact, oriented manifold $M$ and the Lie algebra $\Lg=\mathfrak{gl}(N)$.
However, we may assume that Assumption 1 holds true for odd-dimensional compact, oriented manifolds; this is in analogy with the assumption made by Axelrod and Singer in \cite{AS} in the introduction of \cite{AS}. This forces us to exclude principal bundles $P$, whose structure group $G$ possesses nontrivial $0$-th cohomology group with coefficients in $\Lg$. For the even dimensional case, we may consider special even-dimensional manifolds arising as products of two odd-dimensional manifolds $M_1$ and $M_2$, one of which (say $M_1$) is the base space of a principal bundle $P$ with Lie group $G$, satisfying Assumption 1, with flat acyclic connection $A_0$. We consider then on $M_1$ the complex $\left(\Lambda^* T^*M_1\otimes \ad P,\dd_{A_0}\right)$ and on $M_2$ the complex $\left(\Lambda^* T^*M_2,\dd\right)$; both are elliptic complexes, and the first one is acyclic by assumption. We take then the exterior tensor product of the two complexes defined on $M_1\times M_2$, with indu!
ced differential; this is again an elliptic complex, and, by the Kuenneth Theorem and the acyclicity of the complex on $M_1$, it is acyclic.
So, the existence of odd-dimensional manifolds, for which Assumptions 1 holds true implies the existence of even-dimensional manifolds, for which Assumption 1 is valid. So, we have found some algebraic-topological obstructions to the existence of odd-dimensional compact, oriented manifolds for which Assumption 1 is valid, but we are still not able to produce a definitive criterion for the existence of such manifolds.
We work therefore under the hope that there are Lie groups $G$ and odd-dimensional compact, oriented manifolds, for which Assumption 1 holds. In the case $G=GL(N)$, as we have seen before, there are no such manifolds. So, in this case, i.e.\ in section 8, we choose implicitly $M=\bbr^n$ with the flat trivial connection.
\end{subsection}

\section{$BF$ theories}\label{sec-BFtheories}
%\subsection{The $BF$ action functional}
The fundamental ingredients that we need are a connection $1$-form $A$
on $P$ and an $(m-2)$-form of the adjoint type $B$.
We then construct the curvature $F_A$ of the connection and define the classical
$BF$ action as the functional
\begin{equation} 
\Scl=\int_{M} \braket{B}{F_{A}} \label{BFaction}.
\end{equation}
\begin{Rem}
A more natural setting would be to consider $B$ as a form of the coadjoint 
type. In this case, one would not have to introduce an invariant bilinear
form on $\Lg$, and Assumption~\ref{Ass3} could be discarded. Instead one would
use the canonical pairing between $\Lg^*$ and $\Lg$. We will call these
theories {\em canonical $BF$ theories}\/ and 
will comment
more on them in subsection~\ref{sec-dual}.
However, for the main purposes of this paper (namely, to define loop
observables), one needs anyway to consider $B$ of the adjoint type
(or to introduce an isomorphism between $\Lg$ and its dual). So we will
stick for most of the paper
to the setting described in this section.
\label{rem-cBF}
\end{Rem}
Let us first compute the Euler--Lagrange
equations of motion for the $BF$ action; they are given by the couple of equations
\begin{equation} \label{eqmotion}
F_{A}=0,\quad \covar{A}{B}=0.
\end{equation}
In the following, by ``on shell'' we will refer to the space of solutions with
the extra condition that the connection $1$-form is as in 
Assumption~\ref{Ass1}.
Next we turn to the symmetries of this action:
\[
A\mapsto A^{g},\qquad B\mapsto \Ad(g^{-1})B +\dd_{A^g} \tau_1,
\]
where by $A^g$ we denote the right action of the gauge group element $g$
on the connection $A$, and $\tau_1$ is an element of $\Omega^{m-3}(M,\ad P)$.
The symmetries under which the $BF$ action is invariant can be interpreted as the action of the semidirect product $\mathcal{G}\rtimes_{\Ad}\Omega^{m-3}(M,\ad P)$ on $\mathcal{A}\times\Omega^{m-2}(M,\ad P)$, where $\mathcal{A}$ denotes
the space of connections on $P$.
In infinitesimal form we obtain
\begin{equation} \label{gauge}
A\mapsto A +\epsilon\, \dd_{A} c,
\qquad B\mapsto B+\epsilon([B,c]+\dd_{A} \tau_1),
\end{equation}
where $c$ is in $\Omega^0(M,\ad P)$ (the Lie algebra of $\mathcal{G}$).

These symmetries are reducible on shell, i.e.\ each solution $(A_0;B_0)$ 
with $A_0$ as in Assumption~\ref{Ass1} has as isotropy group the semidirect product $\{e\}\rtimes_{\Ad} \{\tau_1\in \Omega^{m-3}(M,\ad P):\covar{A_0}{\tau_1}=0\}$.
This isotropy group is isomorphic to $\Omega^{m-4}(M,\ad P)/\covar{A_0}{\Omega^{m-5}(M,\ad P)}$, because of Assumption $1$; there are in this quotient nontrivial isotropy groups isomorphic to $\Omega^{m-5}(M,\ad P)/\covar{A_0}{\Omega^{m-6}(M,\ad P)}$, and so on until we arrive at $\Omega^0(M,\ad P)$ which acts freely on $\Omega^1(M,\ad P)$.
Therefore, we have to adopt the extended BRST procedure to consistently fix all the symmetries, by introducing a hierarchy of ghosts for ghosts.
Unfortunately the isotropy groups off shell are different from the above groups; so we have to resort to the BV formalism which generalizes BRST and works
also in this case;
see the next subsections for more details on both procedure.

\subsection{The BRST procedure}
For the sake of simplicity,
let us restrict ourselves for the moment to the special case $m=4$.

We first promote the $0$-form $c$ and the $1$-form $\tau_1$
appearing in the infinitesimal gauge transformations \eqref{gauge}
to anticommuting fields of ghost number $1$;
$A$ (and every variation of $A$ which is a $1$-form) and $B$ will be given ghost number $0$.
We then define the BRST operator $\dBRST$ for the $4$-dimensional $BF$ theory by the rules
\begin{equation}
\dBRST A = \covar{A}{c},\quad \dBRST B = [B,c]+\covar{A}{\tau_1},
\quad \dBRST c=-\frac{1}{2}[c,c],
\label{BRSTABc}
\end{equation}
and
\[
\dBRST \tau_1=-[\tau_1,c]+\covar{A}{\tau_2},\quad \dBRST \tau_2 = [\tau_2,c],
\]
where $\tau_2$ is a form in $\Omega^0(M,\ad P)$ to which we assign ghost 
number two. Then
$\dBRST$ is an odd operator of ghost number $1$ and a differential for the Lie bracket.
By the graded Leibnitz w.r.t.\ the ghost number, it follows that
\[
\delta_{\mathrm{BRST}}^2 B = [F_{A},\tau_2] \neq 0,
\]
while for the other fields, $\delta_{\mathrm{BRST}}^2=0$. We notice that a sufficient condition for $\dBRST$ to be a differential is $F_{{A}}=0$; this is exactly the first equation in (\ref{eqmotion}). 
Otherwise the BRST quantization procedure fails, but the BRST operator closes on shell; we can therefore apply to this situation another formalism to quantize the $BF$ theory, namely the BV quantization procedure which works well for such a theory.
A similar problem arises for any $m\ge4$.

In general, however, 
because of the on-shell reducibility discussed in the last subsection, we have to introduce more ghosts for ghosts $\tau_k$ with values in $\Omega^{m-2-k}(M,\ad P)$, $k=1,\dots,m-2$, and ghost number $k$. The BRST operator 
is defined by \eqref{BRSTABc} and by the rules:
\begin{equation}
\begin{aligned}
\dBRST \tau_k&=(-1)^k\Lie{\tau_k}{c}+\covar{A}{\tau_{k+1}},\\
\dBRST \tau_{m-2}&=(-1)^m \Lie{\tau_{m-2}}{c}.
\end{aligned}
\label{BRSTtau}
\end{equation}
It is then easy to see that $\delta_{\mathrm{BRST}}^2=0$ mod $F_{A}$.

The case $m=2$  and $m=3$ are the only ones in which the BRST formalism works, but one can apply the BV formalism there as well obtaining equivalent results.

\subsection{Classical observables}\label{class-oss}
We start by considering
$\tr_{\rho}\holo{A}01$,
where by $\holo{A}01$ we denote the inverse of the usual holonomy w.r.t.\ the connection $A$ viewed as a $G$-valued function on $\lo M$
(see Remark~\fullref{rem-hol} for technical details).
By considering a representation $(\rho,V)$ we then obtain an $\Aut V$-valued
function, which under the trace then yields an ordinary function.
This function depends also on the choice of a connection $A$, but its very definition implies that it is invariant w.r.t.\ the action of the gauge group $\mathcal{G}$ on the space $\mathcal{A}$ of connections on $P$, so it defines a function on $\mathcal{A}/\mathcal{G}\times \lo M$
($\mathcal{A}/\mathcal{G}$ is the moduli space of $G$-connections).
We notice that in a local trivialization, the inverse of the holonomy possesses a representation in terms of a formal series of iterated integrals.
In the case $P$ trivial, the holonomy becomes a function on $\lo M$ with values in $G$.

Next we define
\begin{equation} \label{Taylorholon}
h_{n,\rho}(A,B):=\tr_{\rho}\left\{\pi_{n*}\left[\widehat{B}_{1,n}\wedge\dots\wedge\widehat{B}_{n,n}\right]\holo{A}01\right\}\ ,
\end{equation}
where the notations are as in Appendix~\ref{app-hol}.

{}From now on, we will omit the wedge product between forms.
Notice that we have already omitted to write $\rho$\  before all forms in the definition of $\widehat{B}_i$; for all $i$, $\widehat{B}_{i,n}$ is a form on $\lo M\times \triangle_n$ with values in $\en{V}$.
It follows from the definition that $h_{n,\rho}(A,B)$ for all $n$ is a differential form of degree $(m-3)n$ on $\lo M$.
\begin{Prop}
Let $A$ and $B$ be on shell; then $h_{n,\rho}(A,B)$ is a closed form for $m$ odd and for all $n$, while for $m$ even and greater than $4$, the $h_{2k+1,\rho}(A,B)$'s are closed.
\end{Prop}
\begin{proof}[Proof (sketch)]
Since $F_A=0$ and $\dd_{A} B=0$, as a consequence of  Theorem~\ref{holonomy1}
the following identities hold:
\begin{align*}
&\dd_{\pi_1^* \ev(0)^* A}\widehat{B}=\widehat{\dd_A B}=0,\ \forall i;\\
&\dd_{\ev(0)^* A}\holo{A}{0}{1}=0
\end{align*}
as a consequence of~\ref{derpar}.
The cyclicity of $\tr_{\rho}$ implies $\tr_{\rho} \dd_{\ev(0)^* A}=\dd \tr_{\rho}$.

We recall now that in the definition \eqref{Taylorholon}
before the products of the $\widehat{B}_{i,n}$'s there is a push-forward; 
thus, in order to compute $\dd h_{n,\rho}$,
we will have to apply the generalized Stokes Theorem \eqref{prop-push}.
The boundary of the $n$-simplex can be written as the union of other $(n-1)$-simplices (the faces of the simplex), corresponding to the collapsing of successive points, plus two other faces, where the first point tends to $0$ or the last tends to $1$.
The faces of the first type give $0$, because they yield terms containing $\widehat{B^2}_{i,n}$, which vanish for dimensional reasons. The remaining two faces give
\begin{align*}
&-(-1)^{m(n-1)}\tr_{\rho}\left\{\ev(0)^*B\pi_{n-1*}\left[\widehat{B}_{1,n-1}\dots\widehat{B}_{n-1,n-1}\right]\holo{A}{0}{1}\right\}\\
&+(-1)^{n+1}\tr_{\rho}\left\{\pi_{n-1*}\left[\widehat{B}_{1,n-1}\dots\widehat{B}_{n-1,n-1}\right]\holo{A}01 \ev(0)^*B\ \right\} ;
\end{align*}
again, for $m=\dim M$ odd, the cyclicity of $\tr_{\rho}$ implies that these terms cancel each other. This also works for $m$ even, in case $n$ is odd.
On the other hand, when both $m$ and $n$ are even, these two terms have the same sign, and therefore they do not cancel each other. 
\end{proof}
Similar computations show that the $h_{n,\rho}(A,B)$'s are observables on shell and modulo exact terms, either if $m$ is odd and greater than $5$ or if $m$ is even and greater than $4$ but $n$ is odd.

Another advantage of the BV formalism that we must introduce for the reasons
explained before is that it allows to deal with observables which are
BRST closed on shell, upon extending them suitably. This will be explained
in the next sections.

\section{The Batalin--Vilkovisky quantization procedure for $BF$ theories}
\label{sec-BV}
We now briefly review the BV formalism \cite{BV}, though in a form already
adapted to $BF$ theories. For a general account on the 
formalism, see 
e.g.\ \cite{Sta} and references therein.

Let us consider all the fields of the theory, i.e.\ the connection one-form $A$
(which we write $A_0+a$, where $A_0$ is a given flat connection on $P$, and
$a$ is an element of $\Omega^1(M,\ad P)$), 
the tensorial $(m-2)$-form B of adjoint type,
the ghost $c$ with values in $\Omega^0(M,\ad P)$ and the
ghosts $\tau_j$, $j \in \{1,\ldots,m-2\}$, for which holds: $\tau_j$ takes 
values in $\Omega^{j}(M, \ad P)$ and has ghost number $\gh \tau_j=j$. 

We then associate to each field $\phi^{\alpha}$ a canonical ``antifield,''
denoted by $\phi^{+}_{\alpha}$, as follows: suppose that the field
$\phi^{\alpha}$ has degree $\dg{\phi^{\alpha}}$ and ghost number
$\gh{\phi^{\alpha}}$; then the antifield $\phi^{+}_\alpha$ is a form on
$M$ with values in $\ad P$, whose degree is set to be equal to $m-\dg{\phi^{\alpha}}$ 
and its ghost number is set to be $-1-\gh{\phi^{\alpha}}$. 

The fundamental ingredients of the BV antibracket are the left and right partial derivatives of a functional $F$, which we are going to define precisely in the following subsection.
 
To simplify the notations from now on we will denote all the fields 
and antifields collectively as ``fields'' and will use the symbols 
$\varphi_{\alpha}$, where $\alpha$ runs from $1$ to $(2m+2)$;
$\mathcal{M}:=\{\varphi_{\alpha}\}_\alpha$. 

\subsection{Functional derivatives}\label{sec-funderiv}
We pick a commutative
algebra $\mathfrak{A}$ (usually, we take $\mathfrak{A}=\mathbb{R}$ or $\mathfrak{A}=\mathbb{C}$, but see Remark~\ref{A-NE}).
We are going to consider (formal) power series of
local functionals in the fields
taking values in $\mathfrak{A}$.
We introduce a grading, which on monomials is defined as the sum of
the ghost numbers of all the fields appearing there and which is then
extended by linearity.
We finally consider the graded commutative algebra $\mathcal{S}(\mathfrak{A})$
generated by such objects.
We then define the left and right functional derivatives of an element $F$ in
$\mathcal{S}(\mathfrak{A})$ w.r.t.\ the field $\varphi_{\alpha}$ by
\begin{displaymath} 
\ddt F(\varphi_{\alpha}+t\rho_{\alpha}) 
=\int_{M} \braket{\rho_{\alpha}}{\frac{\ld F}{\der}}= \int_{M}
\braket{\frac{F \rd}{\der}}{\rho_{\alpha}}.
\end{displaymath} 
It follows from these definitions that the functional derivatives 
are in general
distributional forms. For convenience of notations, however, we will denote
the space of distributional forms with the same symbol $\Omega^*$ used
for smooth forms since this causes no harm. 
So the functional derivatives of $F$
in $\mathcal{S}(\mathfrak{A})$ w.r.t.\ $\varphi_\alpha$
are elements of
$\Omega^{p_\alpha}(M,\ad P\otimes \mathfrak{A})$, with the property that 
\begin{equation} 
p_\alpha := \dg \frac{\ld F}{\der} = \dg \frac{F \rd}{\der} =  
m-\dg \varphi_{\alpha}. 
\label{palpha}
\end{equation} 
As for the ghost number one has 
\begin{equation} 
g_\alpha :=\gh \frac{\ld F}{\der} = \gh \frac{F \rd}{\der} = 
\gh F - \gh \varphi_{\alpha}. 
\label{galpha}
\end{equation} 
{}From the definitions and the above introduced notations, 
we also obtain the following useful identities: 
\begin{equation} \label{relder} 
\frac{\ld F}{\der}=(-1)^{p_{\alpha}\dg{\phi_{\alpha}}+
g_{\alpha}\gh{\phi_{\alpha}}}\ \frac{F \rd}{\der}.
\end{equation}  

Beside the manifold $M$, let us consider another (possibly
infinite dimensional) manifold $N$ (e.g., $\lo M$).            
In the following we will also consider (formal) power series of
functionals taking values in $\Omega^*(N;E)$, for some associative
algebra $E$ (e.g,
$\mathbb{R}$, $\mathbb{C}$, $\Ug$ or $\en{V}$, for some $\Lg$-module
$V$).
On this space we can introduce two gradings: the first is the ghost number
which is defined as in the case of $\mathcal{S}(\mathfrak{A})$;
the second is simply the form degree on $N$.
We denote by
$\mathcal{S}^*(N;E)$ the bigraded
superalgebra generated by such functionals
(this superalgebra is supercommutative if{f} $E$ is).
Let us notice, at last, that for $E$ an $\mathfrak{A}$-module, $\mathcal{S}(\mathfrak{A})$ is a subalgebra of $\mathcal{S}^*(N;E)$.

For the left (resp.\ right) derivative of a functional in
$\mathcal{S}^*(N;E)$, we use the canonical identification of
$\Omega^{p}(M, \ad P\otimes\Omega^{q}(N;E))$ with 
$\Omega^{p,q}(M\times N, \ad P\boxtimes E)$ (respectively with
$\Omega^{q,p}(N\times M, E\boxtimes \ad P)$). 
We next introduce the following notations:
\begin{align*} 
\pi_1:N \times M &\longrightarrow N, &
 (\tilde{x}, x) &\longmapsto \tilde{x}, \\ 
\pi_2:N \times M &\longrightarrow M, &
(\tilde{x}, x) &\longmapsto x, \\ 
\intertext{and}
\tilde{\pi}_1:M \times N &\longrightarrow M, &
(x; \tilde{x}) &\longmapsto x, \\
\tilde{\pi}_2:M \times N &\longrightarrow N, &
(x; \tilde{x}) &\longmapsto \tilde{x}. 
\end{align*} 
\label{A-NE} 
We have used the following useful notation: Let $\mathcal{E}\to M$ and $\mathcal{F}\to N$ vector bundles over $M$, resp.\ $N$. Then we define
\begin{align*}
\mathcal{E}\boxtimes \mathcal{F}&:=\tilde{\pi}_1^*(\mathcal{E})\otimes \tilde{\pi}_2^*(\mathcal{F}),\quad \text{resp.}\\
\mathcal{F}\boxtimes\mathcal{E}&:=\pi_1^*(\mathcal{F})\otimes \pi_2^*(\mathcal{E});
\end{align*}
it follows that they are vector bundles over $M\times N$, resp.\ $N\times M$. 

With these notations we can finally define the functional
derivatives of $F$ in $\mathcal{S}^*(N;E)$:
\begin{align*} 
\ddt F(\varphi_{\alpha}+t\rho_{\alpha})&= \tilde{\pi}_{2 *}
\braket{\tilde{\pi}_1^{*}\rho_{\alpha}}{\frac{\ld
F}{\der}}=\pi_{1*}\braket{\frac{F \rd }{\der}}{\pi_2^{*}\rho_{\alpha}}. 
\end{align*} 
The functional derivatives have now two different form degrees:
one is the form degree w.r.t.\ $M$ and is still given by \eqref{palpha};
the other is the form degree w.r.t.\ $N$ and remains
equal to $\dg F$. 
The ghost number is given by \eqref{galpha} as before.

\subsection{The BV antibracket}\label{subs-BVanti}
We define the BV antibracket for two elements $F$, $G$ in
$\mathcal{S}(\mathfrak{A})$ as the functional:
\begin{displaymath} 
\BV{F}{G}:=\int_{M} \braket{\frac{F \rd}{\da}}{\frac{\ld G}{\dad}}
- (-1)^{(m+1) \dg{\phi^{\alpha}}} \braket{\frac{F\rd}{\dad}}{\frac{\ld G}{\da}}. 
\end{displaymath} 
We note that this functional is again in $\mathcal{S}(\mathfrak{A})$, since 
we integrate over $M$ and since the functional derivatives of an
element of $\mathcal{S}(\mathfrak{A})$ are once again power series; 
it is not difficult to see that the ghost
number of the BV antibracket of two homogeneous elements $F$ and $G$ in
$\mathcal{S}(\mathfrak{A})$, with ghost numbers $\gh F$ and $\gh G$, 
is given by
$\gh F+\gh G+1$.
Next, we define the BV antibracket for two functionals $F$ and $G$ in
$\mathcal{S}^*(N;E)$ by the formula:
\begin{equation}\label{BVE}
\BV{F}{G}:=\pi_{13*}\braket{\pi_{12}^* \rfuna{F}}{\pi_{23}^*
\lfunad{G}} - (-1)^{\dg{\phi^\alpha} (m+1)} \pi_{13*}\braket{\pi_{12}^* \rfunad{F}}{\pi_{23}^* \lfuna{G}},
\end{equation}
where we use the projections
\begin{align*}
&\pi_{12}:N\times M\times N\to N\times M,\quad (n_1; m; n_2)\mapsto (n_1; m);\\
&\pi_{23}:N\times M\times N\to M\times N,\quad (n_1; m; n_2)\mapsto (m; n_2);\\
&\pi_{13}:N\times M\times N\to N\times N,\quad (n_1; m; n_2)\mapsto (n_1; n_2).
\end{align*}
This formula needs some explanations. Let us suppose that $F$ and $G$ are
homogeneous as elements of $\Omega^{*}(N;E)$, with degrees $\dg F$,
resp.\ $\dg G$. 
We consider the case of a trivial algebra bundle $\mathcal{E}=N\times E$
over $N$; in this case, the left functional derivatives are elements of
$\Omega^{\dg F, p}(N\times M, E\boxtimes \ad P)$, while the right ones
are elements of $\Omega^{q, \dg F}(M\times N, \ad P\boxtimes E)$. 
The product that we write in this case denotes two operations: the first
consists in the usual wedge multiplication of the form parts, while the second
is the multiplication in $E$ of the algebra part. (We refer to the beginning of
Appendix~\ref{app-sign} for further details.) 
Therefore, in this special
case, the BV antibracket of two homogeneous functionals $F$, $G$, in
$\mathcal{S}^*(N;E)$ gives as a result a homogeneous element of
$\mathcal{S}^*(N;E)$, with degree in $N$ equal to $\dg F+\dg G$ and
ghost number $\gh F+\gh G +1$.

We last define the BV antibracket for two special functionals, for we 
will often consider this case in the following: namely, we pick a functional $F$
in $\mathcal{S}(\mathfrak{A})$ and a
functional $G$ in $\mathcal{S}^*(N;E)$, where $E$ is an $\mathfrak{A}$-module:
\begin{equation*}
\BV{F}{G}:=\tilde{\pi}_{2 *} \braket{\tilde{\pi}_1^{*} \rfuna{F}}{\lfunad{G}} - (-1)^{\dg{\phi^{\alpha}} (m+1)} \tilde{\pi}_{2*} \braket{\tilde{\pi}_1^{*} \rfunad{F}}{\lfuna{G}}.
\end{equation*}
It is clear that in this case the BV antibracket of $F$ and $G$ is an element of
$\mathcal{S}^*(N;E)$. For homogeneous elements,
the degree of the antibracket is equal to the degree of $G$,
while $\gh \BV{F}{G}=\gh F+\gh G+1$.

\subsection{Properties of the BV antibracket}\label{sec-propBV}
We recall first, in a unified way, the ghost and degree properties
of the antibracket. We denote by $\mathcal{S}$ the algebra of functionals
(which according to the case may be $\mathcal{S}(\mathfrak{A})$ or
$\mathcal{S}^*(N;E)$) and by $\mathcal{S}^{p,g}$
the subspace of homogeneous functionals
of form degree $p$ and ghost number $g$ by  (in the case
of $\mathcal{S}(\mathfrak{A})$, $p$ is necessarily zero).
Then the antibracket is a bilinear operator
\[
\BV{\ }{\ }:\mathcal{S}^{p,g}\otimes\mathcal{S}^{p',g'}\to
\mathcal{S}^{p+p',g+g'+1}.
\]
We list (without proofs) some useful identities for the
BV antibracket.
We begin with the graded commutativity
\[
\BV{F}{G}=-(-1)^{\dg F \dg G +(\gh F +1)(\gh G+1)}(G,F),
\]
which holds whenever one of the two functionals is central.
The next property is the graded Jacobi identity
\[
(-1)^{\dg F\dg H+(\gh{F}+1)(\gh{H}+1)}
\BV{F}{\BV{G}{H}} + \text{cyclic permutations} = 0.
\]
which holds whenever two of the three functionals are central.
The last property is the graded Leibnitz rule
\[
\BV{F}{GH}=\BV{F}{G}H + (-1)^{\dg F \dg G + (\gh F +1)\gh G}G\BV{F}{H},
\]
which holds whenever $F$ or $G$ or both are central. In particular this holds
in the following important cases: $i$) all functionals are in
$\mathcal{S}(\mathfrak{A})$; $ii$) all functionals are in $\mathcal{S}^*(N;E)$
with $E$ a commutative algebra; $iii$) $F$ or $G$ or both are in 
$\mathcal{S}(\mathfrak{A})$ and the remaining functional(s) are in
$\mathcal{S}^*(N;E)$ for $E$ an $\mathfrak{A}$-module.
\begin{Rem}
If we restrict ourselves to $\mathcal{S}(\mathfrak{A})$, then the
above properties hold on the whole algebra. An algebra with a bracket
satisfying the above properties is known as a Gerstenhaber algebra
\cite{Gerst}.
\end{Rem}

The Leibnitz rule will play a key-r\^ole in the following section, where we 
define the BV operator via the BV antibracket; the
functional $F$ will be there the BV action for the $BF$ theory. 
Let us in fact suppose that we have a homogeneous local functional $S$ in
$\mathcal{S}(\mathfrak{A})$ with even ghost
number (usually, $\mathfrak{A}=\mathbb{R}$ and $\gh S=0$).
We can then define the following operator on the
superalgebra $\mathcal{S}^*(N;E)$, with $E$ an $\mathfrak{A}$-module: 
\begin{displaymath} 
\delta_{S} F := \BV{S}{F}. 
\end{displaymath} 
It follows easily from the   $\mathfrak{A}$-linearity of the
BV antibracket that $\delta_{S}$ is
a $\mathfrak{A}$-linear operator on the
algebra $\mathcal{S}^*(N;E)$. 
The most important property of such an operator is an immediate consequence 
of the Leibnitz rule written above; namely,
\begin{displaymath}  
\delta_{S} (FG)=(\delta_{S} F) G + (-1)^{\gh F} F (\delta_{S} G);  
\end{displaymath} 
i.e.\ the operator $\delta_{S}$ is a graded $(0,\gh S+1)$-derivation on
$\mathcal{S}$. (If moreover  $\BV SS=0$, then 
the Jacobi identity implies $\delta_S^2=0$.)
We now list some other useful properties of the derivation $\delta_S$.
\begin{Lem} \label{pull1} 
Suppose that the functional $F$ lies in $\mathcal{S}^*(N;E)$,
and that we have a map $h:H\to N$ from some manifold $H$ to the
manifold $N$, then the following identity holds:
\begin{equation*} 
\delta_{S} [h^{*}(F)] = h^{*}(\delta_{S} F). 
\end{equation*} 
\end{Lem} 
 
\begin{Lem} \label{push1} 
Suppose that we have a functional $F$ in $\mathcal{S}^*(H;E)$,
where $H$ is the total space of a fiber bundle over $N$
with typical fiber some manifold $B$ (possibly with boundaries or corners)
and projection $\pi$.
The integration along the fiber of the functional $F$ 
yields a functional in $\mathcal{S}^*(N;E)$
with degree $\dg{\pi_{*}(F)}=\dg{F} -
\dim B$, if we suppose additionally that $F$ is homogeneous in the
degree.   
Then we obtain the following identity: 
\begin{equation*} 
\delta_{S} [\pi_{*}(F)] = \pi_{*} (\delta_{S} F) 
\end{equation*} 
\end{Lem} 
 
\begin{Lem} \label{extd1} 
Let us suppose that we have a functional $F$ in $\mathcal{S}^*(N;E)$, 
for some manifold $N$ and some algebra $E$. Let us denote by $\dd$
the exterior derivative on $N$. 
Then the following identity holds: 
\begin{equation*} 
\delta_{S}(\dd F)=\dd(\delta_{S} F).  
\end{equation*} 
\end{Lem} 
 
\begin{Lem} \label{trace1} 
Let us suppose that the functional $F$ belong to the superalgebra
$\mathcal{S}^*(N; \mathfrak{g})$; let us suppose additionally that we
have a $\mathfrak{g}$-module $(V,\rho)$. 
The application of $\tr_{\rho}$ to $F$ gives an element of 
$\mathcal{S}^*(N; \mathbb{R})$.  
Then we obtain the following identity: 
\begin{equation*} 
\delta_{S}[\tr_{\rho}(F)] = \tr_{\rho}(\delta_{S} F). 
\end{equation*} 
\end{Lem} 
We will only sketch a few ideas of the proofs of the above Lemmata. 
 
For Lemma~\ref{pull1} we only have to write down explicitly the expressions for
the two BV antibrackets, which in this special case involve the push-forward
w.r.t\ the projection $\Bar{\pi}_1: H\times M \rightarrow H$, resp. $\pi_1:N\times M
\rightarrow M$, and the pullbacks w.r.t.\ the maps $\Bar{\pi}_2:M\times H
\rightarrow H$, resp. $\pi_2:M\times N \rightarrow M$; these maps do appear in 
the definition of the partial functional derivatives of $F$. 
Then we have to consider the following commutative diagram:
\[  
\begin{CD} 
M\times H  @>\text{id}\times h>> M\times N\\
@V\Bar{\pi}_{2}VV              @VV\pi_{2}V\\
H          @>h>>                 N
\end{CD}
\] 
It is easy to see that $\text{id}\times h$
induces an orientation preserving map (namely, the identity map)
 between the fibers $(\Bar{\pi}_2)^{-1}(\{e\})$
($\cong\{e\}\times M$) and $(\pi_2)^{-1}(\pi(e))$ ($\cong\{\pi(e)\}\times M$), for $e\in H$. From Lemma~\ref{commpull} the claim follows. 

For Lemma~\ref{push1}, we have to write down again explicitly the
BV antibrackets on the two sides of the identity.
In this case we use the following commutative diagram:
\[
\begin{CD}
M\times H  @>\text{id}\times \pi>> M\times N\\
@V\tilde{\pi}_{2}VV              @VV\pi_{2}V\\
H          @>\pi>>                 N
\end{CD} 
\]
The commutativity of this diagram implies that the composite bundles $\mathcal{H}_{B\times M}=(M\times H; \pi \circ \tilde{\pi}_2; N; B\times M)$ and $\mathcal{H}_{M\times B}=(M\times H ; \pi_2 \circ (\text{id} \times \pi); N; M\times B)$ possess the same total space and the same base space, but have different fibers; in fact, the fiber of the first is isomorphic to $B\times M$, while the fiber of the second one is $M\times B$. We can go from a bundle to the other via a bundle morphism which is the identity on the total and on the base space, but which reverses the orientation of the fibers, and we know that the orientation of a fiber bundle is induced by the orientations of the base space and of the fiber; this will imply the following identity: 
\begin{displaymath} 
\pi_{*} \circ{\tilde{\pi}_{2*}}= (-1)^{m\dim B}  \ \pi_{2*}
\circ (\text{id} \times \pi)_{*} , 
\end{displaymath} 
and the coefficient $(-1)^{m\dim B}$ comes from the
orientation reversal of the fibers of the two bundles (for the property of the
push-forward, see Lemma~\ref{pushcom}). This identity will imply the claim.
               
For Lemma~\ref{extd1} we simply apply the generalized Stokes' theorem for the
push-forward w.r.t.\ $\tilde\pi_2:M\times N \rightarrow N$; notice that in 
this case the fiber, i.e.\ $M$, has no boundary.
Then we have to
remember that the exterior derivative $\dd_{M\times N}$ on $M\times N$ splits as
$\dd_{N} +  
\sigma \dd_{M}$, where the sign $\sigma$ is given by  
$\sigma=(-1)^{\dg_N(\omega)}$, for a form $\omega$ on $M\times N$ with  
degree over $N$ equal to $\dg_{N}(\omega)$. 
We have to remember that, in the defining equation for the right functional
derivative, the test form is independent of $N$, therefore the exterior
derivative on $N$ applied to (the pullback w.r.t.\ $\tilde\pi_1$ of) the test form
gives $0$ as result; next, we know that the integrand form has maximal degree
w.r.t.\ $M$, so that the exterior derivative w.r.t.\ $M$ of the integrand gives
$0$. Then the result follows from all the above considerations. 
 
Lemma~\ref{trace1} follows easily from the definition of the partial functional
derivatives and from the fact the trace $\tr_{\rho}$ acts only on the $\End V$-part
of the tensor product (remember that the functionals we are considering take
their values in $\en{V}$ for some $\mathfrak{g}$-module $V$).

\subsection{The BV Laplacian}\label{sec-Lapl}
Let us temporarily choose a Riemannian metric on $M$ and let us
denote by $\star$ the induced Hodge star operator. 
Let us pick a field $\phi^{\alpha}$; we denote by $\phi_\alpha^*$ the field (called sometimes the {\em Hodge dual antifield}\/ of $\phi^{\alpha}$) defined by the formula 
\[
\phi_\alpha^*:=\star \phi_\alpha^{+},
\]
where $\phi_\alpha^+$ is the antifield of $\phi^\alpha$.
It follows easily from the definition that the degree of $\phi_\alpha^*$ is given by the degree of the field it is associated to, while its ghost number is given by $-1-\gh \phi^{\alpha}$.
Let $\alpha$, $\beta$ be two elements in $\Omega^p(M, \ad P)$. We define
\[
\sqbr{\alpha}{\beta}:=\int_M \braket{\alpha}{\star \beta}.
\]

Now we define a new type of functional derivatives. We begin with functionals in the space $\mathcal{S}(\mathfrak{A})$. Let us once again denote collectively
by $\varphi_\alpha$
the fields $\phi^\alpha$ and their newly defined antifields $\phi_\alpha^*$.

Let $\rho_\alpha$ be a form with the same degree and ghost number as 
$\varphi_\alpha$. 
Let $F$ be an element in $\mathcal{S}(\mathfrak{A})$; we then define the 
Hodge functional derivatives of $F$ by the formula
\[
\ddt F(\varphi_\alpha+t\rho_\alpha)=\sqbr{\rho_\alpha}{\lmalpha{F}}=\sqbr{\rmalpha{F}}{\rho_{\alpha}}.
\]
It follows from the definition that, for a homogeneous functional $F$, 
the Hodge functional derivatives w.r.t.\ $\varphi_\alpha$ lie in $\Omega^{\dg \varphi_\alpha}(M, \ad P)$ and possess ghost number equal to $\gh F-\gh \varphi_\alpha$.
We have now at our disposal the essential elements to construct the BV Laplacian.
We start  defining the BV Laplacian of an element of $\mathcal{S}(\mathfrak{A})$ by the formula 
\begin{equation}\label{defBVlapl}
\DBV F:=\sum_\alpha (-1)^{\gh \phi^\alpha} \sqbr{\lmf{}}{\lma{F}}.
\end{equation}
The result is again a functional in $\mathcal{S}(\mathfrak{A})$, and, if $F$ is homogeneous, then $\DBV F$ is homogeneous of ghost number $\gh F +1$.

\begin{Rem}
This definition can also be extended to functionals in the space 
$\mathcal{S}^*(N;E)$ in analogy with the construction presented in
the preceding subsection. 
For a homogeneous functional $G$ in $\mathcal{S}^*(N;E)$, 
$\DBV G$ is again a functional in $\mathcal{S}^*(N;E)$, whose ghost number is given by $\gh G +1$ and whose degree is unchanged.
\end{Rem}
\begin{Rem}\label{rem-DBV}
Turning to a unified notation $\mathcal{S}$, we have in general
\[
\DBV:\mathcal{S}^{p,g}\to\mathcal{S}^{p,g+1}.
\]
Notice however that $\DBV$ is not well-defined for all functionals in
$\mathcal{S}$. It is particularly ill-defined on local functionals.
The correct definition would involve some regularization.
{\em We assume however that, independently of the regularization,
$\DBV F=0$ whenever $F$ depends only on one element in each 
pair field--antifield, as the formal definition of $\DBV$ suggests.}
\end{Rem}

The properties of the BV Laplacian $\DBV$ are:
\begin{itemize}
\item the BV Laplacian is a coboundary operator, i.e.\ 
\[
\Delta_{BV}^2=0;
\]
\item the BV antibracket measures the failure of the BV Laplacian to be a derivation, i.e.\
\begin{equation}\label{fail-Lapl}
\DBV (F G)=(\DBV F) G +(-1)^{\gh F}\BV{F}{G}+(-1)^{\gh F} F(\DBV G),
\end{equation}
where one of the functionals must be central.
\end{itemize}
The latter property in particular implies that the BV Laplacian is well-defined
on the subalgebra generated by those local functionals which are killed
by $\DBV$ (e.g., those described in the previous remark).
\begin{Rem}
If we restrict ourselves to $\mathcal{S}(\mathfrak{A})$,
then the above properties hold on the whole algebra.
A Gerstenhaber algebra 
with an operator $\Delta$ satisfying the above properties is known
as a BV algebra.
\end{Rem}
\begin{Rem}
We note that we can define (independently of the dimension) the BV antibracket
by
\[
\BV{F}{G}:=\sqbr{\rmf{F}}{\lma{G}}-\sqbr{\rma{F}}{\lmf{G}}.
\]
This is the definition of the BV antibracket in its original setting \cite{BV}.
This expression depends in general on the Riemannian metric on $M$,
but in the case of $BF$ theories
the antibracket is actually independent thereof
since it is equal to the one defined in subsection~\ref{subs-BVanti}.
\end{Rem}

\subsection{BV cohomology and observables}
We have introduced the BV Laplacian in order to deal with the quantum version of the
BV formalism, which is needed when considering functional integrals with weight $\exp(\ii/\hslash) S$, where 
$S$ should be a solution of the {\em{quantum master equation}}
\[
\BV{S}{S}-2\ii\hslash \DBV S=0.
\]

The main consequence of the quantum master equation is that the operator
\begin{equation}\label{OBV}
\OBV:=\dBV-\ii\hslash \DBV
\end{equation}
is a coboundary operator of ghost number $1$; it is {\em{not}}\/ a differential, because of (\ref{fail-Lapl}).
This operator is fundamental in the BV formalism; namely, all the meaningful observables in 
the BV formalism lie in the $0$-ghost number cohomology of $\OBV$.
This means (at least formally) that the vacuum expectation values of
$\OBV$-cohomology classes, weighted by $\exp (\ii/\hslash) S$,
are independent of the choice of gauge fixing. 
In turn, the v.e.v.s of trivial $\OBV$-cohomology classes or
of classes of ghost number different from zero
vanish.

We will show that the BV action $S$ of $BF$ theories,
to be introduced in \eqref{SBF},
satisfies separately the equations
\[
\DBV S=0\quad\text{and}\quad\BV{S}{S}=0,
\]
which imply that $S$ satisfies the quantum master equation.

\section{The BV superformalism for $BF$ theories}\label{sec-BVsuper}
The aim of this section is to define a new type of BV antibracket, 
which will allow us to obtain
the BV action for $BF$ theories in a simple way and to write it in
a compact form.

{}From now on we consider a new grading on the space of functionals
$\mathcal{S}$ called
the total
degree, which is defined as the sum of the form degree and the
ghost number; we will denote the total degree of a form $\alpha$ with degree
$\dg\alpha$ and ghost number $\gh \alpha$ by $\abs{\alpha}:=\dg \alpha+\gh
\alpha$; by homogeneous we will mean homogeneous w.r.t.\ the total
degree.
 
We note now that all the fields $\Big\{c^{+}; a^{+}; B;
\tau_1;\dots;\tau_{m-2}\Big\}$ have total degree $m-2$, while all the
remaining fields $\Big\{\tau_{m-2}^{+};\dots;\tau_{1}^{+};
B^+; a; c\Big\}$ have total degree equal to $1$. 
Here $a$ is the difference between $A$ and a given background connection
$A_0$ as in Assumption~\ref{Ass1}; for notational consistency, we denote
by $a^+$ the associated antifield.
We can cast all the fields in two homogeneous
superforms which we will denote by $\sfB$ and $\sfA$: 
\begin{align} 
\sfB&:=\sum_{k=1}^{m-2} (-1)^{\frac{k(k-1)}{2}} \tau_k +B+(-1)^m a^++ c^+,
\label{sfB}\\
\sfA&:=(-1)^{m+1} c+A+(-1)^m B^++\sum_{k=1}^{m-2} (-1)^{\frac{k(k-1)}{2}+m(k+1)} \tau_k^+.\label{sfA}
\end{align} 
Further, we define $\sfa:=\sfA-A_0$.

We refer to Appendix~\ref{app-sign}, for the definitions
of the dot product $\cdot$ and of the dot Lie bracket $\lb{\ }{\ }$.
We only recall that the dot structures make the algebra $\mathcal{S}$
into a superalgebra w.r.t.\ the total degree.
Analogously, we define the dot version  $\dbraket{\ }{\ }$ of 
the bilinear form $\braket{\ }{\ }$ on $\Omega^*(M,\ad P)$ by
\[
\dbraket\alpha\beta:=(-1)^{\gh\alpha\dg\beta}\braket\alpha\beta.
\]
It satisfies
\[
\dbraket\beta\alpha = (-1)^{\abs\alpha\abs\beta}\dbraket\alpha\beta.
\]

\subsection{The space of functionals $\mathcal{S}_{\sfA,\sfB}$}\label{sec-ser}
As in the previous section, we consider the algebra generated by
local functionals in the fields taking values in a
commutative algebra $\mathfrak{A}$
or in a de~Rham complex $\Omega^*(N;E)$.
However, from now on we will restrict ourselves only to those
functionals which depend on the linear combinations $\sfA$ and $\sfB$
(and not on the component fields). We will denote these algebras by
$\pser(\mathfrak{A})$, resp.\ $\pser(N;E)$, or generically by $\pser$.

We give $\pser$ a unique grading, by defining
the degree of a monomial in the superfields $\sfA$ and
$\sfB$ to be the sum of the total degrees of its factors.

Since the superform $\sfa$ has total degree $1$ and lies in 
$\Omega^*(M,\ad P)$, we
can consider $\sfA$ as a superconnection in the sense of Mathai and 
Quillen \cite{MQ}.  
With the help of the dot Lie bracket (see Appendix~\ref{app-sign}), 
we can then define the
covariant derivative of $\sfB$ w.r.t.\ the superconnection $\sfA$
and the supercurvature $\sfF_\sfA$ by: 
\begin{align*} 
\dA \sfB &:= \dd_{A_0} \sfB + \lb{\sfa}{\sfB}, \\ 
\sfF_{\sfA} &:=\dd_{A_0} \sfa + \frac{1}{2}
\lb{\sfa}{\sfa}.
\end{align*} 
Notice that the supercurvature would contain the extra term 
$F_{A_0}$ if the background connection $A_0$ were not chosen to be flat.
Note that in this new context the exterior and covariant derivatives are
operators of total degree $1$.

\subsubsection{The super functional derivatives}\label{sec-supfder}
We begin by introducing the  super test forms
$\rho_{\sfa}$ and $\rho_{\sfB}$: the super test form
$\rho_{\sfa}$ is defined to be the sum of the test forms corresponding to
the fields that appear in the superform $\sfa$, with the same sign
convention as in \eqref{sfA}; analogously we define the super test form $\rho_{\sfB}$. By
definition, the super test forms have then total degree $1$, resp.\ $m-2$.  
We then define the super functional derivatives of an element $F$ in $\pser(\mathfrak{A})$ by:
\begin{align*}
\ddt F(\sfa+ t \rho_{\sfa};\sfB) =\int_{M} \dbraket{\rA{F}}
{\rho_{\sfa}}&=\int_{M} \dbraket{\rho_{\sfa}}
{\lA{F}}, \\
\ddt F(\sfa;\sfB+ t \rho_{\sfB}) =\int_{M} \dbraket{\rB{F}}
{\rho_{\sfB}}&=\int_{M} \dbraket{\rho_{\sfB}}
{\lB{F}}. 
\end{align*}
It is then easy to determine the
total degree of the super functional derivatives of $F$; in fact, the following
identities hold:
\begin{equation}
\begin{aligned} 
\Abs{\lA{F}}&=\Abs{\rA{F}}=\abs{F}+m -1,\\   
\Abs{\lB{F}}&=\Abs{\rB{F}}=\abs{F}+2. 
\end{aligned}\label{totdeg}
\end{equation}
It will be also useful to express the right derivative of the functional $F$
in terms of the left one, and vice versa. The result of this 
computation is given by: 
\begin{align*} 
\lA{F} &= (-1)^{\abs{F}+m-1} \rA{F},\\ 
\lB{F} &= (-1)^{\abs{F}m} \rB{F}.
\end{align*} 
 
We next define the super functional derivatives of an element $F$ of
$\pser(N;E)$ by:
\begin{align*} 
\ddt F(\sfa+ t
\rho_{\sfa};\sfB)&=\tilde{\pi}_{2*}\dbraket{\tilde{\pi}_1^{*}                    
\rho_{\sfa}}{\lA{F}}=\pi_{1*}\dbraket{\rA{F}}{\pi_2^{*} \rho_{\sfa}}; \\ 
\ddt F(\sfa;\sfB+ t\rho_{\sfB})&=\tilde{\pi}_{2*}\dbraket{\tilde{\pi}_1^{*}\rho_{\sfB}}{\lB{F}}=\pi_{1*}\dbraket{\rB{F}}{\pi_2^{*}\rho_{\sfB}}.  
\end{align*} 
Their total degrees are still given by \eqref{totdeg}.  

\subsubsection{The super BV antibracket} 
Let us pick two functionals $F$ and
$G$ in $\pser(\mathfrak{A})$; then the super BV antibracket is defined
by: 
\begin{equation}\label{SBV} 
\sbv{F}{G}:= \int_{M} \dbraket{\rB{F}}{\lA{G}} -(-1)^{m}
\dbraket{\rA{F}}{\lB{G}}. 
\end{equation}
Note that the BV antibracket of $F$ and $G$ is again a functional in $\pser(\mathfrak{A})$.

Next we consider a functional $F$ in $\pser(\mathfrak{A})$ 
and a functional $G$ in $\pser(N;E)$, with $E$ an $\mathfrak{A}$-module;
we define the BV antibracket of $F$ and $G$ by:
\begin{equation}
\sbv{F}{G}=\tilde{\pi}_{2*}\dbraket{\tilde{\pi}_1^{*}\Big(\rB{F}\Big)}{\lA{G}} -(-1)^{\dim
M}\tilde{\pi}_{2*}\dbraket{\tilde{\pi}_1^{*}\Big(\rA{F}\Big)}{\lB{G}}.  
\label{SBVmixed}
\end{equation}
In this case the BV antibracket of $F$ and $G$ is  a
functional in $\pser(N;E)$.

We finally define the
BV antibracket of two functionals $F$ and $G$ in $\pser(N;E)$ by: 
\begin{equation*}
\sbv{F}{G}:=\pi_{13*}\dbraket{\pi_{12}^{*}\Big(\rB{F}\Big)}
{\pi_{23}^{*} \Big(\lA{G}\Big)} -(-1)^{m} \pi_{13}^{*} \dbraket{\pi_{12}^{*}\Big(\rA{F}\Big)}
{\pi_{23}^{*} \Big(\rB{G}\Big)}.
\end{equation*}
In this case we obtain that $\sbv{F}{G}$ is a functional in $\pser(N,
E)$.

The antibracket, in all the above cases, has total degree $1$; i.e.,
if we denote generically by $\ser$ the space of functionals 
and by $\ser^k$ the subspace of homogeneous elements of total degree $k$,
then
\[
\sbv{\ }{\ }:\ser^k\otimes\ser^l\to\ser^{k+l+1}.
\]
{}From now on we will use the short notation given in (\ref{SBV}) 
for all types of functionals that we have
discussed until now, and we omit in all cases the specific notation, 
leaving to the reader the understanding of the real
meaning of the formula. 
 
\subsection{Main properties of the super BV antibracket}
One could now wonder if there is a relationship between the super
BV antibracket defined in the previous subsection
and the BV antibracket defined in~\ref{subs-BVanti}
that we have discussed in the previous
subsection. 
We begin by explaining this relationship for the case of functionals
 in $\ser(\mathfrak{A})$.           
\begin{Lem}\label{SBV-BV} 
Suppose that we have two functionals $F$ and $G$ in $\ser(\mathfrak{A})$; then the following 
identity holds: 
\begin{equation} 
\sbv{F}{G}=\BV{F}{G}. 
\end{equation}   
\end{Lem} 

\begin{proof}
We prove the identity for homogeneous functionals; the general case follows by linearity.
We begin by computing the functional derivatives of $F$ and $G$:
\begin{align*}
\ddt F(\sfa+t \rho_{\sfa};\sfB)=\ddt F(a+t\rho_{a};c+t\rho_{c};\dots;\sfB) =\int_{M} \dbraket{\rA{F}}{\rho_{\sfa}}. 
\end{align*}
Next, we note that the integral selects the part of the integrand whose form degree in $M$ is equal to $m$, and that the super test form $\rho_{\sfa}$ is the sum of the usual test forms (with some signs to be considered). We write $\rho_{\sfa}$ as 
\begin{align*}
\rho_{\sfa}=\sum_{i=0}^{m} \sigma_{\sfa_i} \rho_{\sfa_i},
\end{align*} 
where by  $\rho_{\sfa_i}$ we denote the degree $i$ part of $\rho_\sfa$; i.e.,
$\rho_{\sfa_0}=\rho_c$, $\rho_{\sfa_1}=\rho_a$ and so on.
The signs $\sigma_{\sfa_i}$ are the same as in the definition \eqref{sfA} of $\sfA$;
namely, $\sfa=\sum \sigma_{\sfa_i}\sfa_i$. Similarly we introduce signs $\sigma_{\sfB_j}$
as in $\sfB=\sum\sigma_{\sfB_j}\sfB_j$.
We can then write:
\begin{equation}\label{deriv}
\begin{aligned}
\integ{\rA{F}}{\rho_{\sfa}}&=\integw{\rfunc{F}{c}}{\rho_c}+\integw{\rfunc{F}{a}}{\rho_a}+\dots\\
&=\sum_{j=0}^{m} \sigma_{\sfa_j}\integ{\Big(\rA{F}\Big)_{m -j}}{\rho_{\sfa_j}},
\end{aligned}
\end{equation}
where the subscript denotes the restriction to the term of the indicated form degree.
We recall that $\gh{\rho_{\sfa_j}}=1-j$; then we obtain e.g.\ for the $j$-th term in the last expression of the above identity (recalling the definition of the total degree of the functional derivative of $F$ w.r.t.\ $\sfa$):
\[
\integ{\Big(\rA{F}\Big)_{m -j}}{\rho_{\sfa_j}}
=(-1)^{(\abs{F}+j-1)j} \integw{\Big(\rA{F}\Big)_{m -j}}{\rho_{\sfa_j}}
\]
By confronting the two expressions in (\ref{deriv}), and doing similar computations
in the other cases, we obtain for $j=0,\dots,m$:
\begin{align*}
\Big(\rA{F}\Big)_{m -j}&=\sigma_{\sfa_j}(-1)^{(\abs{F}+j-1)j}\rfunc{F}{\sfa_j},
&\quad \Big(\rB{F}\Big)_{m -j}&=\sigma_{\sfB_j}(-1)^{(m-2-j)(m-j)}\rfunc{F}{\sfB_j},\\
\Big(\lA{F}\Big)_{m -j}&=\sigma_{\sfa_j}(-1)^{(1-j)(m-j)}\lfunc{F}{\sfa_j},
&\quad \Big(\lB{F}\Big)_{m -j}&=\sigma_{\sfB_j}(-1)^{(\abs{F}-m+2+j)j}\lfunc{F}{\sfB_j},
\end{align*}
We cast then all these expressions in the definition of the super BV antibracket
\[
\dbraket{\rB{F}}{\lA{G}}. 
\]
Then we use the above expressions, and, after rewriting $\dbraket{}{}$ as $\braket{}{}$, we compute the products $\sigma_{\sfB_{m-j}} \sigma_{\sfa_{j}}$, separately for the case $m$ even and $m$ odd. In order for the superbracket to coincide
with the ordinary bracket, these products must be 
\begin{align*}
\sigma_{\sfB_{m-i}} \sigma_{\sfa_{i}}&=
\begin{cases}
-1 & \text{if $i=0$,}\\
1 &\text{otherwise}
\end{cases}\\
\intertext{for $m$ even, and}
\sigma_{\sfB_{m-i}} \sigma_{\sfa_{i}}&=
\begin{cases}
(-1)^i &\text{for $i=0,1$,}\\
(-1)^{i+1} &\text{otherwise}
\end{cases}
\end{align*}
for $m$ odd.
It can be readily computed that the choice of signs made in \eqref{sfB} and in \eqref{sfA} is consistent  with the above rules; therefore, the proof then follows.
\end{proof}

For the general case of elements of $\ser(N;E)$, 
the above rule must be slightly modified.  
We begin by noting that any homogeneous $F$ of total degree $\abs F$
in this algebra can be written in the form
\begin{displaymath}                    
F=\sum_l F_l, 
\end{displaymath} 
where $F_l$ is an element of $\mathcal{S}^{|F|-l,l}(N;E)$. This is obtained by expanding
the superfields in their components.
We are now ready to state the following 
\begin{Lem}\label{SBV-BVII} 
Let $F$ and $G$ be homogeneous elements of $\ser(N;E)$.
If we expand them according to the above rule
\begin{align*} 
F=\sum_k F_k\ \ \text{and} \ \ G=\sum_l G_l, 
\end{align*}   
then the following identity holds: 
\begin{equation}\label{bvsbv} 
\sbv{F}{G}=\sum_{k,l} (-1)^{(\abs F - k+1)l} \BV{F_k}{G_l}. 
\end{equation} 
\end{Lem} 
 
\begin{proof}
The proof of this identity is similar to the proof of Lemma~\ref{SBV-BV}; in fact, we have to compute the functional derivatives of $F$ and $G$ w.r.t.\ $\sfa$ and $\sfB$, and express them via the functional derivatives w.r.t.\ the usual fields of the theory.
We therefore recall the formulae for the functional derivatives, and we apply them to $F$, obtaining:
\begin{multline}\label{derivfun}
\ddt F(\sfa+t\rho_{\sfa};\sfB)=\push{\rA{F}}{\rho_{\sfa}}=\\
=\sum_{j=0}^{m} \sigma_{\sfa_j}\push{\Big(\rfunc{F}{\sfa}\Big)_{m-j}}{\rho_{\sfa_j}}
=\sum_l \ddt F_l(a+t\rho_{a};c+t\rho_{c};\dots)=\\
=\sum_l \sum_{j=0}^{m} \ddt F_l(\sfa_j+t\rho_{\sfa_j})
=\sum_l 
\sum_{j=0}^{m} \pushw{\rfunc{F_l}{\sfa_j}}{\rho_{\sfa_j}}=\\
=\sum_{j=0}^{m} \pushw{\sum_l \rfunc{F_l}{\sfa_j}}{\rho_{\sfa_j}}.
\end{multline}
Then the following holds, if we go from the dot product to the ordinary product:
\begin{align*}
(-1)^{(\abs F -l-1+j)j}\braket{\rfunc{F_l}{\sfa_j}}{\pi_2^{*}\rho_{\sfa_j}}=\dbraket{\rfunc{F_l}{\sfa_j}}{\pi_2^{*}\rho_{\sfa_j}}.
\end{align*}
By confronting the terms in (\ref{derivfun}), and operating similarly for the other cases,
we obtain the following identities for $j=0,\dots,m$:
\begin{align*}
\Big( \rfunc{F}{\sfa}\Big)_{m-j}&=\sum_l \sigma_{\sfa_j}(-1)^{(\abs F -l-1+j)j} \rfunc{F_l}{\sfa_j},\\ 
\Big( \rfunc{F}{\sfB}\Big)_{m-j}&=\sum_l \sigma_{\sfB_j}(-1)^{(\abs F -l-m+j)j} \rfunc{F_l}{\sfB_j},\\
\Big( \lfunc{F}{\sfa}\Big)_{m-j}&=\sum_l \sigma_{\sfa_j}(-1)^{(1-j)(l-m+j)} \lfunc{F_l}{\sfa_j},\\ 
\Big( \lfunc{F}{\sfB}\Big)_{m-j}&=\sum_l \sigma_{\sfB_j}(-1)^{(m-j)(l-m+j)} \lfunc{F_l}{\sfB_j}.
\end{align*}
We can finally cast all these expressions in the explicit formula for the super BV antibracket, and, by recalling the explicit values of the chosen signs $\sigma_{\sfa_j}$ and $\sigma_{\sfB_j}$, we can finally obtain the desired identity (recall the form degree selection rule imposed by the 
pushforwards). 
\end{proof}

We note that for the case in which $F$ is in $\ser(\mathfrak{A})$  
and $G$ is in $\ser(N;E)$ for an $\mathfrak{A}$-module $E$, 
then the following identity holds: 
\begin{equation}
\sbv{F}{G}=\sum_{l} (-1)^{(\abs{F}+1)l} \BV{F}{G_l}; 
\label{SBV-BViii}
\end{equation}
this formula will play a special r\^ole in some later computations (we skip the proof of this identity, because it is in principle the same as for the two previous Lemmata). 

Let us now extend the super BV antibracket $\sbv{\ }{\ }$ to the whole of
$\mathcal{S}$ by the following rule  
\[
\sbv{\alpha}{\beta}:=(-1)^{(\gh \alpha+1)\dg \beta}\BV{\alpha}\beta,
\] 
with $\alpha$ and $\beta$ homogeneous elements of $\mathcal{S}$.
Recalling the graded commutativity rule, the graded Leibnitz rule and the graded Jacobi rule for $\BV{\ }{\ }$, one can then show the following properties
of the super BV antibracket $\sbv{\ }{\ }$:
\begin{itemize}
\item $\sbv{\alpha}{\beta}=-(-1)^{(\abs{\alpha}+1)(\abs{\beta}+1)}\sbv{\beta}{\alpha}$,\\
 whenever one of the two elements is central in $\mathcal{S}$.
\item $\sbv{\alpha}{\beta \gamma}=\sbv{\alpha}{\beta} \gamma+(-1)^{(\abs{\alpha}+1)\abs{\beta}}\beta\sbv{\alpha}{\gamma}$,\\
whenever $\alpha$ or $\beta$ or both are central in $\mathcal{S}$.
\item $(-1)^{(\abs{\alpha}+1)(\abs{\gamma}+1)}\sbv{\alpha}{\sbv{\beta}{\gamma}}+\text{cyclic permutations}=0$,\\ 
whenever two of the three elements are central in $\mathcal{S}$.
\end{itemize}
Here we have used the previous notational convention for the total degree.
In particular if we restrict to $\ser$, by linearity
the previous identities hold if we replace $\alpha$, $\beta$ and $\gamma$
with elements $F$, $G$ and $H$ of $\ser$.

For central elements in $\pser$ we can take e.g.\ any functional $F$ in $\pser(\mathfrak{A})$, while considering as more general elements in $\pser(N;E)$, for an $\mathfrak{A}$-module $E$.
(We have omitted the products between elements in $\pser$, but it is understood that we are dealing with the shifted dot product.)
Let us now pick a central element $S$ of $\pser$ with even total degree; we then define an operator $\boldsymbol{\delta}$ on the superspace $\pser$ by 
\[
\boldsymbol{\delta}:=\sbv{S}{\ };
\]
since $S$ has even total degree, $\boldsymbol{\delta}$ is an odd derivation by the above identities.
{}From Lemma~\ref{SBV-BVII},~\ref{push1},~\ref{pull1},~\ref{extd1} and~\ref{trace1} we can derive the useful properties of $\delta_{S}$:              
\begin{Cor}\label{pull2} 
Suppose that the functional $F$ lies in $\pser(N;E)$,
and that we have a map $h$ from some manifold $H$ to the
manifold $N$, then the following identity holds:
\begin{equation*} 
\boldsymbol{\delta} [h^{*}(F)] = h^{*}(\boldsymbol{\delta} F). 
\end{equation*} 
\end{Cor} 
 
\begin{Cor} \label{push2} 
Suppose that we have a homogeneous functional $F$ in $\pser(H;E)$, where
$E$ is a real or complex algebra and H is the total space $H$ of a bundle over $N$ with typical fiber $B$. The integration along fiber of the functional $F$ gives a functional of the same type,
defined on the manifold $N$ and with total degree $\abs{\pi_{*}(F)}=\abs{F} -\dim B$.   
Then we obtain the following identity: 
\begin{equation*} 
\boldsymbol{\delta} [\pi_{*}(F)] = (-1)^{\dim B} \pi_{*} (\boldsymbol{\delta} F) 
\end{equation*} 
\end{Cor} 
 
\begin{Cor} \label{extd2} 
Let us pick a functional $F$ in $\pser(N;E)$, for $N$ and $E$ as in the preceding Lemma. 
Let us denote by $\dd$ the exterior derivative on the manifold $N$. 
Then the following identity holds: 
\begin{equation*} 
\boldsymbol{\delta}(\dd F)=-\dd(\boldsymbol{\delta} F).  
\end{equation*} 
\end{Cor} 
 
\begin{Cor} \label{trace2} 
Let us suppose that the functional $F$ belong to the superalgebra
$\pser(N,\Lg)$; let us suppose additionally that we
have a $\mathfrak{g}$-module $V$.
The application of the trace to $F$ gives an element of $\pser(N;\mathbb{R})$ (or of $\pser(N;\mathbb{C})$, depending on whether $V$ is a real or complex module).  
Then we obtain the following identity: 
\begin{equation*} 
\boldsymbol{\delta}[\tr_{\rho}(F)] = \tr_{\rho}(\boldsymbol{\delta} F). 
\end{equation*} 
\end{Cor} 

\subsection{The super BV Laplacian}
In analogy with what we have done for the BV antibracket,
let us introduce a ``twisted'' version $\sDBV$ of the BV Laplacian
on the superalgebra $\mathcal{S}$, 
endowed with the two usual gradings (the form degree and the ghost number).
We define the {\emph{super BV Laplacian}} by
\[
\sDBV \alpha:=(-1)^{\dg \alpha} \DBV \alpha,
\]
for $\alpha\in \mathcal{S}$.
Since the BV Laplacian is nilpotent, it follows immediately that the super BV Laplacian is nilpotent, too.
Let us take $\alpha$ and $\beta$ in $\mathcal{S}$, and let us suppose that at least one of the two elements is central in $\mathcal{S}$. It follows then from (\ref{fail-Lapl})
that
\begin{equation}
\sDBV (\alpha\cdot\beta)= (\sDBV \alpha)\cdot\beta +(-1)^{\abs{\alpha}}\sbv{\alpha}\beta 
+(-1)^{\abs{\alpha}}\alpha\cdot (\sDBV\beta),
\label{sfail-Lapl}
\end{equation}
where $\alpha$ or $\beta$ must be central.
Restricting to the super algebra $\pser$, it follows easily that the super BV Laplacian is a 
coboundary operator
\[
\sDBV:\pserr{k}\to \pserr{k+1}
\]
which satisfies \eqref{sfail-Lapl} with $\alpha$ and $\beta$ in $\ser$.
The BV operator $\OBV$ defined in \eqref{OBV} is replaced in the superformalism
by the operator
\[
\sOBV = \sfdelta -\ii\hslash\sDBV.
\]
As a consequence of the general case, $\sOBV$ is a coboundary operator.%, and
%$\sOBV$-closed functionals have
%v.e.v.s independent of the gauge-fixing, while $\sOBV$-exact functionals
%have vanishing v.e.v.s.

\subsection{Twists}\label{ssec-twist}
Let $\sfO$ be an even element of $\ser$. We define the twisted BV coboundary
operator by
\[
\widetilde\sOBV = \exp\left(-\frac\ii\hslash\sfO\right)\,
\sOBV\,\exp\left(\frac\ii\hslash\sfO\right)=
\sOBV + \bde_\sfO + \frac\ii\hslash\Phi_\sfO,
\]
with $ \bde_\sfO:=\sbv\sfO{\ }$, and
\[
\Phi_\sfO=\sOBV\sfO+\frac12 \sbv\sfO\sfO
\]
as a multiplication operator.

\begin{Def}\label{def-flobs}
We call {\em flat}\/ an even functional $\sfO$ with
$\Phi_\sfO=0$, {\em flat observable}\/ an $\sOBV$-closed flat functional,
and {\em flat invariant observable}\/ a $\sfdelta$-closed flat observable.
\end{Def}
A basic fact that we will need in Sections~\ref{sec-gWl} and
\ref{sec-lobs} is expressed by the following
\begin{Lem}\label{lem-flobs}
If $\sfO$ is a flat observable, then so is $\lambda\sfO$ for any constant 
$\lambda$; moreover,
$\bde_\sfO$ is a superdifferential
(of degree $\abs\sfO+1$)
which anticommutes with $\sOBV$. If we also assume that $\sfO$ is invariant, 
then 
$\bde_\sfO$ anticommutes with $\sfdelta$, so $\sfdelta_\lambda:=
\sfdelta + \lambda\bde_\sfO$ is an odd differential for all $\lambda$.
\end{Lem}
\begin{proof}
By definition, a flat observable $\sfO$ satisfies separately
\[
\sOBV\sfO=0\qquad\text{and}\qquad
\sbv\sfO\sfO=0.
\]
This implies that $\Phi_{\lambda\sfO}=0$ for all $\lambda$.
The second equation above together with the Jacobi identity implies that
$\bde_\sfO$ is a coboundary operator.
Since $\widetilde\sOBV$ and $\sOBV$ square to zero and $\Phi_\sfO=0$, we obtain
\[
\sOBV\bde_\sfO + \bde_\sfO\sOBV + \bde_\sfO^2=0.
\]
The second claim then follows since $\bde_\sfO$ squares to zero.
For $\sfO$ invariant, we also have $\sbv S\sfO=0$. So by Jacobi we obtain
$\sfdelta\bde_\sfO + \bde_\sfO\sfdelta=0$.
\end{proof}

\section{The BV action for $BF$ theories}\label{sec-BVaction}
We have now at our disposal all the tools needed to write down the correct BV action for $BF$ theories.
Namely, we claim that it is given by
\begin{equation} 
\SBV:=\int_{M} \dbraket{\sfB}{\sfF_{\sfA}}.          
\label{SBF}
\end{equation}
\begin{Rem}
Earlier versions of this form for the BV action of $BF$ theories can be found
in \cite{Wall, Ikem}, where however proofs were not given and, 
in particular, there was no explicit treatment of the sign conventions
(i.e., our ``dot'' structures).
Special cases (with explicit signs) were also considered in \cite{C2,BF7}.
In particular, the structure of the BV action in terms of superfields is an agreement with the
general pattern described in \cite{DG}. See also \cite{W3} and \cite{AKSZ}
for the case of Chern--Simons theory.
\end{Rem}
This form of the BV action holds not only for the $BF$ theories described in the previous
sections but also for the ``canonical $BF$ theories'' pointed out in Remark~\ref{rem-cBF}
(observe that the two-dimensional case has already been considered in \cite{CF}).

We divide the proof, for the ordinary case,
in two steps: $i)$ we show that the above functional is a solution of the master equation
corresponding to the $BF$ action \eqref{BFaction} with infinitesimal 
symmetries \eqref{BRSTABc} and \eqref{BRSTtau}
(subsection~\ref{sec-mastereq});
$ii)$ we show that it is $\DBV$-closed (subsection~\ref{sec-LaplBF}).
In subsection~\ref{sec-dual} we will then give the proof in the case of canonical
$BF$ theories.

\subsection{The master equation}\label{sec-mastereq}
We begin with the statement of the main Theorem, and we devote the rest of the section to its proof and to some important consequences.

\begin{Thm} \label{SME}  
The following identity holds: 
\begin{align*} 
\sbv{S}{S}=0. 
\end{align*} 
\end{Thm} 
 
\begin{Rem} 
By lemma~\ref{SBV-BV}, the above result is equivalent to the
statement that the action $S$ satisfies the usual ME w.r.t.\ the usual 
BV antibracket.
\end{Rem} 

\begin{proof}
We begin by computing the left and right partial derivatives w.r.t.\ $\sfa$ and $\sfB$; e.g.\ the left partial derivative of $S$ w.r.t.\ $\sfa$ is given by
\begin{multline*}
\ddt S(\sfa+t\rho_{\sfa};\sfB)
=\integ{\sfB}{\dd_{\sfA}\rho_{\sfa}}=\\
=(-1)^{m-1}\integ{\dd_{\sfA} \sfB}{\rho_{\sfa}}=\integ{\rho_{\sfa}}{\dd_{\sfA} \sfB}.
\end{multline*}
It follows that:
\begin{align}
\rA{S}&=(-1)^{m-1}\lA{S}= (-1)^{m-1}\dd_{\sfA} \sfB;
\label{dSdA}\\
\intertext{similarly}
\rB{S}&=\lB{S}=\sfF_{\sfA}.\label{dSdB}
\end{align}
If we now insert the above functional derivatives in the formula for the BV antibracket, 
we obtain
\begin{align*}
\sbv{S}{S}&=2 \integ{\sfF_{\sfA}}{\dd_{\sfA}\sfB}=2\int_M \dd\dbraket{\sfF_{\sfA}}{\sfB}-2\integ{\dd_{\sfA} \sfF_{\sfA}}{\sfB},
\end{align*}
by the invariance of $\braket{\ }{\ }$ ($\sfA$ is a superconnection).
The first term vanishes by Stokes' Theorem, and the second by the super Bianchi identity
\[
\dd_{\sfA}\sfF_{\sfA}=0.
\]
So the claim follows.
\end{proof}
Since $S$ satisfies the ME, the Leibnitz rule and the Jacobi identity for the super BV antibracket imply that the operator
\[
\boldsymbol{\delta}:=\sbv{S}{\ },
\]
defined on $\pser(N;E)$,
is an odd differential. In many of the forthcoming computations we need the following
\begin{Prop}
The action of $\sfdelta$ on the superfields 
$\sfa$ and $\sfB$ is given by:
\begin{align}\label{deltaa}
\boldsymbol{\delta} \sfa&=(-1)^m \sfF_{\sfA}\\
\intertext{and }
\label{deltab}
\boldsymbol{\delta} \sfB&=(-1)^m \dd_{\sfA}\sfB.
\end{align}
\end{Prop}
\begin{proof}
The above formulae follow from \eqref{SBVmixed}. Let us begin with $\sfdelta\sfa$:
\[
\sfdelta\sfa(x) = \int_M \dbraket{\rB S}{\lA\sfa(x)}.
\]
By definition however
\[
\int_M\dbraket{\rho}{\lA\sfa(x)} = \ddt[\sfa(x)+t\rho(x)]=\rho(x),
\]
provided $\rho$ is a test form of total degree $1$. Since $\rB S$ has total
degree $2$, we cannot apply the above formula directly. We use then the following trick.
Let $\epsilon$ be a scalar of total degree $-1$. Then
\[
\epsilon\cdot\sfdelta\sfa(x) = (-1)^m \int_M \dbraket{\epsilon\cdot\rB S}{\lA\sfa(x)}=
(-1)^m\epsilon\cdot\rB S(x).
\]
Thus,
\[
\sfdelta\sfa(x)=(-1)^m\rB S(x)=(-1)^m\sfF_\sfA(x),
\]
where we have used \eqref{dSdB}.
Similarly, we have
\[
\sfdelta\sfB(x)=-\rA S(x)=(-1)^m
\dd_{\sfA}\sfB,
\]
by \eqref{dSdA}.
\end{proof}
Recalling the formula \eqref{SBV-BViii}
that expresses the super BV antibracket in term of the usual BV antibracket, we can now
recover the action of
the usual $\dBV$ operator defined by $\BV S{\ }$.
Namely,
\[
\boldsymbol{\delta} \sfa=\sum_{j=0}^{m} \sigma_{\sfa_j}(-1)^j \dBV \sfa_j.
\] 
By decomposing the expression for $\boldsymbol{\delta} \sfa$ in its homogeneous components and by confronting the two expressions, we get the action of the BV operator $\dBV$ on the fields. Similarly, we can recover the action of $\dBV$ on the components of $\sfB$.
By setting the antifields to zero, one obtains then that $\dBV$ on the
fields $\{A,B,c,\tau_1,\dots,\tau_{m-2}\}$ coincides with the BRST operator
given in \eqref{BRSTABc} and \eqref{BRSTtau}.
Moreover, 
it follows easily from the definition of $S$ and of the superforms $\sfa$ and $\sfB$ 
that the action reduces to the classical $BF$ action, if we set all antifields to $0$.
Thus, we have proved the following
\begin{Thm}
$\SBV$ is a solution of the master equation for the $BF$ theory.
\end{Thm}

\subsection{The $\DBV$-closedness of the BV action}\label{sec-LaplBF}
We now turn to the proof of the identity
\begin{equation}\label{laplace}
\DBV S=\boldsymbol{\Delta} S=0.
\end{equation}
First, we recall that 
$\frg$ is endowed with a nondegenerate, symmetric, invariant bilinear form $\braket{}{}$. 
We now choose a basis $\{X_k\}$ of $\frg$ such that $\braket{X_i}{X_j}=s_i\delta_{ij}$,
$s_i=\pm1$; in this basis we have the structure constants $f_{ij}^k$ given by
the relation
\[
\Lie{X_i}{X_j}=\sum_{k=1}^{\dim\Lg}f_{ij}^k X_k.
\]
We then introduce the symbols $\tilde f_{ij}^k$ as $s_kf_{ij}^k$.
Thus,
\[
\tilde f_{ij}^k=\braket{\Lie{X_i}{X_j}}{X_k}.
\]
{}From the non-degeneracy of $\braket{}{}$ one then gets the useful relation
\begin{equation}\label{strconst}
\tilde f_{ij}^k=-\tilde f_{ik}^j=-\tilde f_{kj}^i.
\end{equation}
If we write the BV action as a sum of local terms in the fields, 
we see from the very definition of the BV Laplacian $\DBV$ 
(see Remark~\ref{rem-DBV})
that the only terms in this sum which are not automatically $0$ have the form
\[
\sqbr{\phi_{\alpha}^*}{[\phi^\alpha,c]}\ ,
\]
for all $\alpha$ in the index set of the fields (this the only way to pair a field and its antifield that is allowed by the integration over $M$); we can rewrite it in the form (up to some sign)
\begin{align}\label{Hodge-BV}
\sqbr{\phi_{\alpha}^*}{[\phi^\alpha,c]}=
\tilde f_{jk}^i \cirbr{\phi_{\alpha}^{*,i}}{\phi^{\alpha,j}c^k}\ ,
\end{align}
with 
\begin{equation}\label{cirbr}
\cirbr{\alpha}{\beta}:=\int_{M}\alpha\wedge\star \beta,\qquad \alpha,\ \beta\in\Omega^*(M),
\end{equation}
and $\phi^{\alpha}=\sum\phi^{\alpha,i}X_i$ and similarly for $\phi_{\alpha}^{*}$ and $c$.
Now, by \eqref{strconst}, one sees that in the above formula no field component
is paired to the corresponding antifield component. So, by 
Remark~\ref{rem-DBV}, it is annihilated by the BV Laplacian.

\subsection{Canonical $BF$ theories}\label{sec-dual}
We start here a digression about the version of $BF$ theories mentioned
in Remark~\ref{rem-cBF}. The material covered in this subsection
is not essential for the rest of the paper and can be safely skipped.
Though, this kind of $BF$ theories is interesting by itself (and appears
in two-dimensions as a particular case of the Poisson sigma model
\cite{Ik,SS}). 

We recall now the basic idea:
since the curvature $F_{A}$ is a tensorial form of the adjoint type,
the most natural way to define a $BF$ theory is to choose $B$ of
the coadjoint type and to use the canonical
pairing between $\Lg^*$ and $\Lg$.
We consider then $B$ as a form in $\Omega^{m-2}(M,\ad^* P)$.
Observe that since we do not introduce a bilinear form on $\Lg$ anymore, 
Assumption~\ref{Ass3} is in this case meaningless. For simplicity we will
retain in this case as well Assumptions~\ref{Ass1} and~\ref{Ass2}.
We begin with some notations:
\begin{itemize}
\item {\em By $\braket{}{}$ we will denote in this subsection the canonical pairing between $\Lg^*$ and $\Lg$}\/; 
it can be naturally extended to a pairing between forms in $\Omega^p(M, \ad^* P)$ and forms in $\Omega^q(M,\ad P)$, and we will denote this pairing by the same symbol.
\item By $\{X_i\}$ we denote a basis of $\Lg$, while by $\{X^j\}$ we denote its dual basis: $\braket{X^i}{X_j}=\delta^i_j$.
\item By $f_{ij}^k$ we denote the structure constants w.r.t.\ the basis $\{X_l\}$, i.e.\ 
\[
f_{ij}^k=\braket{X^k}{\Lie{X_i}{X_j}}.
\]
\item By $\Ad^*$ we denote the coadjoint action of $G$ on $\Lg^*$;
i.e.\ $\braket{\Ad^*(g)\xi}X:=\braket{\xi}{\Ad(g^{-1})X}$.
\item by $\ad^*$ we denote the coadjoint action of $\Lg$ on $\Lg^*$; 
i.e., $\braket{\ad^*(X)\xi}{Y}=-\braket{\xi}{\ad(X)Y}$.
The coadjoint action can be extended to an action of forms in $\Omega^p(M,\ad P)$ on $\Omega^q(M,\ad^* P)$ in the usual way.
We only notice the sign rules for this extended coadjoint action
\begin{align*}
\ad^*(\Lie{\alpha}{\beta})\gamma&=\ad^*(\alpha)\ad^*(\beta)\gamma-(-1)^{\deg \alpha \deg \beta +\gh \alpha \gh \beta}\ad^*(\beta)\ad^*(\alpha)\gamma;\\
\braket{\ad^*(\alpha)\gamma}{\beta}&=-(-1)^{\deg \alpha \deg \gamma+\gh \alpha \gh \gamma}\braket{\gamma}{\Lie{\alpha}{\beta}} ,
\end{align*}
for $\alpha, \beta \in \Omega^*(M, \ad P)$ and $\gamma \in \Omega^*(M,\ad^* P)$, where we have implicitly supposed to consider forms with additional ghost number gradation.
\item Finally, we denote (improperly) by $\dd_A$ the covariant derivative acting on $\Omega^*(M,\ad^* P)$; it satisfies
\[
\dd \braket{\alpha}{\beta}=\braket{\dd_A \alpha}{\beta}+(-1)^{\deg \alpha} \braket{\alpha}{\dd_A \beta},
\]
where $\alpha\in \Omega^*(M,\ad^* P)$ and $\beta\in \Omega^*(M, \ad P)$, and
\[
\dd_A (\dd_A \alpha)=\ad^*(F_A)\alpha.
\]
\end{itemize}
With these conventions in mind, we define the canonical $BF$ action functional by
\[
S^{\mathrm{can}}:=\int_M \braket{B}{F_A}.
\]
The Euler--Lagrange equations are still given by \eqref{eqmotion},
where now the covariant derivative is understood to operate on 
$\Omega^{m-2}(M,\ad^* P)$.

We let the group $\Omega^0(M, G)\rtimes \Omega^{m-3}(M,\ad^* P)$ operate 
(from the right)
on $\mathcal{A}\times \Omega^{m-2}(M, \ad^* P)$ by the rule
\[
(A, B)(g, \tau):=(A^g, \Ad^*(g^{-1}) B+\dd_{A^g} \tau).
\]
It is then easy to verify that $S^{\mathrm{can}}$ is invariant under this
action.
The infinitesimal transformations then read
\[
\delta A = \dd_A c;\qquad \delta B= -\ad^*(c) B +\dd_A \tau.
\]
These symmetries present the same reducibility problems as in 
Section~\ref{sec-BFtheories}; 
therefore, we have to resort to the BV formalism here as well.

\subsubsection{The BRST and the BV formalism}
\label{sec-BRSTcan}
The BRST transformations corresponding to the reducible infinitesimal symmetries
in this case read
\begin{align*}
\dBRST A &= \dd_A c; &\qquad \dBRST B&= -\ad^*(c) B + \dd_A \tau_1;\\
\dBRST c &= -\frac{1}{2} \Lie{c}{c}; &\qquad
\dBRST \tau_k &=-\ad^*(c)\tau_k +\dd_A \tau_{k+1},\quad k=1,\dots,m-3;\\
 & &\qquad \dBRST \tau_{m-2}&=-\ad^*(c) \tau_{m-2}.
%\dBRST c &= -\frac{1}{2} \Lie{c}{c};\quad \dBRST A = \dd_A c;\quad \dBRST B= -\ad^*(c) B + \dd_A \tau_1;\\
%\dBRST \tau_k &=-\ad^*(c)\tau_k +\dd_A \tau_{k+1},\quad k=1,\dots,m-3;\\
%\dBRST \tau_{m-2}&=-\ad^*(c) \tau_{m-2}.
\end{align*}
Here $c$ denotes the Faddeev--Popov ghost, i.e.\ a form on the space of fields with values in $\Omega^0(M,\ad P)$ with ghost number $1$, and by $\tau_k$ we denote the ghosts for ghosts taking values in $\Omega^{m-2-k}(M, \ad^* P)$ and with ghost number $k$.
These BRST transformations present the same problems as in 
Section~\ref{sec-BFtheories}, namely $\dBRST$ is a differential only modulo terms containing the curvature of $A$. We have therefore to switch to the BV formalism.
We associate to each field $\phi^{\alpha}\in \{A,B,c,\tau_1,\dots,\tau_{m-2}\}$ a canonical antifield $\phi_{\alpha}^+$ following the rules
\begin{itemize}
\item if $\phi^{\alpha}$ takes values in $\Omega^{p^{\alpha}}(M, \ad P)$, resp.\ $\Omega^{p^{\alpha}}(M, \ad^* P)$, then its canonical antifield takes values in $\Omega^{m-p^{\alpha}}(M, \ad^* P)$, resp. $\Omega^{m-p^{\alpha}}(M, \ad P)$;
\item the ghost number of $\phi_{\alpha}^+$ is set to be equal to $-1-\gh \phi^{\alpha}$.
\end{itemize}
We define the total degree of a form $\alpha$ with degree $\dg \alpha$ and ghost number $\gh \alpha$ by
\[
\abs{\alpha}:=\dg \alpha +\gh \alpha.
\]
Accordingly to what we have done in Section~\ref{sec-BVsuper}, 
we define the {\emph{dot duality}} by the rule
\[
\dbraket{\alpha}{\beta}:=(-1)^{\gh \alpha \dg \beta} \braket{\alpha}{\beta},
\] 
for $\alpha$ an element of $\Omega^*(M,\ad^* P)$ with ghost number $\gh \alpha$ 
and $\beta$ in $\Omega^*(M,\ad P)$ with form degree $\deg\beta$.
The dot Lie bracket $\lb{\ }{\ }$ is defined analogously as in Appendix~\ref{app-sign}, and it enjoys the same sign rule.
We define additionally the {\emph{super coadjoint action}} of $\Omega^*(M, \ad P)$ on $\Omega^*(M, \ad^* P)$ by the rule
\[
\sad^*(\alpha)\beta:=(-1)^{\gh \alpha \dg \beta}\ad^*(\alpha)\beta.
\]
Without proof we write down some useful formulae, which are analogous to the formulae displayed in Appendix~\ref{app-sign} 
\begin{align*}
\sad^*(\lb{\alpha}{\beta})\gamma&=\sad^*(\alpha)\sad^*(\beta)\gamma-(-1)^{\abs{\alpha}\abs{\beta}} \sad^*(\beta)\sad^*(\alpha)\gamma,\\
\dbraket{\sad^*(\alpha)\gamma}{\beta}&=-(-1)^{\abs{\alpha}\abs{\gamma}}\dbraket{\gamma}{\lb{\alpha}{\beta}},
\end{align*}
for $\alpha, \beta \in \Omega^*(M, \ad P)$ and $\gamma\in \Omega^*(M, \ad^* P)$.
If $A$ is a connection on $P$, we also have
\[
\dd \dbraket{\gamma}{\alpha}=\dbraket{\dd_A \gamma}{\alpha}+(-1)^{\abs{\gamma}}\dbraket{\gamma}{\dd_A \alpha}.
\]
Finally, it is useful to write the duality pairing also in the opposite
order; as usual, one defines $\braket X\xi=\braket\xi X$ for $X\in\Lg$ and
$\xi\in\Lg^*$. When we extend the pairing to forms and then consider the dot 
version, we obtain the rule
\[
\dbraket{\gamma}{\alpha}=(-1)^{\abs{\gamma}\abs{\alpha}} \dbraket{\alpha}{\gamma}.
\]

We then define the functional derivatives w.r.t.\  all the fields of the theory and the BV antibracket $\BV{}{}$ by the same formulae as in 
subsection~\ref{sec-funderiv} (where the invariant, nondegenerate bilinear form $\braket{}{}$ is replaced by the duality pairing). 
This antibracket enjoys all usual properties of a BV antibracket.
Analogously, provided we have a solution $S$ of the master equation $\BV{S}{S}=0$, 
we define the BV differential $\delta$ by the rule
$\delta :=\BV{S}{\ }$.
This operator has the same properties of the previously introduced BV differential
(see subsection~\ref{sec-propBV}). 

Finally, we define the Hodge duals of the fields, the Hodge functional derivatives and the BV Laplacian in this case by the same formulae as in Section~\ref{sec-Lapl} 
(with the only difference that by $\braket{}{}$ we mean here the duality pairing between $\Lg$ and $\Lg^*$). 
We will denote all these objects by the same symbols as in the previous sections.

\subsubsection{The BV superformalism and the BV action}
We choose a background flat connection $A_0$, and we write a general connection $A$ as $A=A_0+a$, with $a$ in $\Omega^1(M,\ad P)$, with ghost number $0$.
We are now ready to define in this case the analogues 
of the superforms introduced in Section~\ref{sec-BVsuper}, namely
\begin{align*}
\sfB&:=\sum_{k=1}^{m-2} (-1)^{\frac{k(k-1)}{2}} \tau_k +B+(-1)^n a^+ + c^+,\\
\sfA&:=(-1)^{n+1} c+A+(-1)^n B^++\sum_{k=1}^{m-2} (-1)^{\frac{k(k-1)}{2}+n(k+1)} \tau_k^+.
\end{align*}
We also define $\sfa:=\sfA-A_0$.
We notice that $\sfB$ is a superform of total degree $m-2$ with values in $\ad^* P$, 
while $\sfA$ can be interpreted once again as a superconnection on $M$.
It is not difficult to see that the curvature of the superconnection $\sfA$ is given by the formula
\[
\sfF_\sfA=\dd_{A_0} \sfa +\frac{1}{2}\lb{\sfa}{\sfa}.
\]
We go on, as in subsection~\ref{sec-ser}, to define the functional derivatives 
w.r.t.\ $\sfB$ and $\sfa$ and the super BV antibracket; they enjoy the same properties as the previously introduced ones, and we will denote them by the same symbols.

Finally, we claim that the BV action for the canonical $BF$ theory
on $M$ is given by the formula
\begin{equation*}
\SBV:=\int_M \dbraket{\sfB}{\sfF_{\sfA}}.
\end{equation*}
In order to prove the claim, we show once again separately that $\SBV$ satisfies the master equation and that it is (at least formally) $\DBV$-closed.

{}~

\noindent{\em The master equation}.
The proof that $\SBV$ satisfies the master equation is analogous to the proof of the corresponding claim in Section~\ref{sec-mastereq}; we therefore omit it.
We will only write down the action of the super BV differential 
$\sfdelta:=\sbv{\SBV}{\ }$ on the super fields
$\sfa$ and $\sfB$
\begin{align}
\sfdelta \sfa&=(-1)^m \sfF_{\sfA},\label{deltaacan}\\
\sfdelta \sfB&=(-1)^m \left(\dd_{A_0}\sfB+\sad^*(\sfa)\sfB\right);
\label{deltabcan}
\end{align}
the action of the usual BV differential on all the fields (fields {\emph{and}} antifields) is encoded in the two previous equations and, upon switching off the antifields, gives
back the BRST operator defined at the beginning of subsubsection~\ref{sec-BRSTcan}.
It is also easy to verify that $\SBV$ reduces to $S^{\mathrm{can}}$ if
we set the antifields to zero.

{}~

\noindent{\em The $\DBV$-closedness of the BV action}.
The proof that $\SBV$ satisfies the equation
\[
\Delta_{BV} \SBV=0
\]
is a little bit different from the proof of the same identity in Section~\ref{sec-LaplBF}; it relies on a formal argument similar to that used in \cite{CF}.

As noted before, the main property of the BV Laplacian lies in the fact that it contracts each field with the corresponding antifield at the same point
(see Remark~\ref{rem-DBV}); 
therefore, the only terms in the BV action that are not trivially annihilated by the BV Laplacian are of the form
\[
\int_M \braket{\phi_\alpha^+}{\Lie{c}{\phi^{\alpha}}},
\]
for some field $\phi^{\alpha}$.
More precisely, they are (independently of the dimension of $M$) given by the combination
\begin{multline*}
I=
\frac{1}{2}\sqbr{c^*}{\Lie{c}{c}}-\sqbr{a^*}{\Lie{c}{a}}+\\
-\sqbr{B^*}{\ad^*(c)B}
+\sum_{l=1}^{m-2}(-1)^{l+1}\sqbr{\tau_l^*}{\ad^*(c)\tau_l}.
\end{multline*}
This is obtained from the formula for the BV action after rewriting the dot duality, 
the super coadjoint action and the dot Lie bracket in terms the usual ones, 
and recalling that the integral selects only the top form degree part of the integrand.

W.r.t.\ the bases $\{X_i\}$ and $\{X^j\}$, we can write a field $\phi^{\alpha}$ with values in $\Omega^*(M,\ad P)$, resp.\ in $\Omega^*(M,\ad^* P)$, as 
$\phi^{\alpha}=\phi^{\alpha\ i} X_i$, resp.\ $\phi^{\alpha}=\phi^{\alpha}_j X^j$.
For any two real-valued forms on $M$ with the same degree we define
\[
\cirbr\alpha\beta = \int_M \alpha\wedge\star\beta,
\]
where $\star$ is the star Hodge operator w.r.t.\ some chosen metric on $M$.
We therefore obtain 
\begin{multline*}
I=
-\frac{1}{2}f_{jk}^i\cirbr{c^*_i}{c^k\ c^j} - f_{jk}^i\cirbr{a^*_i}{a^k c^j} +\\
+f_{ji}^k\cirbr{(B^{*})^i}{B_k c^j}
+\sum_{l=1}^{m-2}f_{ji}^k \cirbr{(\tau^*_l)^i}{(\tau_{l})_{k} c^j};   
\end{multline*}
we have used here the identity
$\braket{X_i}{\ad^*(X_j) X^k}=-\braket{X^k}{\Lie{X_j}{X_i}}=-f^k_{ji}$.

Finally, we apply the BV Laplacian to the above expression and get (independently of the dimension of $M$)
\begin{multline*}
\DBV\SBV=\DBV I =\\
=C \int_M \dd \mathit{vol}\ 
f_{ji}^i c^j \left[\binom{m}{m}-\binom{m}{m-1}+\binom{m}{m-2}-
\dots+(-1)^l \binom{m}{m-l}+\dotsb\right]=\\
=C \int_M \dd\mathit{vol}\  f_{ji}^i c^j (1-1)^m=0,
\end{multline*}
where $\dd\mathit{vol}$ is the Riemannian volume element and
$C$ is an infinite constant (explicitly, a Dirac distribution evaluated in $0$). 
The binomial coefficients appear as the number of components of the forms 
$\phi^{\alpha}_j$; e.g.\ , $B_k$ is an $m-2$ form on the $m$-dimensional manifold $M$, so
it has $\binom{m}{m-2}$ components in local coordinates. The signs before the binomial coefficients are determined by the ghost numbers of the fields $\phi^{\alpha}$ (recall the explicit definition \eqref{defBVlapl} of the BV Laplacian).

Of course the previous computation should be performed with a regularization
in order to avoid the infinite constant $C$. If the regularization is such that
the above formal manipulations still hold, then $\SBV$ is BV harmonic.

\subsection{Sigma-model interpretation}\label{sec-sm}
The BV action \eqref{SBF} can be seen as the action functional
for a supersymmetric sigma model
with source $\Sigma:=\Pi TM$, where $M$ is our original $m$-dimensional 
manifold and $\Pi$ indicates that the fiber has to be taken with reversed 
Grassmann parity; the target $N$ has to be chosen among 
the following possibilities:
\[
\begin{tabular}{c|c|c}
 & ordinary $BF$ & canonical $BF$\\ \hline
$m$ odd & $\Pi\frg \times\Pi\frg$  & $\Pi\frg\times\Pi\frg^*$  \\ \hline
$m$ even & $\Pi\frg\times\frg$  & $\Pi\frg\times\frg^*$
\end{tabular}
\]
where $\Pi$ again reverses the Grassmann parity. To encompass all cases, we will
write $N=V_1\times V_2$ with $V_1$ and $V_2$ as in the above table.
The superfields $\sfa$ and $\sfB$ are then {\em related to}\/
the $1$ and $2$ components
of a map $f\colon\Sigma\to N$.
Also recall that there is a pairing $\braket{\ }{\ }$ of $V_2$
with $V_1$ which is induced from the bilinear form of Assumption~\ref{Ass3}
resp.\ from the canonical pairing in the case of ordinary resp.\ canonical
$BF$ theories.

Following \cite{AKSZ} one can give the BV bracket and the BV action for
$BF$ theories (to begin with in the case when
the background connection $A_0$ is trivial)
a beautiful interpretation in terms of a natural $QP$-structure
on the space $\calE$ of maps $\Sigma\to N$.
Let us recall that a $P$-manifold is a supermanifold endowed with an odd
non-degenerate closed $2$-form (shortly, an odd symplectic form);
a $Q$-manifold is a supermanifold endowed with an odd vector field $Q$ that
has vanishing Lie bracket with itself; finally, a $QP$-manifold is a 
supermanifold that has both structures in a compatible way, i.e., 
such that the odd symplectic form is
$Q$-invariant. Notice that an odd symplectic form defines a BV bracket;
moreover, an even solution of the master equation defines a compatible
$Q$-structure and vice versa.

The $P$-structure on $\calE$
is defined in terms of the
following constant symplectic form on $N$:
\begin{gather*}
\omega(v_1\oplus v_2,w_1\oplus w_2):=
\braket{v_2}{w_1}-(-1)^m\braket{v_1}{w_2},\\
v_1\oplus v_2,\ w_1\oplus w_2\in T_{(\xi_1,\xi_2)}N\simeq
V_1\oplus V_2,\qquad \forall(\xi_1,\xi_2)\in N.
\end{gather*}
Observe that $\omega$ is odd (even) when $m$ is even (odd); 
i.e., $\omega$ defines an ordinary symplectic structure---though on an 
odd vector space---when $m$ is odd and a $P$-structure when $m$ is even.
This induces the following constant {\em odd}\/ symplectic form on $\calE$:
\begin{gather*}
\tilde\omega(\phi,\phi'):=\int_\Sigma \omega(\phi,\phi')\;\dd\mu,\\
\phi,\ \phi'\in T_f\calE\simeq\Gamma(\Sigma, f^*TN),
\qquad \forall f\in\calE.
\end{gather*}
Here we have denoted by $\int_{\Sigma} \dd\mu$ the canonical density associated 
to the supermanifold $\Pi TM$. It is defined as follows: since $\Pi TM=(M,\Omega^*(M))$, 
then every function on $\Pi TM$ can be identified with a sum of forms on $M$ of all degrees, so there is a canonical density given by the usual integral of forms (which selects the top degree part of the integrand). Locally, 
\[
\int_{\Sigma\arrowvert_U} \dd\mu=\int_U \dd x^1\dotsm\dd x^m\ ,
\]
where the $x$'s are local coordinates on $M$. 
%and the $\theta$s are the
%induced odd local coordinates on the fiber. Observe that $\dd\mu$ is
%odd (even) when $m$ is odd (even). So $\tilde\omega$ is odd in any case.
%The BV bracket defined by this $P$-structure on $\calE$ is the one
%considered in the previous sections.

Next we come to the $Q$-structure.
Observe first that any flow on $\Sigma$
or on $N$ defines (by composition on the right resp.\ on the left) a flow
on $\calE$ and that flows of the two types commute. Infinitesimally, any 
vector field on $\Sigma$ or on $N$ defines a vector field on $\calE$ and
vector fields of the two types commute. Moreover, nilpotency is preserved.
In conclusion, any $Q$-structures $Q_\Sigma$ on $\Sigma$ and $Q_N$ on
$N$ determine $Q$-structures $\widetilde{Q}_\Sigma$ and $\widetilde{Q}_N$ on
$\calE$; moreover, any linear combination of the two is still a $Q$-structure.
On $N$ we consider the $Q$-structure given by the Hamiltonian vector field 
associated to the function
\[
\sigma(\xi_1,\xi_2)=\frac12 \braket{\xi_2}{\Lie{\xi_1}{\xi_1}}.
\]
Observe that this function is odd (even) for $m$ odd (even), so the
corresponding vector field is always odd. The corresponding $Q$-structure
on $\calE$ yields the following action on the superfields
$\sfa$ and $\sfB$:
\[
\sfdelta_N\sfa=\frac12 \lb\sfa\sfa,
\qquad \sfdelta_N\sfB =
\begin{cases}
\lb\sfa\sfB & \text{ordinary $BF$,}\\
\sad^*(\sfa)\sfB & \text{canonical $BF$.}
\end{cases}
\]
This $Q$-structure is automatically compatible with the $P$-structure
defined above. 
On $\Sigma$ we consider instead the canonical $Q$-structure which in local
coordinates reads
\[
Q_\Sigma=\theta^i\,\frac\partial{\partial x^i}.
\]
The induced $Q$-structure on $\calE$ acts on the superfields by
\[
\sfdelta_\Sigma\sfa=\dd\sfa,
\qquad\sfdelta_\Sigma\sfB=\dd\sfB.
\]
This $Q$-structure is also compatible with the $P$-structure defined by
$\tilde\omega$ as follows by an explicit computation: in fact, it is not difficult to check that the odd vector field $Q_{\Sigma}$ has zero-divergence w.r.t.\ the density specified above. 
Since $M$ has no boundary, this guarantees automatically that the $P$-structure on $\mathcal{E}$ is compatible with the $Q$-structure defined by $Q_\Sigma$.

Finally, one considers a linear combination with nonvanishing coefficients
of the above vector fields. This yields an entire family of $QP$-structures
on $\calE$. Notice however that rescaling $\sfa$ with a parameter $\lambda$
and $\sfB$ with $1/\lambda$ ($\lambda\not=0$)
is a canonical transformation. So, up to equivalence, one can always
set the coefficients to have the same ratio as in \eqref{deltaa} and
\eqref{deltab} (or \eqref{deltaacan} and \eqref{deltabcan}
for canonical $BF$ theories). Given the $P$-structure, there is
a unique (up to an additive constant) action functional generating
the given $Q$-structure. Choosing the additive constant appropriately,
the action functional is then a multiple of our $S$ in \eqref{SBF}. 
Finally, the 
remaining multiplicative constant can be absorbed in $\hslash$ (or taken
as a definition thereof).

In order to take into account nontrivial background connections (or
even nontrivial bundles $P\to M$), one has to modify a little bit the above
construction. First one has to introduce a vector bundle $E\to\Sigma$ with
fiber $N$, with $\Sigma$ and $N$ as above. If the original bundle $P$ is 
trivial, so will be $E$ (otherwise it will be constructed by using the
transition functions of $\ad P$ and $\ad^* P$). The space $\calE$ will be
now the space of sections of $E$. The $P$- and $Q_N$-structures are introduced
as above. The $Q_\Sigma$-structure instead requires the choice of a 
connection $A_0$
in order to lift to $E$ the vector field on $\Sigma$; this connection
has moreover to be flat to ensure the nilpotency of $\tilde Q_\Sigma$.
The rest of the construction is the same as above.

\subsection{Gauge fixing}\label{gf}
We conclude this section giving a brief account on the gauge fixing necessary
to start a perturbative expansion of the theory. (For simplicity we restrict
ourselves to ordinary $BF$ theories, the modifications needed for the
canonical ones being obvious.)

The first step is to extend the space of fields by introducing antighosts
and Lagrange multipliers. Along with the usual ghost $c$ one introduces
an antighost $\bar c$ (of ghost number $-1$) and a Lagrange multiplier
$\lambda$ (of ghost number $0$); both are chosen to take values in
the sections of $\ad P$.
Similarly, along with the ghost $\tau_1$ one introduces an antighost
$\bar\tau_1$ and a Lagrange multiplier $\lambda_1$ with values in
$\Omega^{m-3}(M,\ad P)$. As for the ghosts-for-ghosts $\tau_k$, one needs an
entire family of $k$ antighosts and $k$ Lagrange multipliers (\cite{BlT}).
Namely, we denote by $\bar\tau_{k,i}$ and $\lambda_{k,i}$
($i=1,\dots,k$) the antighosts and the Lagrange multipliers
corresponding to $\tau_k$, 
all of which take values in $\Omega^{m-2-k}(M,\ad P)$.
As for the ghost number, one sets
\[
\gh(\bar\tau_{k,i})=2i-k-2,\qquad
\gh(\lambda_{k,i})=2i-k-1.
\]
We will denote by $\Phi$ the collection of the fields including the newly
introduced ones.

Next one has to consider antifields for the antighosts and the Lagrange 
multipliers. They will be denoted by $\bar c^+$, $\lambda^+$,
$\bar\tau_1^+$, $\bar\lambda_1^+$, $\bar\tau_{k,i}^+$ and 
$\lambda_{k,i}^+$ ($k=2,\dots,m-3$; $i=1,\dots,k$)
with the usual rule; i.e., each antifield takes values
in the space of forms of complementary degree of
the corresponding field and its ghost number is minus
the ghost number of the corresponding field, minus one. We will denote
by $\Phi^+$ the collection of all the antifields including the new ones.
Finally, we extend the BV structure by declaring that each of the new antifield
is BV-canonically conjugated to its field.

The newly introduced fields are there only to write down a gauge fixing fermion
(see later). From the point of view of BV cohomology their complex must be
trivial; i.e., one sets
\[
\delta\bar\tau_{k,i}=\lambda_{k,i},\quad
\delta\lambda_{k,i}=0,\qquad
k=2,\dots,m-3;\ i=1,\dots,k.
\]
This is achieved by defining the extended BV action:
\[
S^{\mathrm{ext}}=S + \Sigma,
\]
with $S$ given in \eqref{SBF} and
\begin{equation}\label{Sigma}
\Sigma=\int_M\left(
-\bar c^+\lambda - \bar \tau_1^+\lambda_1
+\sum_{k=2}^{m-3}\sum_{i=1}^k (-1)^k \bar\tau_{k,i}^+\lambda_{k,i}
\right).
\end{equation}
The gauge-fixed action, which is a function of $\Phi$ only,
is then defined by
\[
S^{\mathrm{g.f.}}=S^{\mathrm{ext}}_{\big|_{{\Phi^+=
\frac{\ld \Psi}{\partial\Phi}}}},
\]
where $\Psi$ (the {\em gauge-fixing fermion})
is a function of $\Phi$ of ghost number $-1$ and has to
be chosen so that the Hessian of $S^{\mathrm{g.f.}}$ at a critical
point is not degenerate.
In case one wants to expand around a given flat connection $A_0$,
a suitable gauge-fixing fermion (in accordance with Assumption~\ref{Ass1})
is
\begin{equation}\label{gff}
\Psi=\int_M \bar c\; \dd_{A_0}\star a+
\bar\tau_1\;\dd_{A_0}\star B +
\sum_{k=1}^{m-4} \bar\tau_{k+1,1}\;\dd_{A_0}\star\tau_k +
\sum_{k=1}^{m-4}\sum_{i=2}^{k+1} \bar\tau_{k+1,i}\;\dd_{A_0}\star
\bar\tau_{k,k+2-i},
\end{equation}
where $\star$ is the Hodge star operator induced from a Riemannian
metric on $M$. The BV formalism ensures in particular
that the partition function
and the expectation values of BV closed observables do not depend
on the chosen metric.

\subsection{Superpropagator}
The perturbative expansion of the partition function or of the expectation
values of observables is obtained in terms of propagators, i.e., expectation
values of the fields w.r.t.\ the quadratic part of the action $S^{\rm ext}$.
We will briefly describe this computation in the case of ordinary $BF$ 
theories.

Since the interaction terms and the observables that we will introduce in the
next sections depend only on the superfields $\sfa$ and $\sfB$, it is
enough to compute the ``superpropagator''
\[
\mathrm{i}\hslash\,\eta=
\vevo{\pi_1^*\sfa\ \pi_2^*\sfB}:=\frac1{Z}
\int_{\Phi^+=\frac{\ld \Psi}{\partial\Phi}}
\exp\left(\int_M \dbraket\sfB{\dd_{A_0}\sfa}+\Sigma\right)
\ \pi_1^*\sfa\ \pi_2^*\sfB,
\]
where $Z$ is the partition function, $A_0$ is the chosen background
flat connection, $\Sigma$ is the extension defined in \eqref{Sigma},
and $\pi_1$ and $\pi_2$ are the projections
from $M\times M$ to $M$. So $\eta$ is a distributional $(m-1)$-form
on $M\times M$ with values in $\ad P\boxtimes\ad P$.
This superpropagator with the gauge fixing
\eqref{gff} can be computed by generalizing Axelrod and Singer's
construction \cite{AS} to higher dimensions.
Another possibility is to compute the main properties of the superpropagator
and then construct a form that satisfies them generalizing
the construction of \cite{BC2}. The first property relies on the symmetry
$\sfa\leftrightarrow\sfB$ of the quadratic part of the action:
$\int_M \dbraket\sfB{\dd_{A_0}\sfa}$. This implies
\begin{equation}\label{Teta}
T^*\eta = (-1)^m\,\eta
\end{equation}
where $T$ is the automorphism of $\ad P\boxtimes\ad P$ that acts
on the basis by exchanging the points and at the same time exchanges
the corresponding fibers (in a local trivialization 
$T(x,x';\xi,\xi')=(x',x;\xi',\xi)$, with $x,x'\in U\subset M$ 
and $\xi,\xi'\in\Lg$).
A subsequent computation shows that
\[
\mathrm{i}\hslash\,\dd_{A_0}\eta=(-1)^m\,\vevo{\sfdelta_0(\pi_1^*\sfa\ \pi_2^*\sfB)},
\]
where $\sfdelta_0$ is the linear part of $\sfdelta$. By the main properties
of the BV formalism, one then gets the Ward identity
\[
(-1)^m\,\dd_{A_0}\eta=\vevo{\sfDelta(\pi_1^*\sfa\ \pi_2^*\sfB)}.
\]
The right-hand side is a delta form localized on the diagonal 
$\mathit{Diag}$
of $M\times M$ tensorized with the section $\phi$ of $\ad P\otimes \ad P
\to \mathit{Diag}$ determined by the invariant form $\braket{\ }{\ }$;
that is, $\phi$ is the section corresponding to the constant
equivariant map $\tilde\phi\colon P\to\Lg\times\Lg$, 
$p\mapsto\sum_i \sigma_i\,e_i\otimes e_i$, where $\{e_i\}$ is 
a pseudo-orthonormal basis of $\frg$:
$\braket{e_i}{e_j}=\sigma_i\delta_{ij}$, $\sigma_i=\pm1$.
Thus, if we define $\braket{\ }{\ }_{13}$ on $\Lg\otimes\Lg\otimes\Lg$
as acting on the first and third components
and define consequently $\dbraket{\ }{\ }_{13}$,
we may interpret $\eta$ as a distributional
form such that the 
linear operator $\sP\colon\Omega^*(M,\ad P)\to\Omega^{*-1}(M,\ad P)$,
\[
\sP\gamma := \pi_{2*}\dbraket\eta{\pi_1^*\gamma}_{13},
\qquad\gamma\in\Omega^*(M,\ad P),
\]
satisfies the equation
\begin{equation}\label{dP}
\dd_{A_0}\sP + \sP\dd_{A_0} = \mathit{id}.
\end{equation}

A regularized version of $\eta$ consists in finding a smooth
form (which we will continue to denote by $\eta$) on the configuration space
$C_2(M):=M\times M\setminus\mathit{Diag}$
so that $\sP$ defined as above (with the obvious understanding
that the projections $\pi_1$ and $\pi_2$ are now from $C_2(M)$
to $M$)
satisfies the same equation.
Notice however that $C_2(M)$ is not compact; so one has to replace it
by its compactification 
$\overline C_2(M)$ which is obtained from
$M\times M$ by replacing the diagonal with its differential-geometric
blowup. Observe that $\overline C_2(M)$ is a manifold with boundary
the spherical normal bundle $SN\mathit{Diag}$
of $\mathit{Diag}$ in $M\times M$.
Since we have removed the diagonal, we require now that the
superpropagator should be an $A_0$-covariantly closed form 
$\eta\in\Omega^{m-1}(\overline C_2(M), \ad P\boxtimes\ad P)$, where, by abuse of notation, we have denoted by $\ad P\boxtimes \ad P$ the pulled-back bundle of $\ad P\boxtimes \ad P$ w.r.t.\ the projection $\\overline{C}_2(M)\to M\times M$.
Observe that the generalized Stokes formula implies that
the left-hand side of \eqref{dP} applied to a form $\gamma$
is $\pi^\partial_*\dbraket{\iota ^*\eta}{\pi^{\partial*}\gamma}_{13}$,
with $\iota$ the inclusion of $SN\mathit{Diag}$ in $\overline C_2(M)$
and $\pi^\partial$ the projection $SN\mathit{Diag}\to\mathit{Diag}$.
Thus, for \eqref{dP} to hold, one has to require that 
the restriction of $\eta$ to the boundary should be
\begin{equation}\label{iotaeta}
\iota^*\eta=
\vartheta\otimes\pi^{\partial*}\phi + \pi^{\partial*}\beta
\end{equation}
where $\vartheta$
is the global angular form of $SN\mathit{Diag}$ and 
$\beta\in\Omega^{m-1}(\mathit{Diag},\ad P\otimes\ad P)$
is so far undetermined. Recall that a global angular form
$\vartheta$ on a sphere bundle $S\xrightarrow{\pi^\partial}B$
is a form satisfying $\pi^\partial_*\vartheta=1$ and 
$\dd\vartheta=-\pi^{\partial*}e$, where $e$ is a representative
of the Euler class of the bundle.
In our case, since $N\mathit{Diag}\simeq TM$, $e$ will be a representative
of the Euler class of $M$. The first property of $\vartheta$ is what we need
for \eqref{dP} to hold; the second property, of which one cannot dispose,
implies
$\dd_{A_0}\beta=e\otimes\phi$. 

This is a very strong constraint in
even dimensions, as it requires the Euler class of $M$ to be trivial.
In fact, multiply both sides of the equation by $\phi$ and contract
the first $\Lg$-component with the third and the second with the fourth;
this yields
\[
\dd\dbraket\phi\beta_{13,24}=e\,\dim\Lg.
\]
% This in particular means that
%in even dimensions, where $e$ cannot be chosen to vanish, a nonvanishing
%non--covariantly closed
%$\beta$ is necessary.
Notice finally that we also want $\eta$ to satisfy \eqref{Teta},
with $T$ now the corresponding involution on 
$\ad P\boxtimes\ad P\to\overline C_2(M)$. In particular, this implies that
one has to choose $\vartheta$ to be even (odd) under the antipodal
map on the fibers if $m$ is even (odd); in odd dimensions this also
implies that one must choose $e=0$.
Moreover, $\beta$ has to be an element of 
$\Omega^{m-1}(\mathit{Diag},\mathfrak{S}^2\ad P)$ in even dimensions and
of $\Omega^{m-1}(\mathit{Diag},\bigwedge^2\ad P)$ in odd dimensions.

Such a form $\eta$ 
can be obtained generalizing the construction of \cite{BC2}:
\begin{Thm}
Under Assumptions~\ref{Ass1} and \ref{Ass3}, there exists a covariantly
closed  element
$\eta$ of $\Omega^{m-1}(\overline C_2(M), \ad P\boxtimes\ad P)$,
satisfying \eqref{Teta} and \eqref{iotaeta}.
Moreover, in odd dimensions $\beta$ will be automatically
covariantly closed, while in even dimensions---where the above
Assumptions imply $[e]=0$---this will be true only if one chooses $e=0$.
Finally, $\beta$ may be chosen to vanish if 
$H^{m-1}_{\dd_{A_0}}(M,\bigwedge^2\ad P)$ is trivial in odd dimensions
and if $H^{m-1}_{\dd_{A_0}}(M,\mathfrak{S}^2\ad P)$ is trivial in 
even dimensions.
\end{Thm}
\begin{proof}
One first builds
a global angular form $\vartheta$ on $SN\mathit{Diag}$ 
with the correct behavior under the antipodal map on the fibers:
one may construct it as in Appendix~\ref{app-gaf} using 
the Levi-Civita connection for a given Riemannian metric,
which also allows to identify $SN\mathit{Diag}$ with the unit
sphere bundle $SO\mathit{Diag}\times_{SO(m)}S^{m-1}$.
Next one extends $\vartheta$ to
the complement of the zero section of $N\mathit{Diag}$ and multiplies
it by a function $\rho$ that is identically one in a neighborhood $U_1$ of
the zero section and identically zero outside a second neighborhood 
$U_2\supset U_1$. One then defines 
$\eta_0\in\Omega^{m-1}(\overline C_2(M), \ad P\boxtimes\ad P)$ 
as the extension by zero of $\rho\,\vartheta\otimes\pi^{\partial*}\phi$.
Since $\dd_{A_0}\phi=0$, 
$\dd_{A_0}\eta_0$ is the extension by zero of
$\dd\rho\,\vartheta\otimes \pi^{\partial*}\phi-\rho\,\pi^{\partial*}(e\otimes\phi)$. 
The last form may be extended to the zero section of $N\mathit{Diag}$; hence, the extension by zero of $\dd_{A_0} \eta_0$ can be seen as a covariantly closed element of $\Omega^m(M\times M,\ad P\boxtimes \ad P)$. 
The general K\"unneth theorem implies $H^*_{\dd_{A_0}}(M\times M, \ad P\boxtimes\ad P)\cong
H^*_{\dd_{A_0}}(M,\ad P)^{\otimes2}$. So Assumption~\ref{Ass1} implies
that there is a form $\alpha\in\Omega^{m-1}(M\times M, \ad P\boxtimes\ad P)$
such that $\dd_{A_0}\pi^*\alpha=\dd_{A_0}\eta_0$, with $\pi$
the projection $\overline C_2(M)\to M\times M$. Also observe that
one may choose $\alpha$ to satisfy $T^*\alpha=(-1)^m\,\alpha$.
Finally, define $\eta:=\eta_0-\pi^*\alpha$. An easy check shows
that it satisfies all properties above
(with $\beta$ determined by the restriction
of $\alpha$ to the diagonal).
\end{proof}

\begin{Rem}
There are a couple of interesting cases when $M$ does not satisfy
Assumption~\ref{Ass1}, but one can define the superpropagator anyway.
First, when $M=\bbR^m$ (see subsection~\ref{ssec-Rm}) all boils down
to looking for (the higher-dimensional generalization of) 
Bott and Taubes's \cite{BT} tautological forms, as described in \cite{CR}.
Second, when $M$ is a rational homology sphere,
one can generalize the construction of \cite{BC}
(which does not yield a closed $\eta$, so that extra diagrams
must be introduced to correct for it)
or alternatively remove one point, as suggested in \cite{K}, and essentially
reduce to the previous case.
\end{Rem}

\section{Generalized Wilson loops in odd dimensions}\label{sec-gWl}
In this section we display some observables for 
odd-dimensional $BF$ theories which in some sense generalize the 
classical observables \eqref{Taylorholon}, i.e.\ the iterated-integral
expansions of Wilson loops.
In the first subsection we construct a flat invariant observable (see 
Definition~\fullref{def-flobs})
$\sfS_3$ 
which represents a sort of ``cosmological term''  (although it does not have the correct ghost number, except for the case $\dim M=3$).
We next define in subsection~\ref{ssec-gWlBV}
a ``generalized holonomy'' constructed via iterated integrals 
by means of $\sfA$ and $\sfB$, and we show that it defines a cohomology
class w.r.t.\  the super BV coboundary operator twisted with $\sfS_3$
which takes values in $H^*(\ImbfM)$. From this we then derive a true BV 
observable.

\subsection{The ``cosmological term''}
We define the local functional 
\[
\sfS_3:=\frac{1}{6} \integ{\sfB}{\lb{\sfB}{\sfB}}
\]
which is an element of $\pser(\mathbb{R})$ of total degree $2m-6$.
We want to show that $\sfS_3$ is a flat invariant observable in the sense of
definition~\ref{def-flobs}. This is expressed by the following
\begin{Lem}
\begin{align}
\boldsymbol{\delta} \sfS_3 &=0,\label{dO}\\
\boldsymbol{\Delta} \sfS_3 &=0,\label{DO}\\
\sbv{\sfS_3}{\sfS_3}&=0.\label{OO}
\end{align} 
\end{Lem}

\begin{proof}
First of all, we write down the left partial derivatives of $\sfS_3$: 
\[
\lA{\sfS_3}=0,\quad \lB{\sfS_3}=\frac{1}{2}\lb{\sfB}{\sfB}.
\]
With the help of~\eqref{dSdA} and by the definition of the super BV antibracket, we get
\[
\boldsymbol{\delta} \sfS_3=\sbv{S}{\sfS_3}=  \frac{1}{2} \integ{\dd_{\sfA}\sfB}{\lb{\sfB}{\sfB}};
\]
By the invariance of $\braket{\ }{\ }$ it follows
\[
\frac12 \integ{\dd_{\sfA}\sfB}{\lb{\sfB}{\sfB}}=\frac16 \int_M \dd \dbraket{\sfB}{\lb{\sfB}\sfB}=0
\]
by Stokes' theorem.
So we have proved \eqref{dO}.

Eqns.\ \eqref{DO} and \eqref{OO} follow from the definitions of the super 
BV antibracket and of the super 
BV Laplacian $\sDBV$ and 
from the fact that $\sfS_3$ depends only on $\sfB$.
\end{proof}

It follows from Lemma~\ref{lem-flobs} 
that not only $\sfS_3$ but any of its multiples
is a flat observable. So we introduce the ``cosmological constant''
$\kappa$ and consider a twisting by $\kappa^2\sfS_3$ (the reason for putting
$\kappa^2$ instead of $\kappa$ will be clear in the next subsection). 
We then define
\[
\sfdelta_{\kappa^2}:= \sfdelta + \kappa^2\sbv{\sfS_3}{\ }.
\]
and, again by Lemma~\ref{lem-flobs}, $\sfdelta_{\kappa^2}$ is an odd differential
for any $\kappa$.
Its action  on the fundamental superfields is easily
computed:
\begin{equation}\label{derivo}
\sfdelta_{\kappa^2}\sfa = -\sfF_\sfA
-\frac{\kappa^2}2\lb\sfB\sfB,
\qquad \sfdelta_{\kappa^2}\sfB=-\dd_{\sfA}\sfB.
\end{equation}

\subsection{The generalized Wilson loop in the BV superformalism}\label{ssec-gWlBV}
We want to define an object that generalizes the observable introduced 
in~\cite{CCFM} for the $3$-dimensional $BF$ theory with cosmological term.
We shall realize this proposal by introducing the new superform 
\[
\sfC_{\kappa}:=\sfa+\kappa \sfB.
\]
Observe that $\sfC_{\kappa}$ is not a homogeneous element in $\pser(M,\ad P)$
w.r.t.\ the total degree, but it is homogeneous of degree one 
with respect to its reduction modulo $2$. 
By recalling (\ref{derivo}), it is easy to see that
\[
\boldsymbol{\delta}_{{\kappa^2}}\sfC_{\kappa}=-\dd_{A_0} \sfC_{\kappa}-\frac{1}{2}\lb{\sfC_{\kappa}}{\sfC_{\kappa}}.
\]
The previous equation suggests that we may interpret the superform $\sfC_{\kappa}$ as a ``variation'' of the flat connection $A_0$, and therefore $\boldsymbol{\delta}_{{\kappa^2}}\sfC_{\kappa}$ can be interpreted as its curvature.
Observe that, since $\sfC_\kappa$ is of odd degree, all the formulae
of Appendix~\ref{app-hol} are basically the same as if $\sfC_\kappa$ were
an ordinary variation of $A_0$. We exploit then this analogy to define
the $n$-th iterated integral of $\sfC_{\kappa}$ as
\begin{align*}
\pi_{n*}\left(\widehat{\sfC_{\kappa}}_{1,n}\cdots \widehat{\sfC_{\kappa}}_{n,n}\right)\holo{A_0}{0}{1}.
\end{align*}
We refer from now on to Appendix~\ref{app-hol} for the main notations (simplices, evaluation maps, etc.). 
We recall the definition of $\widehat{\sfC_{\kappa}}$:
We have written 
\[
\widehat{\sfC_{\kappa}}:=\holo{A}0\bullet\ev_1^* \sfC_{\kappa}\left(\holo{A}0\bullet\right)^{-1},
\]
and $\widehat{\sfC_{\kappa}}_{i,n}:=\pi_{i,n}^*\widehat{\sfC_{\kappa}}$.
We have suppressed $\rho$ before all the $\widehat{\sfC_{\kappa}}$'s in the above product; the forms considered in the $n$-th iterated integral take values in the associative algebra $\en V$.
We then define the generalized holonomy of $\sfC_{\kappa}$ from $0$ to $1$ via the path-ordered exponential
\[
\sfholo{\sfC_{\kappa}}:=\holo{A_0}{0}{1}+\sum_{n\geq 1}\pi_{n*}\Big(\widehat{\sfC_{\kappa}}_{1,n} \cdots \widehat{\sfC_{\kappa}}_{n,n}\Big)\holo{A_0}{0}{1};
\]
it defines an element in $\ser(\lo M, \en V)$, and since $\dim \triangle_{n} =n$, it follows that it has even total degree. 
We now pick a finite-dimensional representation $\rho$ and
define the generalized Wilson loop
\begin{equation}
\boldsymbol{\mathcal{H}}_{\rho}(\kappa;\sfA,\sfB)=\tr_{\rho} \sfholo{\sfC_{\kappa}}.
\end{equation}
{}From the previous considerations, it is an element of $\ser (\lo M,\mathbb{R})$, with even total degree.
We are now ready to state the main Theorem of this Section.
\begin{Thm}\label{gen-hol3}
The generalized Wilson loop is 
$(\boldsymbol{\delta}_{{\kappa^2}}+\dd)$-closed:
\begin{equation*}
(\boldsymbol{\delta}_{{\kappa^2}}+\dd) \boldsymbol{\mathcal{H}}_{\rho}(\kappa;\sfA,\sfB)=0.
\end{equation*}
\end{Thm}

\begin{proof}
By above reasonings, we can consider $\sfC_{\kappa}$ as a variation of the (flat) connection $A_0$. 
The cyclicity of the trace allows to replace the exterior derivative by the covariant derivative $\dd_{\ev(0)^*A_0}$. $\sfholo{\sfC_{\kappa}}$ has the same form as $\holo{A+a}01$ of Appendix~\ref{app-hol}, where we have set $A_0=A$, and we have replaced $a$ by $\sfC_\kappa$ and the wedge product by the dot product.
Accordingly to the sign rules for the dot product and repeating almost verbatim the arguments used in the proof of Theorem~\ref{holonomy1}, we get 
\[
\dd \boldsymbol{\mathcal{H}}_{\rho}(\kappa;\sfA,\sfB)=\sum_{m\ge 1} \sum_{i=1}^m (-1)^{m+i} \tr_{\rho} \left\{ \pi_{m*}\left[\widehat{\sfC_{\kappa}}_{1,m} \cdots \left(\widehat{\boldsymbol{\delta}_{\kappa^2}\sfC_\kappa}\right)_{i,m} \cdots \widehat{\sfC_{\kappa}}_{m,m}\right]\holo{A_0}{0}{1} \right\}.
\] 
Recalling Lemma~\ref{push2},~\ref{pull2} and~\ref{trace2} and the Leibnitz rule, it is  then not difficult to verify that
\[
\boldsymbol{\delta}_{{\kappa^2}}\boldsymbol{\mathcal{H}}_{\rho}(\kappa;\sfA,\sfB)=\sum_{m\ge 1}\sum_{i=1}^m (-1)^{m+i+1}\tr_{\rho}\left\{\pi_{m*}\left[\widehat{\sfC_{\kappa}}_{1,m}\cdots \widehat{\left(\boldsymbol{\delta}_{{\kappa^2}}\sfC_{\kappa}\right)}_{i,m}\cdots \widehat{\sfC_{\kappa}}_{m,m}\right]\holo{A_0}{0}{1}\right\} 
\]
which yields the desired identity.
\end{proof}

%Since $\sfS_3$ has even total degree, it follows that $\exp \frac{\ii}{\hslash}\kappa^2 \sfS_3$ has even total degree.
%Let us suppose that we have another functional in $\mathcal{S}_{\sfA,\sfB}(N,\mathbb{R})$, say $G$, for which holds $\boldsymbol{\Delta} G=0$.
%Then we obtain  
%\begin{Lem}\label{q-oss3}
%Under the above assumption on $G$, the product $(\exp\frac{\ii}{\hslash}\kappa^2 \sfS_3)\ G$ satisfies the identity
%\[
%(-\ii \hslash \boldsymbol{\Delta}+\boldsymbol{\delta}+\dd) \big[\big(\exp \frac{\ii}{\hslash}\kappa^2\sfS_3\big)\ G\big]=0
%\]
%if and only if $G$ satisfies 
%\[
%(\boldsymbol{\delta}_{{\kappa^2}}+\dd)G=0\quad,
%\]
%with the above introduced notations.
%\end{Lem}

%\begin{proof}
%We shall simply apply the main properties of super BV Laplacian $\boldsymbol{\Delta}$ and the Leibnitz rule for the BV--operator $\delta$, along with (\ref{dO}) and (\ref{DO}) and the fact that $\boldsymbol{\Delta} G=0$.
%\end{proof}

We would like a stronger assertion than what we proved in the above Theorem;
namely, that ${\mathcal{H}}_{\rho}$ is 
$(-\ii\hslash\sfDelta + \sfdelta_{\kappa^2} + \dd)$-closed. So
we need
\[
\boldsymbol{\Delta}\boldsymbol{\mathcal{H}}_{\rho}(\kappa;\sfA,\sfB)=0.
\]
If a loop has transversal self-intersections, the above identity is certainly
false since on the two intersecting strands appear complementary components
of a field and its antifield. If the loop has non-transversal intersections
or cusps, it is not even clear what the action of the BV Laplacian should be.
However, even restricting to imbeddings might not be enough since
in the computation of the BV Laplacian
there are ill-defined terms coming from subsequent fields in the iterated
integrals as the evaluation points come together.
To establish the validity of the above identity, we can choose the following 
\begin{Reg*}\label{conv-reg}
We only consider elements of
$\ImbfM$, the space of framed imbeddings of $S^1$ into $M$. 
For each element
we then consider a tubular
neighborhood of the imbedding and use the framing to select a companion
imbedding on the boundary.
Finally we put each component of $\sfA$ appearing in the iterated integrals
on the imbedding and  each component of $\sfB$ on its companion
(following a procedure introduced in \cite{CCM}). 
\end{Reg*}
Since the cosmological term is a flat invariant observable we then
obtain, under the above assumption, the following
\begin{Cor}
\[
(\sOBV+\dd)\left[\exp\left(\frac{\ii}{\hslash}\kappa^2 \sfS_3\right) \boldsymbol{\mathcal{H}}_{\rho}(\kappa;\sfA,\sfB)\right]=0. 
\]
\end{Cor}
As a consequence, the $\dd$-cohomology class of the
above functional are BV observables. This implies 
Theorem~2 of \cite{CR}, which states that the above functional
defines an $H^*(\ImbfM)$-valued BV observable.

\begin{Rem}\label{Rem-flat}
We notice finally that the v.e.v.s of the generalized Wilson loops together with the cubic cosmological term do not depend on the representative of
flat connection $A_0$.
Let in fact $g\in\mathcal{G}$ be a gauge transformation viewed as a section
of $\Ad P$. Then one can verify that $\holo{A_0}{t_1}{t_2}$ is sent to 
$g^{-1}(\gamma(t_1))\holo{A_0}{t_1}{t_2}g(\gamma(t_2))$
(see Remark~\ref{rem-hol} for the precise definitions in the general case). 
This implies that the superfields $\sfa$ and $\sfB$ in the generalized Wilson loops are acted upon by $\Ad_g$ (this is a consequence of the definition of the generalized Wilson loops and of the cyclicity of the trace). 
This can be compensated by a change of variables in the functional integral
whose formal measure is constructed using the bilinear form $\braket{\ }{\ }$
and hence formally $\Ad$-invariant.
Therefore, the v.e.v.s of the generalized Wilson loops 
are functions on the moduli space of flat connections.
\end{Rem}

\section{Other loop observables}\label{sec-lobs}
We now generalize the ideas of the previous Section along two 
directions: $i$) consider variations of the connection which
are not necessarily of odd degree; $ii$) introduce interaction terms
with higher powers of $\sfB$. Both generalizations require the following
\begin{Ass}\label{trodd}
Throughout this section we work with a Lie algebra $\frg$, coming from an associative algebra endowed with a trace $\tr$ 
(e.g., we may take $\frg=\mathfrak{gl}(N)$ with the usual trace of matrices).
Furthermore, we define the $\ad$-invariant symmetric bilinear form 
$\braket{\eta}{\xi}$ on $\frg$ by $\tr \eta \xi$ and assume that it is nondegenerate (as required
by Assumption~\ref{Ass3}).
Finally, we will only consider representations $\rho$ of $\frg$ as an associative algebra.
\end{Ass}

\subsection{The odd-dimensional case}
\subsubsection{Higher-order $\sfB$-interactions}
We define, for $k\in \mathbb{N}$, the following even element of $\ser$:
\[
\obs_{2k+1}=\frac{1}{2k+1} \int_M \tr\sfB^{2k+1}.
\]
Observe that even powers of $\sfB$ would vanish by the cyclicity of the 
trace  
\begin{Lem}
The following identities hold for the functional $\obs_{2k+1}$:
\begin{align}
\boldsymbol{\delta} \obs_{2k+1}&=0,\forall k\in \mathbb{N}\label{deltaO} \\
\boldsymbol{\Delta} \obs_{2k+1}&=0,\forall k\in \mathbb{N}\label{DeltaO} \\
\sbv{\obs_{2k+1}}{\obs_{2l+1}}&=0,\quad\forall k,l \in \mathbb{N}.\label{BVOO}
\end{align}
It follows in particular that, $\forall k\in \mathbb{N}$, the functional $\obs_{2k+1}$ is a flat invariant observable (see subsection~\ref{ssec-twist}).
\end{Lem}

\begin{proof}
{}From the definition of the super BV antibracket, we get
\[
\sbv{S}{\obs_{2k+1}}=\int_M \tr\left[\dd_{\sfA}\sfB \cdot \sfB^{2k}\right]
=\frac{1}{2k+1} \int_M \dd \tr \sfB^{2k+1}=0.
\]
(\ref{DeltaO}) and (\ref{BVOO}) follow respectively from the definition of the BV Laplacian and of the super BV antibracket, and from the fact that
the functionals $\obs_{2k+1}$ do not depend on $\sfa$. 
\end{proof}
Let us now choose $n\in\mathbb{N}$ and
a sequence of real numbers
$\boldsymbol{\mu}=\boldsymbol{\mu}(\boldsymbol{\lambda})=
\{\mu_2,\mu_4,\dots,\mu_{4n+2}\}$. Then
we define the following even element of $\pser(\bbr)$:
\[
\obs_{\boldsymbol{\mu}}:=\sum_{k=1}^{2n+1} \mu_{2k} \obs_{2k+1}.
\]
{}From the Lemma it follows that $\obs_{\boldsymbol{\mu}}$ is a flat invariant
observable for any ${\boldsymbol{\mu}}$. So, as in subsection~\ref{ssec-twist},
we can introduce the following odd differential:
\[
\sfdelta_{\boldsymbol{\mu}}=
\sfdelta+\sum_{k=1}^{2n+1}\mu_{2k} \bde_{\sfO_{2k+1}}.
\]
Its action on the fundamental superfields is easily computed:
\begin{equation}\label{derivodispari}
\sfdelta_{\boldsymbol{\mu}}\sfa=
-\sfF_\sfA - \sum_{k=1}^{2n+1}\mu_{2k} \sfB^{2k},\qquad
\sfdelta_{\boldsymbol{\mu}}\sfB=-\dd_{\sfA}\sfB.
\end{equation}

\subsubsection{Extended generalized Wilson loops}
Let $\boldsymbol{\lambda}:=\{\lambda_1,\lambda_3,\dots,\lambda_{2n+1}\}$ be 
another sequence of real numbers with the same $n$ as above.
We then define the odd superform
\[
\sfC_{\boldsymbol{\lambda}}:=\sfa+\sum_{k=0}^{n} \lambda_{2k+1} \sfB^{2k+1}.
\]
If the sequences $\boldsymbol{\mu}$ and $\boldsymbol{\lambda}$ are
related by
\begin{equation}\label{bmulambda}
\mu_{2k}:=\sum_{\substack
{0\le i,j\le n\\
i+j=k-1}}
\lambda_{2i+1}\lambda_{2j+1},
\end{equation}
then \eqref{derivodispari} implies
\[
\boldsymbol{\delta}_{\boldsymbol{\mu}} \supercon =-\dd_{A_0} \supercon -\frac{1}{2} \lb{\supercon}{\supercon}.
\]
The above expression has again the form of a curvature; we can therefore view the superform $\supercon$ as a variation of the connection $A_0$. 
So, analogously to what we did in subsection~\ref{ssec-gWlBV},
we define the path-ordered integral
\[
\sfholo{\supercon}=\holo{A_0}{0}{1}+\sum_{m\ge 1} \pi_{m*}\left(\widehat{\supercon}_{1,m}\cdots\widehat{\supercon}_{m,m}\right)\holo{A_0}{0}{1}.
\]
We next define accordingly 
\[
{\boldsymbol{\mathcal{H}}}_{\rho}(\boldsymbol{\lambda};\sfA,\sfB):=\tr_{\rho}\sfholo{\supercon}.
\]
Repeating the arguments used in the proof of (\ref{gen-hol3}), we can state the following
\begin{Thm}\label{gen-holodd}
If  $\boldsymbol{\mu}$ and $\boldsymbol{\lambda}$ are
related by \eqref{bmulambda}, then
\[
(\boldsymbol{\delta}_{\boldsymbol{\mu}}+\dd)
{\boldsymbol{\mathcal{H}}}_{\rho}(\boldsymbol{\lambda};\sfA,\sfB)
=0.
\]
\end{Thm}
Since $\obsmu$ is a flat invariant observable, this implies the following
\begin{Cor} 
With the same hypothesis as above and with the regularization procedure
defined on page~\pageref{conv-reg}, we obtain
\[
(\sOBV+\dd)\left[\exp\left(\frac\ii\hslash\obsmu\right)
{\boldsymbol{\mathcal{H}}}_{\rho}(\boldsymbol{\lambda};\sfA,\sfB)
\right]=0.
\]
\end{Cor}
Again this implies that the $\dd$-cohomology class of the above functional
is a BV observable; from this
Theorem~4 of \cite{CR} follows.

\begin{Rem}
The same reasonings sketched in Remark~\ref{Rem-flat} do hold in this case as well; therefore, we may conclude that the v.e.v.s of the generalized Wilson loops with higher-order $\sfB$-interactions depend only the class $\left[A_0\right]$ in $\left\{A\in\mathcal{A}\colon F_A=0\right\}/\mathcal{G}$.
\end{Rem}

\subsection{The even-dimensional case}
We now turn to the problem of defining generalized Wilson loop observables
for the case $\dim M$ even. Observe that
in even-dimensional $BF$ theories $\sfB$ has even total degree; so $\lb{\sfB}{\sfB}=0$. This implies that it is not possible to define a generalized cosmological term as in Section~\ref{sec-gWl} because we cannot anymore rely on the dot Lie bracket to construct this functional. 
Therefore, in order to define products of $\sfB$ with itself 
(either cubic or not) we must do as in the preceding subsection and, 
in particular, we need Assumption~\fullref{trodd}.

\subsubsection{$\sfB$-interactions}
For a given $k>1$ we define the following even element of $\pser$:
\[
\sfO_k=
\frac{1}{k} \int_M \tr\sfB^{k}.
\]
We now state the following
\begin{Lem}
The functionals $\sfO_k$ satisfy the identities
\begin{align}
\boldsymbol{\delta} \sfO_k&=0, \forall k>1,\label{paridO}\\
\boldsymbol{\Delta} \sfO_k&=0, \forall k>1,\label{pariDO}\\
\sbv{\sfO_k}{\sfO_l}&=0, \forall k,l>1 \label{pariOO}.
\end{align}
\end{Lem}
\begin{proof}
By definition of the super BV antibracket, we can write
\[
\frac 1k\boldsymbol{\delta}\int_M \tr\sfB^{k} 
=\int_M \tr\left[\dd_{\sfA}\sfB\cdot 
\sfB^{k-1}\right]
= \frac{1}{k}\int_M\dd \tr \sfB^{k}=0.
\]
(\ref{pariDO}) and (\ref{pariOO}) are consequences of the fact that the 
$sfO_k$s 
do not depend on $\sfa$ and of the definitions 
of the super BV antibracket and of the super BV Laplacian.
\end{proof}
Again it follows that each linear combination of $\sfO_k$s is a flat invariant
observable (see subsection~\ref{ssec-twist}). So, for a given
positive integer $n$,
we take a sequence of real numbers $\boldsymbol{\mu}:=\{\mu_2,\mu_3,\dots,\mu_{2n}\}$ and define
\[
\obs_{\boldsymbol{\mu}}=\sum_{i=2}^{2n} {\mu_i}\sfO_{i+1}.
\]
Therefore, Lemma~\ref{lem-flobs} implies that 
\[
\boldsymbol{\delta}_{\boldsymbol{\mu}}:=\boldsymbol{\delta}+\bde_{\boldsymbol{\mu}}:=\boldsymbol{\delta}+
\sum_{i=2}^{2n} {\mu_i}\sbv{\sfO_{i+1}}{\ }
\]
is an odd differential for any sequence $\boldsymbol{\mu}$.
Moreover we have
\begin{equation}\label{evendeltamu}
\boldsymbol{\delta}_{\boldsymbol{\mu}}\sfa=\sfF_\sfA+
\sum_{i=2}^{2n} {\mu_i} \sfB^i,
\qquad
\boldsymbol{\delta}_{\boldsymbol{\mu}}\sfB=
\dd_\sfA\sfB,
\end{equation}
using once again arguments similar to those introduced in the proof of (\ref{deltaa}) and (\ref{deltab}).

\subsubsection{The generalized Wilson loop}
For the same $n$ as above, we consider a sequence
$\boldsymbol{\lambda}=\{\lambda_1,\dots,\lambda_n\}$
so that
\begin{equation}\label{mulambdaeven}
\mu_{i}:=\sum_{\substack
{1\le k,l\le n\\
k+l=i}}
\lambda_{k}\lambda_{l}
\end{equation}
for the previously introduced subsequence $\boldsymbol{\mu}$.
Then we define
\[
\sfB_{\boldsymbol{\lambda}}=\sum_{i=1}^n\lambda_i\sfB^i.
\]
{}From \eqref{evendeltamu} it follows that
\begin{equation}\label{derivdiBlam}
\bde_{\boldsymbol{\mu}} \sfa = \sfB_{\boldsymbol{\lambda}}\cdot \sfB_{\boldsymbol{\lambda}},\qquad \bde_{\boldsymbol{\mu}} \sfB_{\boldsymbol{\lambda}}=0.
\end{equation}
Then, in analogy with the superform $\widehat{\supercon}$ of the previous
subsection, we define
\[
\widehat{\sfB_{\boldsymbol{\lambda}}}:=\holo{A_0}{0}{\bullet}\ \ev_1^*\sfB_{\boldsymbol{\lambda}}\left(\holo{A_0}{0}{\bullet}\right)^{-1},\qquad \widehat{\sfa}:=\holo{A_0}{0}{\bullet}\ \ev_1^*\sfa \left(\holo{A_0}{0}{\bullet}\right)^{-1},
\]
and, accordingly with the notations of Appendix~\ref{app-hol}, $\widehat{\sfB_{\boldsymbol{\lambda}}}_{i,n}$ and $\widehat{\sfa}_{i,n}$, which we will write as $\widehat{\sfB_{\boldsymbol{\lambda}}}_{t_i}$ and $\widehat{\sfa}_{t_i}$.
We then define
\[
{\sfh}_{m,\rho}(\boldsymbol{\lambda};\sfA,\sfB):=\tr_{\rho} \pi_{m*} \left[\holo{\widehat{\sfa}}{0}{t_1}\cdot \widehat{\sfB_{\boldsymbol{\lambda}}}_{t_1} \cdot \holo{\widehat{\sfa}}{t_1}{t_2}\cdots \widehat{\sfB_{\boldsymbol{\lambda}}}_{t_m}\cdot \holo{\widehat{\sfa}}{t_m}{1}\right]\holo{A_0}01,
\]
where we have written
\begin{itemize}
\item $\holo{\widehat{\sfa}}{0}{t_1}:=\pi_{1,m}^*\holo{A_0+\sfa}0\bullet$;
\item $\holo{\widehat{\sfa}}{t_m}1:=\pi_{m,m}^*\holo{A_0+\sfa}\bullet1$;
\item $\holo{\widehat{\sfa}}{t_i}{t_{i+1}}:=\pi_{i,i+1,m}^* \holo{A_0+\sfa}\bullet\bullet$,
\end{itemize}
using the notations of Remark~\ref{def-olovar}, where we have set again $A_0=A$, and we have replaced $a$ by $\sfa$ and wedge products by dot products; $\pi_{i,i,+1,m}(\gamma;t_1,\dots,t_m):=(\gamma;t_i,t_{i+1})$, for $i\in\{1,\dots,m-1\}$.
We finally define
\[
\boldsymbol{\mathcal{H}}_{\rho}^o (\boldsymbol{\lambda};\sfA,\sfB)=
\sum_{m=0}^\infty 
{\sfh}_{2m+1,\rho}(\boldsymbol{\lambda};\sfA,\sfB).
\]
We can now state the main theorem of this subsection
\begin{Thm}\label{pariolonomiagen}
The following identity holds:
\begin{equation*}
(\dd-\boldsymbol{\delta}_{\boldsymbol{\mu}}) 
\boldsymbol{\mathcal{H}}_{\rho}^o (\boldsymbol{\lambda};\sfA,\sfB)
 =0.
\end{equation*}
\end{Thm}

\begin{proof}
We begin by computing the exterior derivative of one of the factors of the above sum.
With the help of the generalized Stokes Theorem we obtain
\begin{equation}
\begin{aligned}
\dd {\sfh}_{2m+1,\rho}(\boldsymbol{\lambda};\sfA,\sfB)&=\tr_\rho (-1)^{2m+1} \pi_{2m+1*}\left\{ \dd_{\pi_{2m+1}^*\ev(0)^* A_0} \left[\holo{\widehat{\sfa}}{0}{t_1}\cdot\widehat{\sfB_{\boldsymbol{\lambda}}}_{t_1}\cdots \right]\right\}+\\
&\phantom{=}+\tr_{\rho}(-1)^{2m} \pi_{\partial_{2m+1}*}\left[\holo{\widehat{\sfa}}{0}{t_1}\cdot \widehat{\sfB_{\boldsymbol{\lambda}}}_{t_1}\cdots \right] 
\end{aligned}\label{derextholo}
\end{equation}
We consider the first term on the right-hand side of (\ref{derextholo}); the Leibnitz rule for the dot product implies
\begin{align*}
\dd_{\pi_{2m+1}^*\ev(0)^*A_0}\left[\holo{\widehat{\sfa}}{0}{t_1}\cdot\widehat{\sfB_{\boldsymbol{\lambda}}}_{t_1}\cdots \right]&=\sum_{i=1}^{2m+1}\holo{\widehat{\sfa}}{0}{t_1}\cdot \widehat{\sfB_{\boldsymbol{\lambda}}}_{t_1}\cdots \dd_{\pi_{2m+1}^*\ev(0)^*A_0}\Big[\holo{\widehat{\sfa}}{t_{i-1}}{t_i}\widehat{\sfB_{\boldsymbol{\lambda}}}_{t_i}\Big]\cdots+ \\
&\phantom{=}+\holo{\widehat{\sfa}}{0}{t_1}\cdots \widehat{\sfB_{\boldsymbol{\lambda}}}_{t_{2m+1}}\cdot \dd_{\pi_{2m+1}^*\ev(0)^* A_0} \holo{\widehat{\sfa}}{t_{2m+1}}{1}.
\end{align*}
We recall that $\dd_{\ev(0)^* A_0}\holo{A_0}{0}{1}=0$ by (\ref{derpar}).
We compute the following expression
\begin{align*}
\dd_{\pi_{2m+1}^*\ev(0)^* A_0} \left[\holo{\widehat{\sfa}}{t_{i-1}}{t_i}\widehat{\sfB_{\boldsymbol{\lambda}}}_{t_i}\right]= \left[\dd_{\pi_{2m+1}^*\ev(0)^* A_0} \holo{\widehat{\sfa}}{t_{i-1}}{t_i}\right]\cdot \widehat{\sfB_{\boldsymbol{\lambda}}}_{t_i} + \holo{\widehat{\sfa}}{t_{i-1}}{t_i} \cdot \dd_{\pi_{2m+1}^*\ev(0)^* A_0} \widehat{\sfB_{\boldsymbol{\lambda}}}_{t_i};
\end{align*}
For the second term on the right-hand side of the above equation,
we obtain, repeating (almost) the same arguments used in the proof of Theorem~\ref{holonomy1},
$\holo{\widehat{\sfa}}{t_{i-1}}{t_i} \cdot \widehat{\dd_{A_0} \sfB_{\boldsymbol{\lambda}}}_{t_i}$;
for the first term, we obtain analogously
\begin{align*}
&\Big[\boldsymbol{\delta} \holo{\widehat{\sfa}}{t_{i-1}}{t_i}-\widehat{\sfa}_{t_{i-1}}\cdot\holo{\widehat{\sfa}}{t_{i-1}}{t_i} + \holo{\widehat{\sfa}}{t_{i-1}}{t_i} \cdot \widehat{\sfa}_{t_i}\Big]\cdot \widehat{\sfB_{\boldsymbol{\lambda}}}_{t_i}.
\end{align*}
Summing up all these contributions with the right signs and using repeatedly (\ref{prop-push}), we obtain, for the first term 
on the right-hand side of (\ref{derextholo}), the result
\begin{align*}
-\ev(0)^* \sfa\cdot \pi_{2m+1*} \Big[\holo{\widehat{\sfa}}{0}{t_1}\cdot \widehat{\sfB_{\boldsymbol{\lambda}}}_{t_1}\cdots \Big]&-\pi_{2m+1*} \Big[\holo{\widehat{\sfa}}{0}{t_1}\cdot \widehat{\sfB_{\boldsymbol{\lambda}}}_{t_1} \cdots \Big] \cdot \ev(0)^* \sfa \\
&+\boldsymbol{\delta} \Big\{\pi_{2m+1*}  \Big[\holo{\widehat{\sfa}}{0}{t_1}\cdot \widehat{\sfB_{\boldsymbol{\lambda}}}_{t_1}\cdots \Big]\Big\}.
\end{align*}
By the invariance of $\tr_\rho$, and since $\sfa$ and  $\pi_{2m+1*}  \Big[\holo{\widehat{\sfa}}{0}{t_1}\cdot \widehat{\sfB_{\boldsymbol{\lambda}}}_{t_1}\cdots \Big]$ have odd total degree, we get
\[
(-1)^{2m+1} \tr_\rho \pi_{2m+1*}\Big\{ \dd_{\pi_{2m+1}^*\ev(0)^* A_0}\Big[\holo{\widehat{\sfa}}{0}{t_1}\cdot \widehat{\sfB_{\boldsymbol{\lambda}}}_{t_1}\cdots \Big]\Big\}=\boldsymbol{\delta} \Big\{\tr_\rho \Big[\holo{\sfa}{0}{t_1}\cdot \widehat{\sfB_{\boldsymbol{\lambda}}}_{t_1}\cdots \Big]\Big\}.
\]

We now consider the second term on the right hand side of (\ref{derextholo}).
We recall the orientation choices for the $m$-simplex made in Appendix~\ref{app-hol}; with these in mind
we obtain (once again with the same arguments of the proof of Theorem~\ref{holonomy1})
\begin{multline}
-\tr_{\rho} \ev(0)^* \sfB_{\boldsymbol{\lambda}}\cdot \left\{\pi_{2m*}\left[\holo{\widehat{\sfa}}{0}{t_i}\cdot \widehat{\sfB}_{t_1}\cdots\right]\right\}+\\
+\tr_\rho \left\{ \pi_{2m*}\left[\holo{\widehat{\sfa}}{0}{t_1}\cdot \widehat{\sfB_{\boldsymbol{\lambda}}}_{t_1}\cdots \right]\cdot \ev(0)^{*} \sfB_{\boldsymbol{\lambda}} \right\}+\\ 
+\sum_{j=1}^{2m}(-1)^{2m+j+1}\tr_{\rho}\pi_{2m*}\left[\holo{\widehat{\sfa}}0{t_1}\cdot\widehat{\sfB}_{t_1}\cdots \widehat{(\sfB_{\boldsymbol{\lambda}}\cdot\sfB_{\boldsymbol{\lambda}})}_{t_i}\cdots \right].
\label{bordo-1}
\end{multline}
Since the trace is cyclic in the arguments and $\sfB_{\boldsymbol{\lambda}}$ has even total degree, the first two terms in the above expression cancel each other. In summary, we have obtained 
\begin{align*}
(-1)^{2m}& \tr_\rho \pi_{\partial_{2m+1}*}\Big[\holo{\widehat{\sfa}}{0}{t_1}\cdot \widehat{\sfB_{\boldsymbol{\lambda}}}_{t_1}\cdots \Big]=\\
&=\sum_{j=1}^{2m}(-1)^{2m+j+1} \tr_\rho \pi_{2m*}\Big[\holo{\widehat{\sfa}}{0}{t_1}\cdots\widehat{(\sfB_{\boldsymbol{\lambda}}\cdot \sfB_{\boldsymbol{\lambda}})}_{t_i}\cdots \Big].
\end{align*}
Recalling formulae (\ref{derivdiBlam}), we now apply $\bde_{\boldsymbol{\mu}}$ to $\holo{\widehat{\sfa}}{t_i}{t_{i+1}}$. We obtain by repeating (almost) the same arguments as in the proof of Theorem~\ref{holonomy1}  
\begin{align*}
\bde_{\boldsymbol{\mu}} \holo{\widehat{\sfa}}{t_i}{t_{i+1}}=-\underset{t_{i+1}\ge t\ge t_i}\int \holo{\widehat{\sfa}}{t_i}{t} \cdot \widehat{(\sfB_{\boldsymbol{\lambda}}\cdot \sfB_{\boldsymbol{\lambda}})}_t \cdot \holo{\widehat{\sfa}}{t}{t_{i+1}},
\end{align*}
with the same unifying notation of Remark~\ref{not-olovar}.
After repeated application of Lemma~\ref{commpull} and~\ref{pushcom}, 
we get the following 
result
\begin{align*}
\bde_{\boldsymbol{\mu}} \sfh_{2m-1,\rho}(\boldsymbol{\lambda};\sfA,\sfB)
&= \sum_{k=1}^{2m}(-1)^{2m+k+1} \tr_\rho \pi_{2m*} \Big[\holo{\widehat{\sfa}}{0}{t_1}\cdot \widehat{\sfB_{\boldsymbol{\lambda}}}_{t_1}\cdots \widehat{(\sfB_{\boldsymbol{\lambda}}\cdot \sfB_{\boldsymbol{\lambda}})}_{t_k} \cdots \Big]\\
&=(-1)^{2m} \tr_\rho \pi_{\partial_{2m+1}*}\Big[\holo{\widehat{\sfa}}{0}{t_1}\cdot \widehat{\sfB_{\boldsymbol{\lambda}}}_{t_1}\cdots \Big]\ ;
\end{align*} 
so the claim follows.
\end{proof}

\begin{Rem}
Observe that the statement of Theorem~\ref{pariolonomiagen} does not extend
to $\sfh_{2i,\rho}$. The problem in this case arises in (\ref{bordo-1})
in which the first two terms sum up instead of canceling each other.
This reflects what was already noted in subsection~\ref{class-oss}
about the classical versions of these observables in even dimensions. 
\end{Rem}

Since $\obsmu$ is a flat invariant observable, the results of 
subsection~\ref{ssec-twist} together with
Theorem~\ref{pariolonomiagen} imply
\begin{Cor} If $\boldsymbol{\mu}$ and $\boldsymbol{\lambda}$
are related by \eqref{mulambdaeven}, then 
\begin{equation}\label{ossBV-class}
(\sOBV-\dd)\left[\exp\left(\frac{\ii}{\hslash}\obs_{\boldsymbol{\mu}}\right)
\boldsymbol{\mathcal{H}}_{\rho}^{o}
(\boldsymbol{\lambda};\sfA,\sfB)\right]=0,
\end{equation}
under the assumptions of the regularization procedure on 
page~\pageref{conv-reg}.
\end{Cor}
We notice that this implies
Theorem~3 and Theorem~4 (for $M$ even-dimensional) 
of \cite{CR}.

\begin{Rem}
Let us finally note that, following the same arguments sketched in Remark~\ref{Rem-flat}, we may prove that the v.e.v.s of $\boldsymbol{\mathcal{H}}_{\rho}^{o}(\boldsymbol{\lambda};\sfA,\sfB)$ together with (the exponential of) the polynomial $\sfB$-terms depend only on the $\mathcal{G}$-equivalence class of the flat connection $A_0$.
\end{Rem}

\subsection{The $\DBV$-exactness of the polynomial 
observables}\label{expobs}
We end with a digression devoted to proving the identity
\begin{equation}\label{DBVex}
\obs_n=\DBV\left(\frac1n\obs_n \sfs\right), 
\end{equation}
where we have used the following notation:
\[
\sfs:=\int_M \dbraket{\sfa}\sfB.
\]
Of course, the functional $\sfs$ depends implicitly on a chosen background flat connection $A_0$, because the superfield $\sfa$ is seen as a supervariation of the superconnection $\sfA$, constructed via $A_0$; we do not indicate the dependence on $A_0$ in order to avoid cumbersome notation.
It is immediate to verify that $\sfs$ is an element of $\mathcal{S}$ with ghost number $-1$.
The validity of (\ref{DBVex}) relies on the important identity satisfied by the BV antibracket and by the BV Laplacian, namely the failure of the BV Laplacian $\DBV$ to satisfy the Leibnitz rule \eqref{fail-Lapl}.
We already know that, for all $n$, $\int_M \tr \sfB^n$ is $\DBV$-closed (since it does depend only on $\sfB$).
We want to prove separately the following identities:
\begin{align}\label{DBVex2}
\BV{\obs_n}{\sfs}=n\obs_n,\qquad \DBV \sfs=0.
\end{align}
If we assume the validity of the two previous identities, we can then derive 
(\ref{DBVex})  from (\ref{fail-Lapl}).
We begin with the first identity:
\begin{Thm}
The following identity holds
\begin{equation}
\BV{\obs_n}{\sfs}=n\obs_n
\end{equation}
for all $n\in \mathbb{N}$.
\end{Thm}
\begin{proof}
Since $\obs_n$ does not depend on $\sfa$, we have
\[
\BV{\obs_n}{\sfs}=\sbv{\obs_n}{\sfs}
=\int_M \dbraket{\rB{\obs_n}}{\lA{\sfs}}.
\]
We compute the right super functional derivative of $\sfs$ w.r.t.\ $\sfa$ 
getting 
\[
\lA{\sfs}=\sfB.
\]
The left super functional derivative w.r.t.\ $\sfB$ of $\obs_n$ reads 
\[
\rB{\obs_n}=\sfB^{n-1}.
\]
So it follows by 
\[
\int_M \dbraket{\rB{\obs_n}}{\lA{\sfs}}=\int_M \dbraket{\sfB^{n-1}}{\sfB}=\int_M \tr \sfB^{n-1}\cdot \sfB=\int_M \tr \sfB^n,
\]
that the claim is true.
\end{proof}

We want now to prove the second identity in (\ref{DBVex2}). 
Since $\braket{\ }{\ }$ is nondegenerate by assumption, we can find a basis $X^i$ of $\Lg$, $i=1,\dots,\dim \Lg$, which satisfies $\braket{X^i}{X^j}=\delta^{ij} \sigma_i$, where $\sigma_i=\pm 1$. We can then write 
\[
\phi^\alpha=\phi^\alpha_i X^i,\qquad \phi_\alpha^*=\phi_{\alpha\ j}^* X^j,
\]
where the coefficients $\phi^\alpha_i$ and $\phi_{\alpha\ i}^*$ are
forms on $M$ (of course, sum over repeated indices is understood here). 
By recalling the formulae defining the Hodge dual antifields and the definition of $\dbraket{\ }{\ }$ for forms with ghost number, we may write, despite of the dimension of $M$,
\begin{align*}
\sfs&=-\cirbr{c^*}c-\cirbr{a^*}a+\cirbr{B^*}B +\sum_{k=1}^{m-2} \cirbr{\tau_k^*}{\tau_k}=\\
\phantom{\sfs}&=\sum_{i=1}^{\dim \Lg} \sigma_i \left[-\cirbr{c^*_i}{c_i}-\cirbr{a^*_i}{a_i}+\cirbr{B^*_i}{B_i}+\sum_{k=1}^{m-2} \cirbr{\tau^*_i}{\tau_i}\right],
\end{align*}
where $\cirbr{\ }{\ }$ is defined in \eqref{cirbr}.
We now apply the BV Laplacian to the above expression, and we get the following result
\begin{align*}
\DBV \sfs&=\sum_{i=1}^{\dim \Lg} \sigma_i C \left[\binom{m}{m}-\binom{m}{m-1}+\binom{m}{m-2}-
\dots+(-1)^l \binom{m}{m-l}+\dots\right]=\\
& =\sum_{i=1}^{\dim \Lg}\sigma_i C \left(1-1\right)^m=0,
\end{align*}
where $C$ is an infinite constant (in fact, it is given by the Dirac $\delta$ distribution evaluated in $0$, multiplied by the volume of the manifold $M$). 
This argument is very similar to that used in the proof of the $\DBV$-closedness of the BV action for canonical $BF$ theories (see subsection~\ref{sec-dual}). The binomial coefficients take into account the number of components of $\phi^{\alpha}_i$ (recall that they are forms on $M$), while the signs come from the ghost numbers of the fields. In an appropriate regularization procedure, the above expression vanishes before one applies the regularization procedure.
So the claim follows.

\section{Generalizations}\label{sec-gen}
At the beginning of this paper, see page~\pageref{Ass1},
we have made three assumptions. As we have seen,
Assumption~\ref{Ass3} is necessary for the construction of loop observables
in odd dimensions
(one can get rid of it, if one just considers the action, see 
subsection~\ref{sec-dual}), and actually one needs the stronger
Assumption~\fullref{trodd} for the even-dimensional case (or for more
general loop observables).
Assumption~\ref{Ass1} is on the other hand needed if one wants to avoid
extra symmetries than those displayed in \eqref{gauge} or extra
reducibility than that described thereafter. Some modifications
are needed when $M$ is not compact, as we briefly describe
in subsection~\ref{ssec-Rm} for the case $M=\bbr^m$. 
Finally, Assumption~\ref{Ass2} is there only for the sake of simplicity.
We sketch in subsection~\ref{ssec-nontrivial}
some ideas for the generalization 
of the constructions in the paper if we get rid of it.

\subsection{The case $M=\bbr^m$}\label{ssec-Rm}
Here one has to require that the superfields should
decay sufficiently fast at infinity, the only flat connection is
the trivial one, and the infinitesimal gauge transformations \eqref{gauge}
have no extra reducibility than that considered in the paper.
Therefore, all the conclusions automatically apply to this case as well.

Let us only remark that in this case one may also consider observables
associated to paths with endpoints at infinity. This simplifies all
proofs in Sections~\ref{sec-gWl} and \ref{sec-lobs}, as one can disregard
the extremal boundary terms in the iterated integrals (which before
had to be proved to cancel each other by the cyclicity of the trace).
In particular, for $m$ even one needs no more the restriction on the definition
of the functional appearing in Theorem~\ref{pariolonomiagen}, as this
holds even if one sums on all the $\sfh$'s and not only on the odd ones.

\subsection{Nontrivial bundles}\label{ssec-nontrivial}
The main features of the paper that we have to generalize
are: $i)$ the BV formalism and $ii)$ the Wilson loops.

\subsubsection{The BV formalism}
The main problem of the generalization of the BV formalism to the case $P$ nontrivial arises because the fields of the $BF$ theories (classical fields, Faddeev--Popov ghost, ghosts for ghosts and associated BV antighosts) are no longer forms on $M$ with values in $\Lg$, but rather forms on $M$ taking values in the nontrivial bundle $\ad P$. We have therefore to generalize for such forms the notion of functional derivatives
of elements of $\mathcal{S}(\mathfrak{A})$ (for an $\mathbb{R}$-algebra
$\mathfrak{A}$).
They are easy to be carried out and have to
be considered as ``distributional forms'' on $M$ taking values in the bundle $\ad P\otimes \mathfrak{A}$.
The main difficulty lies in the generalization of functional derivatives for elements of $\mathcal{S}^*(N;\mathcal{E})$ (where now $\mathcal{E}\to N$ is a bundle over the manifold $N$). 
We consider them as elements of $\Omega^{*,*}(M\times N,\ad P\boxtimes \mathcal{E})$ where, by abuse of notation, we denote by $\Omega^{*,*}$ smooth forms as well as distributional forms (which is the main case) on $M\times N$.
The external tensor product $\boxtimes$ 
of two bundles on two different manifolds was already defined on
page~\fullref{A-NE}.
Then, the pairing induced by $\braket{\ }{\ }$ between the (pull-back of) forms on $M$ with values in $\ad P$ and the so defined functional derivatives gives forms on $M\times N$ with values in the pull-back of the bundle $\mathcal{E}$. It is a well-known fact (see for example~\cite{GHV1} and~\cite{GHV2}) that, given a (smooth) map $f\colon N_1\to N_2$ between two manifolds and a bundle $\mathcal{N}_2\to N_2$, forms on $N_1$ taking values in the pulled-back bundle $f^*\mathcal{N}_2\to N_1$ are generated by pull-backs of sections of $\mathcal{N}_2$ as an $\Omega^*(N_1)$-module. Usually, we consider the map $f$ to be a fibration, in order to perform a push-forward w.r.t.\ it.

The generalization of the notion of functional derivatives (as sketched above) leads to the generalization of the BV antibracket as well.
Analogously one can generalize the BV Laplacian, the super BV antibracket
and the super BV Laplacian. Then all constructions described in 
Sections~\ref{sec-BV}, \ref{sec-BVsuper} and \ref{sec-BVaction} hold
in the general case.

%We notice that a similar notion of functional derivatives was already (implicitly) done in~\cite{BGV} in the second chapter, where they defined the Heat Kernel for generalized Laplace operators on sections of some vector bundle on some Riemannian manifold $M$.

\subsubsection{The generalized Wilson loops}
We face the following problems:  we have to construct iterated integrals consistently with the fact that $\ad P$ is no more trivial. (We refer to Appendix~\ref{app-hol} for the main notations that we use in the next paragraphs.)
We sketch now this generalization which naturally works only in the case
when one restricts oneself to representations of $\Lg$ coming from representations of the corresponding Lie group $G$.

{}~

\noindent{\em The representation $\rho$}.
The Wilson loops defined in the paper depend explicitly by construction on some finite-dimensional representation $\rho:\Lg \to \End(V)$: in fact, we apply to the fields (which, in the case when $P$ is trivial, are forms on $M$ with values in $\Lg$) the representation $\rho$, and we obtain forms on $M$ with values in the algebra $\End(V)$.
This we can do no longer in the case $P$\ nontrivial; the forms are in fact 
elements of $\Omega^*(M,\ad P)$. We shall therefore construct a bundle morphism from $\ad P$ to some associated bundle, which must be related to the representation $\rho$. A natural way to do this consists in taking a representation $\widetilde{\rho}:G\to \Aut(V)$, for some finite-dimensional vector space $V$; this induces a representation $\widehat{\rho}$ on $\End(V)$ by conjugation. 
Therefore, we can define the associated bundle $\End_{P}(V):=P\times_{\widehat{\rho}} \End(V)$.
The derivative at the identity of $\widetilde\rho$ 
is an equivariant morphism $\rho$\ from $\Lg$ to $\End(V)$ and induces
a morphism from $\ad P$ to $\End_{P}(V)$ (which, by abuse of notation, we still denote by $\rho$). 

{}~

\noindent{\em The iterated integrals for the Wilson loops}.
 The generalized Wilson loops (in the case $P$ trivial) are constructed via iterated integrals that involve pull-backs of forms on $M$ with values in $\End(V)$ w.r.t.\ evaluation maps from $\lo M\times \triangle_n$ to $M$.
The construction of the generalized Wilson loops is based on the ``holonomy'' $\holo{A+a}01$ of the connection $A+a$ defined in Appendix~\ref{app-hol}.
The main object in the definition of $\holo{A+a}01$ is $\widehat{a}$; this was constructed by means of the evaluation map $\ev_1$ and by conjugation with $\holo{A}0\bullet$; by pulling back $\widehat{a}$ w.r.t.\ $\pi_{i,n}$ ($i=1,\dots,n$), we obtain a form on $\lo M\times \triangle_n$ with values in $\End(V)$.

In the nontrivial case, we work as follows: first, we take the image under $\rho$ of the form $a$, obtaining a form on $M$ with values in $\End_{P}(V)$. 
Next we take its pull-back w.r.t.\ the map $\ev_{1}$ 
and  obtain a form on $\lo M\times [0,1]$ taking values in $\End_{\ev_{1}^*P}(V)$ ($\ev_{1}^* P$ inherits the structure of a principal bundle over $\lo M\times [0,1]$ as a pull-back of $P$).
We then take a flat background connection $A_0$; by means of it, we can construct a $G$-equivariant isomorphism from $\ev_{1}^* P$ to $(\ev(0)\circ\pi_1)^* P$ that induces in turn an isomorphism between $\Omega^*(\lo M\times[0,1] ,\End_{\ev_{1}^* P}(V))$ and $\Omega^*(\lo M\times[0,1],\End_{(\ev(0)\circ\pi_1)^* P}(V))$. 
We still denote by 
$\widehat{a}\in\Omega^1(\lo M\times[0,1],\End_{(\ev(0)\circ\pi_1)^* P}(V))$ 
the result of these operations on the form $a$. Finally, we define
$\widehat{a}_{i,n}\in\Omega^1(\lo M\times\triangle_n,\End_{(\ev(0)\circ\pi_n)^* P}(V))$ as $\pi_{i,n}^*\widehat{a}$.
We can now multiply $\widehat{a}_{i,n}$ by $\widehat{a}_{j,n}$, for different $i, j\le n$, since all these forms take now values in the same algebra bundle $\End_{(\ev(0)\circ\pi_n)^* P}(V)\to \lo M\times \triangle_n$.
In this way we may define the generalization of all the functionals
appearing in Sections~\ref{sec-gWl} and \ref{sec-lobs}, and the related
theorems still hold.

The isomorphism described above depends 
explicitly on the connection $A_0$. The gauge group $\mathcal{G}$ operates (not necessarily freely) by ``conjugation'' on the set of equivariant bundle morphisms from $\ev_{1}^* P$ to $(\ev(0)\circ \pi_1)^* P$, which we denote by $\Hom_{G}(\ev_{1}^* P;(\ev(0)\circ \pi_1)^* P)$.
It is well-known that $\mathcal{G}$ operates on the space $\mathcal{A}(P)$ of connections on $P$, making it into a $\mathcal{G}$-principal bundle (modulo reducible connections and analytical technicalities, which we have skipped here). It can be proved that the isomorphism described above defines
an equivariant map from $\mathcal{A}(P)$ to $\Hom_{G}(\ev_{1}^* P;(\ev(0)\circ \pi_1)^* P)$.
This implies by construction that the v.e.v.s of the generalized Wilson loops
do not depend on the flat connection $A_0$, but rather on its equivalence class in $\left\{A\in \mathcal{A}(P):F_A=0\right\}/\mathcal{G}$ only. (This is in analogy with Remark~\ref{Rem-flat}.) 

\appendix
\section{Definition and main properties of the pushforward}\label{app-A}
Let $M$ be a manifold and  $\mathcal{E}\xrightarrow{\pi} M$ a smooth fiber bundle with typical fiber $F$, where $F$ is an oriented compact manifold possibly with
boundaries and corners.
Let $m$, resp.\ $e$, resp.\ $f$, denote the dimensions of $M$, resp.\ $\mathcal{E}$, resp.\ $F$ (so $e=f+m$). 

We pick a form $\omega$ in $\Omega^{p}(\mathcal{E})$, where $p\geq f$; we then define the pushforward $\pi_{*}\omega$ of the form $\omega$ w.r.t.\ $\pi$ 
as the form in $\Omega^{p-f}(M)$ which satisfies the following identity:
\begin{equation}
\int_{M} \pi_{*}\omega \wedge \eta = \int_{\mathcal{E}} \omega \wedge \pi^{*}\eta\quad,\ \forall \eta \in \Omega^{m+f-p}(M). 
\end{equation}
In the case $p<f$ we define $\pi_{*}\omega=0$.
We now list without proof the main properties of the push-forward:
\begin{equation}\label{prop-push}
\begin{array}{ll}
\pi_{*}(\pi^{*}\alpha \wedge \beta)=(-1)^{f \dg \alpha} \alpha \wedge \pi_{*} \beta, &\forall \alpha \in \Omega^{*}(M),\ \forall \beta \in \Omega^{*}(\mathcal{E}),\\
\pi_{*}(\alpha \wedge \pi^{*}\beta)=\pi_{*}\alpha\wedge \beta,&\forall \alpha \in \Omega^{*}(\mathcal{E}) ,\ \forall \beta \in \Omega^{*}(M) ,\\ 
\dd \pi_{*} \alpha = (-1)^{f} \pi_{*} \dd \alpha - (-1)^f \pi_{\partial *} \iota^{*}\alpha,&\forall \alpha \in \Omega^{*}(\mathcal{E}) ,
\end{array}
\end{equation}
where $\iota:\mathcal{E}_\partial \to \mathcal{E}$ is the canonical injection of the fiber bundle with typical fiber $\partial F$,
and $\pi_{\partial}:\iota(\mathcal{E}_{\partial})\to M$ is the 
corresponding projection.

Another important property which we use throughout the paper is given by the following Lemma (without proof).
We consider two manifolds $M$ and $N$ and  suppose that $\mathcal{E}\xrightarrow{\pi}  M$, resp.\ $\mathcal{F}\xrightarrow{\tilde\pi} N$, is a fiber bundle over the manifold $M$, resp.\ $N$. Let $\varphi:\mathcal{E} \to \mathcal{F}$ be a bundle morphism with base map $\psi:M\to N$.
We cast all these maps in the following commutative square:
\[
\begin{CD}
\mathcal{E}  @>\varphi>> \mathcal{F}\\
@V\pi VV              @VV\tilde{\pi}V\\
M          @>\psi>>                 N
\end{CD} 
\]
We suppose additionally that $\phi$ induces orientation preserving
diffeomorphisms of the fibers.
\begin{Lem}\label{commpull}
Under the above assumptions, the following identity holds:
\begin{equation}
(\pi_* \circ \varphi^*) \alpha = (\psi^* \circ \tilde{\pi}_*) \alpha, \quad\forall \alpha \in \Omega^p(\mathcal{F}).
\end{equation} 
\end{Lem}
For the proof, see \cite{GHV1}.

Let us suppose that we have a fiber bundle $\mathcal{E}\to \mathcal{F}$, and let additionally suppose that $\mathcal{F}\to M$ is a fiber bundle, too; let us denote the projections by $\pi_1$, resp. \ $\pi_2$.
\[
\begin{CD}
\mathcal{E} @>{\pi_1}>> \mathcal{F} \\
@.  @VV{\pi_2}V\\
 @. M
\end{CD}
\]
If we compose the two projections we obtain a fiber bundle $\mathcal{E}\to M$, with projection $\pi=\pi_2\circ \pi_1$, whose orientation will be determined by the orientation of the resulting fiber, the product manifold of the fibers of the two bundles.
Then  we obtain the following
\begin{Lem}\label{pushcom}
With the above hypotheses, the following identity holds
\begin{equation}
\pi_* \alpha = \pi_{2*} (\pi_{1*} \alpha), \quad\forall \alpha \in \Omega^p(\mathcal{E}). 
\end{equation}
\end{Lem}
This is just Fubini's Theorem for repeated integration, and the definition
of the push-forward is consistent with the orientation choices.

We end this Appendix by defining the push-forward of forms on $\mathcal{E}\to M$ with values in some finite dimensional vector space $W$. This is simply
given by
\[
\pi_{*} (\alpha \otimes \bvec):=\pi_{*} \alpha\otimes \bvec
\]
on generators and extended by linearity.

\section{Sign rules}\label{app-sign}
To introduce the dot product, let us for a moment suppose that we have a 
$\bbZ$-graded superalgebra $E$, 
and let us consider $\Omega^*(M; E)$ with differential
\[
\dd(\omega\otimes e):=\dd\omega\otimes e.
\]
Let us pick an element $\omega\otimes e$ in $\Omega^*(M; E)$; we can assign to it two gradings, namely its degree as a form on $M$ and the degree of its $E$-part; from now on, we will call the degree in $E$ ``ghost number.''
By ``homogeneous'' in $\Omega^*(M;E)$ we mean from now on any element $\alpha$ of given degree {\em{and}}\/ ghost number. 
We then define the product of homogeneous 
elements in $\Omega^*(M; E)$ by the rule
\[
(\omega\otimes e)\ (\omega'\otimes e'):=\omega\wedge\omega'\otimes ee'.
\]
The graded Leibnitz rule reads
\[
\dd(\alpha\ \beta)=\dd\alpha\ \beta + (-1)^{\dg\alpha}\alpha\ \dd\beta,
\quad \forall\alpha,\ \beta\in\Omega^*(M;E).
\]
In the case when $E$ is supercommutative, it also follows that
\[
\alpha\ \beta=(-1)^{\dg\alpha \dg\beta+\gh \alpha \gh\beta}\beta\ \alpha.
\]
In case $E$ is associative, we define the super Lie bracket 
of two homogeneous elements $a,b$ by
\[
\Lie{a}{b}:=a\ b-(-1)^{\gh a\gh b}b\ a,\ \forall a,b\in E;
\]
it satisfies the graded antisymmetry 
\[
\Lie{a}{b}=-(-1)^{\gh a\gh b}\Lie{b}{a}
\]
and the graded Jacobi identity
\[
\Lie{a}{\Lie{b}{c}}=\Lie{\Lie{a}{b}}{c}+(-1)^{\gh a \gh b}\Lie{b}{\Lie{a}{c}},
\]
for all homogeneous $a,b,c \in E$.
The super Lie bracket on $E$ can be extended to $\Omega^*(M; E)$ with the help of the wedge product by the rule
\[
\Lie{\alpha\otimes a}{\beta\otimes b}:=\alpha\wedge \beta\otimes \Lie{a}{b}.
\]

The graded antisymmetry and the graded Jacobi identity imply
\begin{itemize}
\item $\Lie{\alpha}{\beta}=-(-1)^{\dg \alpha\dg\beta+\gh\alpha\gh\beta}\Lie{\beta}{\alpha}$;
\item $\Lie{\alpha}{\Lie{\beta}{\gamma}}=\Lie{\Lie{\alpha}{\beta}}{\gamma}+(-1)^{\dg\alpha\dg \beta+\gh \alpha\gh\beta}\Lie{\beta}{\Lie{\alpha}{\gamma}}$,
\end{itemize}
for all homogeneous forms $\alpha,\beta, \gamma\in \Omega^*(M;E)$.

\begin{Rem}
It is possible to start directly with a super Lie algebra $\mathfrak{H}$
instead of $E$. 
The graded antisymmetry and the graded Jacobi identity in $\Omega^*(M;\mathfrak{H})$ hold as in the previous formulae.
\end{Rem}

\subsection{Dot products}
Since $\Omega^*(M;E)$ has two gradings, each homogeneous element $\alpha$ in the degree and in the ghost number inherits a new grading, the {\em{total degree}}, which is defined by $\abs{\alpha}:=\dg\alpha+\gh\alpha$.

With the help of the total degree, we can define the {\em{dot product}}\/ of two homogeneous forms $\alpha,\beta$ in $\Omega^*(M;E)$ by the rule
\[
\alpha\cdot\beta:=(-1)^{\gh\alpha\dg\beta}\alpha\ \beta,
\]
and accordingly the {\em{dot Lie bracket}}
\[
\lb{\alpha}{\beta}:=(-1)^{\gh\alpha\dg\beta}\Lie{\alpha}{\beta}.
\]
We now list some obvious properties:
Let us suppose that $E$ is supercommutative; then
\begin{align*}
\alpha\cdot\beta&=(-1)^{\abs{\alpha}\abs{\beta}}\beta\cdot\alpha,
&&\text{(graded commutativity).\hfill}\\
\intertext{For the dot Lie bracket holds in general}
\lb{\alpha}{\beta}&=-(-1)^{\abs{\alpha}\abs{\beta}}\lb{\beta}{\alpha}
&&\text{(graded antisymmetry),\hfill}\\
\lb{\alpha}{\lb{\beta}{\gamma}}&=\lb{\lb{\alpha}{\beta}}{\gamma}+(-1)^{\abs{\alpha}\abs{\beta}}\lb{\beta}{\lb{\alpha}{\gamma}},
&&\text{(graded Jacobi identity),\hfill}
\end{align*}
for all homogeneous forms $\alpha,\beta, \gamma$ in $\Omega^*(M;E)$.

Next, we notice that the exterior derivative satisfies the following 
graded Leibnitz rule
\[
\dd(\alpha\cdot\beta)=\dd \alpha\cdot \beta+(-1)^{\abs{\alpha}}\alpha\cdot \dd\beta.
\]

If we consider an (ungraded) algebra bundle (or more generally, a Lie algebra bundle) $\mathcal{B}\to M$, we can consider the space $\Omega^*(M,\mathcal{B})\otimes E$, instead of $\Omega^*(M;E)$; we define accordingly the total degree ($\mathcal{B}$ is ungraded and each fiber is an algebra) and the dot product (and the dot Lie bracket).

We next consider a covariant derivative $\dd_A$, coming from a connection $A$ on $\mathcal{B}$, and define its action on $\Omega^*(M,\mathcal{B})\otimes E$ by the rule
\[
\dd_A(\alpha\otimes a):=\dd_A\alpha\otimes a.
\] 
Then, the Leibnitz rule for $\dd_A$ w.r.t.\ the dot product and the dot Lie 
bracket follows easily.

\subsection{Superderivations}
We can also consider in this setting the BV operator $\delta$ defined by the BV action as a graded derivation on the superalgebra $E$, which we extend to 
$\Omega^*(M;E)$ by the rule
\[
\delta(\alpha\otimes a):=\alpha\otimes \delta a.
\]
It follows:
\begin{itemize} 
\item $\delta(\alpha\ \beta)=\delta\alpha\ \beta+(-1)^{\gh \alpha}\alpha\ \delta{\beta}$ for homogeneous $\alpha, \beta$ in $\Omega^*(M;E)$;
\item $\delta \circ \dd=\dd\circ \delta$ on $\Omega^*(M;E)$.
\end{itemize}
Let us next define $\sfdelta:=(-1)^{\dg}\,1\otimes \delta$, where 
$(-1)^{\dg}$ is the operator which multiplies each homogeneous form on $M$ 
by the parity of its degree. 

{}From its very definition, it follows 
\begin{itemize}
\item $\sfdelta(\alpha\cdot\beta)=\sfdelta\alpha\cdot \beta+(-1)^{\abs{\alpha}}\alpha\cdot\sfdelta{\beta}$  for homogeneous $\alpha, \beta$ in $\Omega^*(M;E)$;
\item $\sfdelta \circ \dd=-\dd \circ \sfdelta$ on $\Omega^*(M;E)$.
\end{itemize}

\begin{Rem}
The same identities can be proved even when we substitute $\Omega^*(M;E)$ with $\Omega^*(M,\mathcal{B})\otimes E$, for an ungraded algebra bundle $\mathcal{B}\to M$ and the exterior derivative with a covariant one, or if we replace $E$ by a super Lie algebra $\mathfrak{H}$. 
\end{Rem}
We can then define the operator 
\[
\mathcal{D}:=\dd\otimes 1+(-1)^{m+1}\,\sfdelta,\qquad
m=\dim M;
\]
it follows easily from all the above results that it is a 
superderivation w.r.t.\ the total degree.
Moreover, if $\delta$ is nilpotent, then so is $\mathcal{D}$, and consequently a differential on $\Omega^*(M;E)$. If we are dealing with $\Omega^*(M,\mathcal{B})\otimes E$, we can replace $\dd$ by a covariant derivative $\dd_A$ and
define $\mathcal{D}_A:=\dd_A\otimes 1+\sfdelta$, which
is then a superderivation. Moreover, if $A$ is flat, $\mathcal{D}_A$ is a 
superdifferential, too.
(Of course any linear combination of $\dd\otimes 1$ and
$\sfdelta$ has these properties. The conventional choice of the factor $(-1)^{m+1}$
is consistent with the choices made in the rest of the paper.)

In the paper, we also consider a flat background connection $A_0$ and its relative covariant derivative, along with a sum of forms, which we denote by $\sfa$, of total degree $1$.
Then $\dd_\sfA=\dd_{A_0}+\lb{\sfa}{\ }$ defines a superconnection on $\Omega^*(M,\ad P)$.
In the setting of this Appendix, this is tantamount to choosing forms on $\Omega^*(M,\ad P)\otimes E$ of total degree $1$; we sum all these forms and obtain a variation of the superconnection $A_0$.
We define accordingly $\mathcal{D}_{\sfA}:=\dd_{\sfA}+(-1)^{m+1}\sfdelta$; it is clear that $\mathcal{D}_{\sfA}$ is a derivation, and its curvature is given by
\[
\mathcal{D}_{\sfA}^2=
\lb{(-1)^{m+1}\sfdelta \sfa+\sfF_{\sfA}}{\ }=:
\lb{\mathcal{F}_\sfA}{\ };
\]
so (\ref{deltaa}) can be interpreted as the vanishing of the curvature 
$\mathcal{F}_\sfA$ of $\sfA$ on $\Omega^*(M,\ad P)\otimes E$; 
thus, $\sfA$ is formally ``superflat.''
Similarly, (\ref{deltab}) implies that the superform $\sfB$ (seen as an element of $\Omega^*(M,\ad P)\otimes E$ of total degree $m-2$) is  $\mathcal{D}_{\sfA}$-closed.

\subsection{Pullbacks and push-forwards}
Finally, let $\pi:\mathcal{E}\to M$ be a fiber bundle.
We then define the pullback, resp.\ push-forward, w.r.t.\ $\pi$ by the rules
\begin{itemize}
\item $\pi^*(\omega\otimes e):=\pi^*\omega\otimes e$, for $\omega\otimes e\in
\Omega^*(M;E)$;
\item $\pi_*(\eta\otimes e):=\pi_*\eta\otimes e$, for $\eta\otimes e\in
\Omega^*(\mathcal{E},E)$.
\end{itemize}
It follows immediately that 
\begin{itemize}
\item $\delta\circ \pi^*=\pi^*\circ \delta$;
\item $\delta\circ \pi_*=\pi_*\circ \delta$.
\end{itemize} 
It is then not difficult to verify that
\begin{itemize}
\item $\sfdelta\circ \pi^*=\pi^*\circ \sfdelta$;
\item $\sfdelta\circ \pi_*=(-1)^{\rk \mathcal{E}}\pi_*\circ \sfdelta$.
\end{itemize}
By definition of the dot product, it follows (in analogy with the first two equations in~(\ref{prop-push}))
\begin{itemize}
\item $\pi_{*}(\pi^{*}\alpha\cdot\beta)=(-1)^{\rk \mathcal{E} \abs{\alpha}} \alpha \cdot\pi_{*} \beta$;
\item $\pi_{*}(\alpha \cdot \pi^*\beta)=\alpha \cdot \pi_{*} \beta$.
\end{itemize}

\section{The universal global angular form}\label{app-gaf}
In this appendix we construct the universal global angular form
by using a fermionic integral representation.
This is analogous
to the construction of the Mathai--Quillen representative \cite{MQ}
of the Thom class  (see 
\cite{BGV} and \cite{CMR} and references therein).
Recall that global angular form on
a sphere bundle $\calS\xrightarrow{p} M$ is a form $\vartheta$
on $\calS$ satisfying $p_*\vartheta=1$ and $\dd\vartheta=-p^*e$,
where $e$ is a representative of the Euler class of the bundle.

Let $Q\to M$ be an $SO(n)$-principal bundle (not necessarily
$SO\left(M\right)$). Let $E$ the associated vector bundle $Q\times_{SO(n)}E_n$
with $E_n$ the $n$-dimensional Euclidean vector space 
We denote
by $\braket{\ }{\ }$ the corresponding scalar product.
Consider then the associated unit sphere bundle 
$\calS=Q\times_{SO(n)}S^{n-1}$ as the base manifold of
$\widehat{\calS}=Q\times S^{n-1}$. Let us summarize all these bundles and the respective projections in the following commutative square:
\[  
\begin{CD} 
Q  @< \widehat{p}  << Q\times S^{n-1}=\widehat{\calS}\\
@V \pi  VV              @V \widehat{\pi}  VV\\
M          @< p << Q\times_{SO(n)} S^{n-1}=\calS.
\end{CD}
\] 
We denote by $\theta$ a connection form on $Q$. By abuse of notation, we denote again by the same symbol its pull-back w.r.t.\ $\widehat{p}$ (which is again a connection on $\widehat{\calS}$), and by $F$ its
curvature. Last, let us denote by $x$ the canonical euclidean coordinates on $\bbr^n$
(with $S^{n-1}$ defined as the locus of $\braket xx=1$). We may consider $x$ as an equivariant function on $\widehat{S}$ with values in $\bbr^n$ (which inherits the canonical representation of $SO(n)$), and by
$\nabla x$ its corresponding covariant derivative, yielding a basic $1$-form on $\widehat{\calS}$ with values in $\bbr^n$. (Here, the right action of $SO(n)$ on $\widehat{\calS}$ is defined by $(q,x)O:=(qO,O^{-1}x)$.) From now, by {\em basic}\/ we will mean every form on $\widehat{\calS}$, which is horizontal and invariant w.r.t.\ the action of $SO(n)$.

Our aim is to construct a global angular form on the trivial
sphere bundle $\widehat{\calS}$ in terms of the monomials
\begin{multline}
\epsilon[x;F,k;\nabla x,l]\doteq
\epsilon_{a_1\dots a_{2k+l+1}}\ x^{a_1}\,
F^{a_2a_3}\dots F^{a_{2k}a_{2k+1}}\,\\
(\nabla x)^{a_{2k+2}}\dots (\nabla x)^{a_{2k+l+1}},
\label{defeps}
\end{multline}
where $2k+l+1=n$,
$\epsilon_{ij\dots n}$ is the totally antisymmetric tensor
and sums over repeated indices are understood. 
Observe that these monomials are basic in $\widehat{\calS}\to\calS$
since $x$, $\nabla x$ and $F$ are horizontal and equivariant.

Our first task is to write a generating function for these
monomials. To do so, we consider $\Pi T\bbr^n$. We go on denoting
by $x$ the (even) coordinates on the base and denote by
$\rho_i$, collectively $\rho$, the $n$ Grassmann coordinates on the fiber.
We introduce then the Berezin integration $\int[D\rho]$ by the rules:
\begin{itemize}
\item $\int[D\rho] P(\rho)=0$, for any polynomial $P$ in the odd variables $\rho_i$ of degree {\em strict}\/ less than $n$;
\item $\int[D\rho] \rho_1\cdots \rho_n=1$.
\end{itemize}
These two rules determine a unique Berezin integral on any polynomial in the Grassmann variables $\rho$ (any smooth function in the variables $\rho$ has the form of a polynomial of maximal degree $n$).

Then the generating function we are looking for reads
\begin{equation}
\Psi = \int[D\rho]\ \braket\rho x\ \exp S ,\label{defPsi}
\end{equation}
where
\begin{equation}
S = \braket\rho{\nabla x} + \frac\lambda2\braket\rho{F\rho},
\label{defS}
\end{equation}
and $\lambda$ is a parameter. For the next discussion
we need to introduce also the following generating function
of basic $n$-forms:
\begin{equation}
\Phi = \int[D\rho]\ \exp S.\label{defPhi}
\end{equation}

To prove that the forms generated by $\Phi$ and $\Psi$ are
actually basic just observe that 
the action of $SO(n)$ on $x$, $\nabla x$ and $F$
can be compensated for by a change of variables corresponding
to the fundamental representation of $SO(n)$ on the vector space
generated by $\{\rho^i\}$. 

\begin{Rem}
The Thom class on $P\times_{SO(n)} \bbR^n$ can be written as
a basic form on $P\times \bbR^n$ as
\[
U = \frac{(-1)^{
\left\lfloor\frac n2 \right\rfloor
}}{(2\pi t)^{n/2}}\ 
\int[D\rho]\ \exp{\left(
-\frac1{2t}\braket xx +
\braket\rho{\nabla x} - \frac t2\braket\rho{F\rho}
\right)},
\]
for any $t>0$ \cite{BGV}. 

So, apart form a mulplicative constant, $\Phi$ is
the restriction of $U|_{t=-\lambda}$ to $P\times S^{n-1}$,
while $\Psi$ is the restriction of the form obtained contracting
$U|_{t=-\lambda}$ with the radial vector field $r\,\frac\de{\de r}$.
\end{Rem}
Now we have the following
\begin{Lem}
$\Phi$ and $\Psi$ obey the equation:
\begin{equation}
d\Psi = (-1)^{n+1}\,\left(
n -2\,\lambda\,\frac\de{\de\lambda} +\frac1\lambda\right)\,\Phi.
\label{dPsi}
\end{equation}
\end{Lem}

\begin{proof}
When differentiating a form given as in \eqref{defPhi} or \eqref{defPsi},
we apply the following rules:
\begin{enumerate}
\item $\rho$ is odd with respect to exterior derivative;
\item $\rho$ behaves ``as if" it were covariantly closed.
\end{enumerate}

To justify the second rule, we first notice that, given any 
$n\times n$ matrix $X$, integration by parts shows that
\begin{equation}
\int[D\rho]\ \braket{X\,\rho}{\frac\de{\de\rho}}\,f =
\tr X\,\int[D\rho]\ f.
\label{tra}
\end{equation}
(With commuting variables we would have the same relation with a minus
sign on the r.h.s.)
As a consequence, 
\[
\int[D\rho]\ \delta f = 0,
\]
with
\begin{equation}
\delta f = \braket{-\theta \rho}{\frac{\de}{\de\rho}}f,
\label{deltarho}
\end{equation}
because $\theta$ takes values in $\mathfrak{so}(n)$.
Therefore,
\[
\dd\int[D\rho]\ f = (-1)^n\,\int[D\rho]\ \ddtilde f,
\]
where the new exterior derivative $\ddtilde$ is defined
by $\dd\pm\delta$. 
Introducing the covariant derivative
\[
\nablat= \ddtilde + \theta\cdot,
\]
we get from \eqref{deltarho} that $\nablat\rho=0$, that is, rule 2.
In particular, we have
\begin{gather*}
\ddtilde\braket\rho x = -\braket\rho{\nablat x} = 
-\braket\rho{\nabla x},\\
\ddtilde\braket\rho {\nabla x} = -\braket\rho{\nablat\nabla x} = 
-\braket\rho{F x},
\end{gather*}
since on $x$-variables $\nablat=\nabla$, and
\[
\ddtilde\braket\rho{F\rho} = 0,
\]
by the Bianchi identity.
Therefore,
\[
\dd\Psi = (-1)^{n+1}\,A + (-1)^n\, B,
\]
with
\begin{align*}
A &= \int[D\rho]\ \braket\rho{\nabla x}\,\exp S,\\
B &= \int[D\rho]\ \braket\rho x\,\braket\rho{F x}\,\exp S,
\end{align*}
and $S$ defined in \eqref{defS}.
Now some simple manipulations and the use of \eqref{tra} show that
\begin{multline}
A = \intrho \braket\rho{\frac\de{\de\rho}}\,\exp S
-\lambda\,\intrho \braket\rho{F\rho}\,\exp S=\\
=n\,\Phi - 2\,\lambda\,\frac\de{\de\lambda}\Phi.
\end{multline}
Similarly, we get
\begin{multline}
B = -\frac1\lambda\,\intrho\braket\rho x\,\braket x{\frac\de{\de\rho}}\,
\exp S +\\
+\frac1\lambda\,\intrho\braket\rho x\,\braket x{\nabla x}\,\exp S =\\
=-\frac1\lambda\,\intrho\left(\braket x{\frac\de{\de\rho}}\,
\braket\rho x\right)\,\exp S =
-\frac1\lambda\,\Phi,
\end{multline}
where we have used the constraint $\braket xx=1$ and the ensuing
identity 
\[
0=\dd\braket xx = 2\braket x{\nabla x}.
\]
\end{proof}
To exploit \eqref{dPsi}, it is convenient to expand $\Phi$ and 
$\Psi$ in powers of $\lambda$:
\begin{align*}
\Phi &= \sum_{k=0}^\infty \lambda^k\,\Phi_k,\\
\Psi &= \sum_{k=0}^\infty \lambda^k\,\Psi_k.
\end{align*}
Notice that these are actually finite sums. 
By performing the integrations we get
\begin{align}
\Phi_k &= \frac{(-1)^{k+\lfloor\frac n2\rfloor}}
{2^k\,k!\,(n-2k)!}\,
\epsilon[F,k;\nabla x,n-2k],\label{Phik}\\
\intertext{for $k=0,1,\dots,\lfloor\frac n2\rfloor$, and}
\Psi_k &= \frac{(-1)^{k+\lceil\frac{n-2}2\rceil}}
{2^k\,k!\,(n-2k-1)!}\,
\epsilon[x;F,k;\nabla x,n-2k-1],\label{Psik}
\end{align}
for $k=0,1,\dots,\lceil\frac{n-2}2\rceil$.
Applying \eqref{dPsi} to the power expansions, we get
\begin{align}
\Phi_0 &= 0,\label{Phi0}\\
\dd\Psi_k &= (-1)^{n+1}\,\left[
(n-2k)\,\Phi_k + \Phi_{k+1}\right].\label{dPsik}
\end{align}
Then we have the following
\begin{Lem}
The form
\begin{equation}
\overline{\vartheta}\doteq \sum_{k=0}^{s} C_k\,\Psi_k
\in\Omega^{n-1}_{\textup{basic}}(Q\times S^{n-1}),
\label{defeta}
\end{equation}
with $s=\lceil\frac{n-2}2\rceil$,
induces a global angular form $\vartheta$ on $S$ if and only if
the coefficients $C_k$ are defined by
\begin{equation}
C_k = \begin{cases}
(-1)^{k+s}\ \frac{(s-k)!}{2^{k+1}\,\pi^{s+1}} 
&\text{for $n=2s+2$}\\
(-1)^{k+s}\ \frac{(2s-2k)!}{2^{s-k+1}\,(2\pi)^s\,(s-k)!}
&\text{for $n=2s+1$}
\end{cases}
\label{defC}
\end{equation}
\end{Lem}

\begin{proof}
The forms $\overline{\vartheta}$ and $\vartheta$ are related by the formula 
\begin{equation}\label{eq-varth}
\overline{\vartheta}=\widehat{\pi}^*\vartheta. 
\end{equation}
The first property a global angular form has to satisfy
is $p_* \vartheta=1$. By the surjectivity of $\pi$ and by (\ref{eq-varth}), it suffices to show that $\widehat{p}_* \overline{\vartheta}=1$.
Since $\widehat{p}_*$ selects the $\theta$-independent
part in
\[
\Psi_0 = \frac{(-1)^s}
{(n-1)!}\ \epsilon_{i_1\dots i_n}\ 
x^{i_1}\,(\nabla x)^{i_2}\dots(\nabla x)^{i_n},
\]
this property is satisfied if and only if
we set the correct normalization:
\begin{equation}
C_0 = \frac{(-1)^s}{\Omega_{n-1}},\label{C0}
\end{equation}
where $\Omega_{n-1}$ is the volume of the unit $(n-1)$-sphere; that is,
\begin{align*}
\Omega_{2s+1} &= \frac{2\,\pi^{s+1}}{s!},\\
\Omega_{2s} &= \frac{2\,(2\pi)^s}{(2s-1)!!}.
\end{align*}
Next we use \eqref{Phi0} and \eqref{dPsik} to get
\[
\dd\overline{\vartheta} = (-1)^{n+1}\,\sum_{k=0}^{s}
\left[(n-2k)\,C_k+C_{k-1}\right]\,\Phi_k +
(-1)^{n+1}\,C_s\,\Phi_{s+1}.
\]
Now recall that the differential of a global angular form must
be basic on $\widehat{\calS}\to Q$ (in particular it has
to be the pullback w.r.t.\ $p$ of a representative of the Euler class). By (\ref{eq-varth}), together with the surjectivity of $\widehat{\pi}$, it is sufficient to show the identity $\dd \overline{\vartheta}=- \widehat{p}^* \pi^* e$, where $e$ is a representative of the Euler class.
All the $\Phi_k$ with $k\le s$ contain a form on $S^{n-1}$, so
they cannot be $\widehat{p}$-basic (i.e., $S^{n-1}$-independent). 
Therefore, we must choose the coefficients
$C_k$ so that the terms in square brackets vanish. This yields
a recursion rule that, once the initial condition is fixed by
\eqref{C0}, has the unique solution \eqref{defC}.

Now observe that the last term $\Phi_{s+1}$ vanishes when $n$
is odd. Therefore, $\vartheta$ is closed in this case, and this
is enough to prove that it is a global angular form.
If $n$ is even, however, 
\[
\Phi_{s+1} = \intrho\exp\left(
\frac12\braket\rho{F\rho}\right) =
\Pfaff F,
\]
with $\Pfaff$ denoting the Pfaffian, and the recursion fixes
\[
C_{s} = \frac1{(2\pi)^{s+1}}.
\]
As a consequence, in this case we get
\[
\dd\vartheta=
\frac{-1}{(2\pi)^{n/2}}\,\Pfaff F.
\]
Since the r.h.s.\ is minus
(a representative of the pullback to $Q\times S^{n-1}$ of) 
the Euler class, the lemma is proved.
\end{proof}

We can rewrite the results of the Lemma and \eqref{Psik} as follows.
In the odd-dimensional case, $n=2s+1$---cf.\ \cite{FHMM}---one has
\begin{equation}
\bar\eta = \frac1{2\,(4\pi)^s}\sum_{k=0}^s \frac1{k!\,(s-k)!}
\,\epsilon[x;F,k;\nabla x,2s-2k].
\label{betaodd}
\end{equation}
In the even-dimensional
case, $n=2s+2$, we get instead
\begin{equation}
\bar\eta = \frac1{2\,\pi^{s+1}}
\sum_{k=0}^s \frac1{4^k}\,
\frac{(s-k)!}{k!\,(2s-2k+1)!}\,
\epsilon[x;F,k;\nabla x,2s-2k+1].
\label{betaeven}
\end{equation}
Also observe that if one denotes by $T$ the antipodal map on the fiber
crossed with identity on the base, one has
\[
T^*\vartheta = (-1)^n\,\vartheta.
\]

\begin{Rem}
>From \eqref{betaodd}, we see that, {\em in the odd-dimensional
case}, $\overline{\vartheta}$ can also be given
the following expression:
\[
\overline{\vartheta} = \frac12\frac1{s!\,(4\pi)^s}
\int[D\rho]\ \braket\rho x\,\widetilde{S}^s =
\frac12\int[D\rho]\ \braket\rho x\,
\exp\left(\frac1{4\pi}\widetilde{S}\right),
\]
with
\[
\widetilde{S} = \braket\rho{F\,\rho}-\braket\rho{\nabla x}^2 =
\braket\rho{(F+\nabla x\,\nabla x)\,\rho}.
\]
This is in accordance with the interpretation given in \cite{BC}
of $\vartheta$ as one half of the Euler class of the
tangent bundle along the fiber $T_{S^{n-1}}\calS$.
\end{Rem}

\section{Parallel transport as a function on $\lo M$}\label{app-hol}
Let us consider a trivial principal bundle $P\to M$; so there exists a global section $\sigma:M\to P$.
Let us now pick a connection $A$ on $M$; we define the covariant derivative on $\Omega^* (M,\ad P)$ by the formula
\begin{equation}
\dd_{A} \mu:=\dd \mu +[\sigma^{*} A, \mu],
\end{equation} 
where $\mu$ is some section on $\ad P$, and $\sigma^{*} A$ is a $1$-form on $M$ with values in $\Lg$.
Since $\dd_{A} :\Gamma(M,\ad P)\to \Omega^1(M,\ad P)$ has all the properties of a covariant derivative, it can be easily extended to forms on $M$ with values in $\ad P$. We pick an element $a$ in $\Omega^1(M,\ad P)$, and we may define another connection starting from $A$, namely $\sigma^*A+a$ (which we write $A+a$). For the sake of simplicity, let us suppose that $A$ is flat.
Then the curvature of $A+a$ is given by
\[
F_{A+a}=\dd_A a+\frac{1}2 \Lie{a}{a}.
\]
We apply to $A+a$ the canonical injection $\iota$ from $\Lg$ to $\Ug$, so as to obtain a $1$-form on $M$ with values in $\Ug$; we omit to write $\iota$ before $A+a$.
Let us then define 
\[
\widehat{a}:=\holo{A}0\bullet \ev_1^*(a) \left(\holo{A}0\bullet\right)^{-1},
\]
where by $\ev_1$ we have denoted the evaluation map $\ev_1(\gamma;t):=\gamma(t)$, a map from $\lo M\times [0,1]$ to $M$, and by $\holo{A}0\bullet$ the (inverse of the) parallel transport w.r.t.\ the connection $A$, viewed as a function on $\lo M\times [0,1]$ with values in $\Aut \Ug$. 
\begin{Rem}\label{rem-hol}
We want to spend here some words on the definition of $\holo{A}0\bullet$ 
(even in the case when $P$ is nontrivial). 
We first consider the product of the pulled-back bundles
$(\ev(0)\circ \pi_1)^*P \to \lo M\times [0,1]$ and
$\ev_1^*P \to \lo M\times [0,1]$ which is then a $G\times G$-bundle
over $\lo M\times [0,1]\times\lo M\times [0,1]$.
Then we denote by 
$\mathcal{P}=(\ev(0)\circ \pi_1)^*P \times_\pi \ev_1^*P \to \lo M\times [0,1]$
the restriction of the product bundle to the diagonal of the base manifold.
We then consider the $G$-valued function $\widetilde{\holo{A}0\bullet}$
defined implicitly by the equation
\[
\widetilde\gamma_{p_1}(t)=p_2\,\widetilde{H(
A;p_1,p_2,\gamma)\arrowvert_{0}^{t}},\qquad
p_1,p_2\in P : \pi(p_1)=\gamma(0), \pi(p_2)=\gamma(t),
\]
where $\widetilde\gamma_{p_1}$ is the unique $A$-horizontal lift of
$\gamma$ starting at $p_1$. It is then clear that this function
is $G\times G$-equivariant if we define the action 
$\phi\colon(g,h;k)\mapsto hkg^{-1}$ of $G\times G$ on $G$.
So we can identify $\widetilde{\holo{A}0\bullet}$ with a section
$\holo{A}0\bullet$ of the associated bundle $\mathcal{P}\times_{\phi} G$.
(In the case when $P$ is trivial, $\holo{A}0\bullet$ can eventually be
identified with a map from $\lo M\times[0,1]$ to $G$.)
Consider next a finite-dimensional representation 
$\rho:G\to \Aut V$. This induces the action 
$\widetilde\rho\colon (g,h;\psi)\mapsto \rho(h)\psi\rho(g^{-1})$
of $G\times G$ on $\Aut V$. So we can also define ${\holo{A}0\bullet}_\rho$
as the corresponding section of the associated bundle 
$\mathcal{P}\times_{\widetilde\rho} \Aut V$. (In the case $P$ trivial,
this can then be identified with a map from $\lo M\times[0,1]$ to $\Aut V$.)
In particular, since $G$ operates on $\Lg$ by the adjoint action, 
it operates on $\Ug$ as well; so the above construction in this
case yields the $\Aut\Ug$-parallel transport with free final point 
${\holo{A}0\bullet}_{\Ad}$.
For the sake of simplicity, we always omit throughout the paper
the index referring to the representation as it is clear from the context.
\end{Rem}
Since $A$ is flat, $\holo{A}{0}{\bullet}$ enjoys the following useful property:
\begin{equation}\label{derpar}
\dd \holo{A}{0}{\bullet}=-\pi_1^*\ev(0)^*A\ \holo{A}{0}{\bullet}+\holo{A}{0}\bullet\ \ev_1^*A,
\end{equation}
where $\ev(0):\lo M\to M$ is defined by $\ev(0)(\gamma):=\gamma(0)$; we define further, for $n\in \mathbb{N}$, the maps $\pi_n:\lo M\times \triangle_n\to \lo M$ by $\pi_n(\gamma;t_1,\dots,t_n):=\gamma$, and by $\triangle_n$ we denote the $n$-simplex
\begin{equation}
\triangle_{n}:=\Big\{(t_1,\dots,t_n)\in [0,1]^{n}:0\leq t_1\leq\dots\leq t_n\leq1 \Big\},
\label{orsimpl}
\end{equation}
with orientation given by $\dd t_n\wedge\dots\wedge\dd t_1$.
By $\holo{A}01$ we denote the (inverse of the) holonomy along the loop $\gamma$, considered as a function on $\lo M$, taking values in $G$.
It follows from its definition that $\widehat{a}$ is a 1-form on $\lo M\times [0,1]$ with values in $\Ug$.
We can now define the {\emph{parallel transport}} w.r.t.\ $A+a$ from $1$ to $0$ as the formal series in $\Ug$
\begin{equation}
\holo{A+a}{0}{1}:=\holo{A}01+\sum_{n\geq 1} \pi_{n*} \left(\widehat{a}_{1,n}\wedge\cdots\wedge \widehat{a}_{n,n}\right)\holo{A}{0}1,
\end{equation}
where $\widehat{a}_{i,n}:=\pi_{i,n}^* \widehat{a}$ and $\pi_{i,n}(\gamma;t_1,\dots,t_n):=(\gamma;t_i)$.
It follows from its very definition that the parallel transport is an element of $\Omega^0(\lo M;\Ug)$.

\begin{Rem}\label{def-olovar}
We can define the parallel transport with free final point w.r.t.\ the connection $A+a$ by
\[
\holo{A+a}{0}{\bullet}:=1+\sum_{n\geq 1} \pi_{n+1,n+1*} \left(\widehat{a}_{1,n+1}\wedge \cdots \wedge \widehat{a}_{n,n+1}\right),
\]
with the same notations as above; it follows from its very definition that this particular parallel transport is a map $\lo M\times [0,1]\to \Ug$. The parallel transport as a function on $\lo M\times [0,1]$ with free initial point is defined analogously by the formula
\[
\holo{A+a}{\bullet}{1}:=1+\sum_{n\ge 1}\pi_{1,n+1*}\left(\widehat{a}_{2,n+1}\wedge\cdots\wedge \widehat{a}_{n+1,n+1}\right).
\]
Further, we can define the parallel transport with free end-points  
as a function on $\lo M\times \triangle_2$:
\[
\holo{{A+a}}{\bullet}{\bullet}:=1+\sum_{n\geq 1} \pi_{1,n,1*} \left(\widehat{a}_{2,n+2}\wedge \cdots \wedge \widehat{a}_{n+1,n+2}\right),
\]
where $\pi_{1,n,1}(\gamma;s_1,s_2,\dots,s_{n+1},s_{n+2}):=(\gamma;s_1,s_{n+2})$.
\end{Rem}

\begin{Thm}\label{holonomy1}
If we denote by $\dd_{\ev(0)^*A}$ the covariant derivative w.r.t.\ the connection $\ev(0)^*A$ on forms on $\lo M$ with values in $\Ug$, then the following identity holds, for any flat connection $A$ on $P$
and any $a\in\Omega^1(M,\ad P)$:
\begin{multline*}
\dd_{\ev(0)^* A} \holo{A+a}{0}{1}=-\ev(0)^*a \wedge \holo{A+a}{0}{1}
+\holo{A+a}{0}{1}\ev(0)^*a+\\
-\left[\underset{1\geq s\geq 0}\int\holo{A+a}{0}{s}\wedge \widehat{F_{A+a}}_s \wedge \holo{A+a}{s}{1}\right]\holo{A}01.
\end{multline*}
\end{Thm}

\begin{Rem}
We have written
\[
\underset{1\geq s\geq 0}\int\holo{A+a}{0}{s}\wedge \widehat{F_{A+a}}_s \wedge \holo{A+a}{s}{1}:=\pi_{1*}\left[\holo{A+a}{0}{\bullet}\wedge \widehat{F_{A+a}}\wedge \holo{A+a}{\bullet}{1}\right],
\]
where we have used again the notations in Remark~\ref{def-olovar}.
\end{Rem}

\begin{proof}
We shall apply Stokes Theorem to the push-forward w.r.t.\ the maps $\pi_{n}$; we note that the $n$-simplex $\triangle_{n}$ has a boundary, and that this boundary can be written as
\[
\partial\triangle_{n}=\bigcup_{\alpha=0}^n (\partial\triangle_n)_\alpha ,
\] 
where each $(\partial\triangle_n)_\alpha \cong \triangle_{n-1}$.
With our choice of orientation of the simplices---see after 
\eqref{orsimpl}---the first face of the boundary comes with opposite orientation, 
while the second has the right one, the third has opposite orientation again, 
and so forth: 
\[
\orient((\partial\triangle_n)_\alpha) =
(-1)^{\alpha+1}\orient(\triangle_{n-1}).
\]
We apply the covariant derivative w.r.t\ $A_0$ to the $n$-th term of the series, and we obtain:
\begin{equation}
\begin{aligned}
\dd_{\ev(0)^*A} {\pi}_{n*} \Big[\widehat{a}_{1,n} \wedge \cdots \wedge \widehat{a}_{n,n}\Big]\holo{A}{0}1=&(-1)^{n} \pi_{n*} \dd_{\pi_n^*\ev(0)^*A} \Big[\widehat{a}_{1,n} \wedge \cdots \wedge\widehat{a}_{n,n} \Big]\holo{A}{0}1+\\
&-(-1)^{n} \pi_{\partial_{n}*}\Big[ \iota_{\partial_{n}}^* \widehat{a}_{1,n}\wedge \cdots \wedge \widehat{a}_{n,n}\Big]\holo{A}{0}1, 
\end{aligned}\label{hol-1}
\end{equation}
where $\pi_{\partial_{n}}:\lo M\times \partial\triangle_{n}\to \lo M$ denotes the projection onto the first factor, while $\iota_{\partial_{n}}:\lo M\times \partial\triangle_{n}\to \lo M\times \triangle_{n}$ is the canonical injection of the boundary of the simplex into the simplex itself.
We have used implicitly the identity
\[
\dd_{\ev(0)^* A}\holo{A}{0}{1}=0,
\]
which follows from \eqref{derpar}.

We now consider the two terms on the right hand side of (\ref{hol-1}) separately, and we begin with the second term, which we call ``the $n$-th boundary term'' from now on.
Since $\partial\triangle_{n}=\bigcup_{\alpha=0}^n \partial\triangle_{n,\alpha}$, we can write
\[
\iota_{\partial_{n}}^* \left[\widehat{a}_{1,n}\wedge \cdots \wedge \widehat{a}_{n,n}\right]=
\sum_{\alpha=0}^{n} \iota_{\partial_{n,\alpha}}^* \Big[\widehat{a}_{1,n}\wedge \cdots \wedge \widehat{a}_{n,n}\Big],
\]
and $\iota_{\partial_{n,\alpha}}:\lo M\times (\partial\triangle_n)_\alpha\to \lo M\times \triangle_{n}$ is the canonical injection of the $\alpha$-th face of the boundary.
Considering the orientations of the faces, we obtain for the $n$-th boundary term the following expression:    \[
\sum_{\alpha=0}^n (-1)^{\alpha+1} \pi_{n-1*} \iota_{\partial_{n,\alpha}}^* \Big[\widehat{a}_{1,n}\wedge \cdots \wedge \widehat{a}_{n,n}\Big].
\] 
We begin with the first face $\alpha=0$; it is not difficult to prove the following identities
\begin{equation*}
\pi_{j,n} \circ \iota_{\partial_{n,0}}=\begin{cases}
\iota(0)\circ \pi_{n-1}& j=1,\\
\pi_{j-1,n-1}& j\neq 1;
\end{cases}
\end{equation*}
similarly, one shows for $\alpha=n$
\begin{equation*}
\pi_{j,n} \circ \iota_{\partial_{n,n}}=\begin{cases}
\pi_{j,n-1}& j\neq n,\\
\iota(1) \circ \pi_{n-1}& j=n.
\end{cases}
\end{equation*}
We have denoted by $\iota(0)$, resp.\ $\iota(1)$, the injection of $\lo M$ into $\lo M\times [0,1]$ given by $\iota(0)(\gamma):=(\gamma;0)$, resp.\ $\iota(1)(\gamma):=(\gamma;1)$.
for $\alpha\neq 0,n$, it holds
\begin{equation*}
\pi_{j,n} \circ \iota_{\partial_{n,\alpha}}=\begin{cases}
\pi_{j,n-1} & j<\alpha,\\
\pi_{\alpha,n-1}& j=\alpha,\alpha+1,\\
\pi_{j-1,n-1}& j>\alpha+1.
\end{cases}
\end{equation*}
It follows therefore ($\iota(0)\holo{A}0\bullet=1$ by its very definition)
\begin{align*}
&\iota_{\partial_{n,0}}^* \left[\widehat{a}_{1,n}\wedge \cdots \wedge \widehat{a}_{n,n}\right]=\pi_{n-1}^{*} \ev(0)^*a \wedge \widehat{a}_{1,n-1}\wedge \cdots \wedge \widehat{a}_{n-1,n-1};\\
&\iota_{\partial_{n,\alpha}}^* \left[\widehat{a}_{1,n}\wedge \cdots \wedge \widehat{a}_{n,n}\right]=\widehat{a}_{1,n-1}\wedge \cdots \wedge \widehat{(a\wedge a)}_{\alpha,n-1} \wedge \cdots \wedge \widehat{a}_{n-1,n-1} ;\\
&\iota_{\partial_{n,n}}^* \left[\widehat{a}_{1,n}\wedge \cdots \wedge \widehat{a}_{n,n}\right]=\widehat{a}_{1,n-1}\wedge \cdots \wedge \widehat{a}_{n-1,n-1}\wedge \pi_{n-1}^* \ev(0)^*a.
\end{align*}
We consider now the first term under the action of the push-forward w.r.t.\ $\pi_{n-1}$:
\begin{align*}
\pi_{n-1*} \iota_{\partial_{n,0}}^* \Big[\widehat{a}_{1,n}\wedge \cdots \wedge \widehat{a}_{n,n}\Big]\holo{A}01&=
\pi_{n-1*} \Big[\pi_{n-1}^{*} \ev(0)^*a \wedge \widehat{a}_{1,n-1}\wedge \cdots \wedge \widehat{a}_{n-1,n-1}\Big]\holo{A}01\\
&=(-1)^{n-1} \ev(0)^*a\wedge \pi_{n-1*} \Big[\widehat{a}_{1,n-1}\wedge\cdots\wedge \widehat{a}_{n-1,n-1}\Big]\holo{A}01.
\end{align*}
Similarly, we obtain for $\alpha=n$
\[
\pi_{n-1*}  \iota_{\partial_{n,n}}^* \Big[\widehat{a}_{1,n}\wedge \cdots \wedge \widehat{a}_{n,n}\Big]\holo{A}01 = \pi_{n-1*}\Big[\widehat{a}_{1,n-1}\wedge\cdots \wedge \widehat{a}_{n-1,n-1}\Big]\holo{A}01\wedge \ev(0)^*a,
\]
and for $\alpha\neq 0,n$ we obtain
\[
\pi_{n-1*} \iota_{\partial_{n,\alpha}}^* \Big[\widehat{a}_{1,n}\wedge \cdots \wedge \widehat{a}_{n,n}\Big]\holo{A}01=\pi_{n-1*}\Big[\widehat{a}_{1,n-1}\wedge \cdots \wedge (A \wedge A)_{\alpha,n-1} \wedge \cdots \wedge \widehat{a}_{n-1,n-1}\Big]\holo{A}01.
\]
Finally, we obtain the following expression for the $n$-th boundary term of (\ref{hol-1}):
\begin{align*}
&\sum_{\alpha=1}^{n-1} (-1)^{\alpha+1} \pi_{n-1*} \Big[\widehat{a}_{1,n-1}\wedge \cdots \wedge \widehat{(a\wedge a)}_{\alpha,n-1}\wedge \cdots \wedge \widehat{a}_{n-1,n-1}\Big]\holo{A}01 -(-1)^{n-1}\ev(0)^*a\wedge \\ 
&\quad \wedge \pi_{n-1*}\Big[\widehat{a}_{1,n-1}\wedge\cdots\wedge\widehat{a}_{n-1,n-1}\Big]\holo{A}01+(-1)^{n+1}\Big[\widehat{a}_{1,n-1}\wedge\cdots\wedge\widehat{a}_{n-1,n-1}\Big]\holo{A}01\wedge \ev(0)^*a.
\end{align*}
We then consider the first term of (\ref{hol-1}), and by the Leibnitz rule we obtain 
\[
\pi_{n*} \dd_{\pi_n^*\ev(0)^* A} \Big[\widehat{a}_{1,n}\wedge\cdots \wedge \widehat{a}_{n,n}\Big]=
\sum_{i=1}^n (-1)^{i+1} \pi_{n*} \Big[\widehat{a}_{1,n}\wedge\cdots\wedge\widehat{\dd_A a}_{i,n}\wedge \cdots \wedge \widehat{a}_{n,n}\Big];
\]
here we have used 
\[
\dd_{\pi_1^*\ev(0)^* A}\widehat{a}=\widehat{\dd_A a},
\]
which is a consequence of (\ref{derpar}).

Summing up all the two contributions to (\ref{hol-1}) with the correct signs, we obtain for the left hand side of (\ref{hol-1})
\begin{align*}
&\sum_{i=1}^n (-1)^{n+i+1} \pi_{n*} \Big[\widehat{a}_{1,n}\wedge\cdots\wedge  \widehat{\dd_A a}_{i,n}\wedge \cdots \wedge\widehat{a}_{n,n}\Big]\holo{A}01+\\ 
&+\sum_{\alpha=1}^{n-1} (-1)^{n+\alpha} \pi_{n-1*} \Big[\widehat{a}_{1,n-1}\wedge \cdots \wedge \widehat{(a\wedge a)}_{\alpha,n-1} \wedge \cdots \wedge \widehat{a}_{n-1,n-1}\Big]\holo{A}01-\\
&-\ev(0)^*a \wedge \pi_{n-1*} \Big[\widehat{a}_{1,n-1}\wedge \cdots \wedge  \widehat{a}_{n-1,n-1}\Big]\holo{A}01+\\
&+\Big[\widehat{a}_{1,n-1}\wedge \cdots \wedge \widehat{a}_{n-1,n-1}\Big]\holo{A}01\wedge \ev(0)^*a.
\end{align*}
We begin by summing up all the terms which contain before them $\ev(0)^*a$, and we obtain $-\ev(0)^*a\wedge \holo{A+a}{0}{1}$; similarly, by summing up all the terms which have $\ev(0)^*a$ on the right, we obtain $\holo{A+a}{0}{1}\wedge \ev(0)^*a$.
By recalling the definition of the curvature of the connection $A+a$, the sum of the remaining terms will give us\[
\sum_{n\geq1} \sum_{i=1}^n (-1)^{n+i+1} \pi_{n*}\Big[\widehat{a}_{1,n}\wedge \cdots \wedge \widehat{F_{A+a}}_{i,n}\wedge \cdots \wedge \widehat{a}_{n,n}\Big]\holo{A}01.
\]
For $1\le i\le n, n\geq 1$, we shall now write the projection $\pi_n$ as the composition of three projections, i.e.\ $\pi_n=\pi_{i,n}\circ \pi_{1,i-1,n}\circ \pi_{i+1,n,n}$, where the projections are defined as follows:
\begin{align*}
\pi_{i+1,n,n}:\ \lo M\times \triangle_n&\to \lo M\times \triangle_i; \\
\qquad (\gamma;s_1,\dots,s_n)&\mapsto (\gamma,s_1,\dots,s_i);\\
\pi_{1,i-1,n}:\ \lo M \times \triangle_i&\to \lo M\times [0;t] ; \\
\qquad (\gamma;s_1,\dots,s_i)&\mapsto (\gamma,s_i);\\
\pi_{i+1,n,n}:\ \lo M \times [0;1]&\to \lo M ; \\
\qquad (\gamma;s_i)&\mapsto \gamma.
\end{align*}
We notice for $j\le i$ the useful identity $\pi_{j,n} =\pi_{j,i}\circ \pi_{i+1,n,n}$, and for $j> i$ holds $\pi_{j,n}= \pi_{j-i,n-i}\circ \Bar{\pi}_{1,i,n}$, for $\Bar{\pi}_{1,i,n}(\gamma;s_1,\dots, s_n)=(\gamma;s_i,\dots,s_n)$. 
We then use the following identity (which follows from $\pi_n=\pi_{i,n}\circ \pi_{1,i-1,n}\circ \pi_{i+1,n,n}$ and Lemma~\ref{pushcom}):
\[
\pi_{n*}=(-1)^{(n-i)i} \pi_{i,n*} \circ \pi_{1,i-1,n*} \circ \pi_{i+1,n,n*};
\]
note the appearance of signs in the above identity: this is due to the fact that the three projections above do reverse the product orientation of the fiber of the trivial bundle over $\lo M$ given by the projection $\pi_n$. 
It is finally useful to introduce the commutative diagram
\[
\begin{CD}
\lo M\times \triangle_n  @>\Bar{\pi}_{1,i,n}>> \lo M\times \triangle_{1+n-i}\\
@V\pi_{i+1,n,n} VV              @VV\pi_{1+n-i}V\\
\lo M\times \triangle_i          @>\pi_{1,i-1,n}>>  \lo M\times [0,1];
\end{CD} 
\]
this diagram allows us to apply Lemma~\ref{commpull} to $\pi_{i+1,n,n*}\circ \Bar{\pi}_{1,i,n}^*$, in order to get the pullback w.r.t.\ $\pi_{1,i-1,n}$ before the pushforward w.r.t.\ $\pi_{1+n-i}$. We can then apply the first identity of (\ref{prop-push}), when we integrate w.r.t.\ $\pi_{1,i-1,n}$. 
We shall use once again such a commutative diagram, after integration w.r.t.\ $\pi_{1,i-1,n}$, along with (\ref{prop-push}) and Lemma~\ref{commpull}, in order to obtain the desired identity, accordingly to the notation introduced in Remark~\ref{def-olovar}.
\end{proof}
\begin{Rem}\label{not-olovar}
Similar identities can be proved for the two other cases in which we consider parallel transports as functions on $\lo M\times [0,1]$, resp.\ on $\lo M\times \triangle_2$.
We obtain for the first case the result 
\begin{align*}
\dd_{\pi^*\ev(0)^* A} \holo{A+a}{0}{\bullet}&=-\pi_1^*\ev(0)^*a \wedge \holo{A+a}{0}{\bullet}+\holo{A+a}{0}{\bullet} \wedge\widehat{a} -\\
&\phantom{=\ }-\pi_{2,2*}\left[\pi_{1,2}^*\left(\holo{A+a}0\bullet \wedge \widehat{F_{A+a}}\right) \wedge\holo{A+a}{\bullet}{\bullet}\right].
\end{align*}
An analogous identity holds for the holonomy as a function depending on the final point:
\begin{align*}
\dd_{\pi^*\ev(0)^* A} \holo{A+a}{\bullet}1&=-\widehat{a} \wedge \holo{A+a}\bullet1+\holo{A+a}\bullet1 \wedge\pi_1^*\left[\holo{A}01\ \ev(0)^*a\ \left(\holo{A}01\right)^{-1}\right] -\\
&\phantom{=\ }-\pi_{1,2*}\left[\holo{A+a}\bullet\bullet \wedge \pi_{2,2}^*\left(\widehat{F_{A+a}} \wedge\holo{A+a}{\bullet}1\right)\right].
\end{align*}
For the second case, we get
\begin{align*}
\dd_{\pi_2^*\ev(0)^* A} \holo{A+a}{\bullet}{\bullet}&=-\pi_{1,2}^*\widehat{a} \wedge \holo{A+a}{\bullet}{\bullet}+\holo{A+a}{\bullet}{\bullet} \wedge \pi_{2,2}^*\widehat{a} -\\
&\phantom{=\ }-\widetilde{\pi}_{2,3*}\left[\widetilde{\pi}_{3,3}^*\holo{A+a}{\bullet}{\bullet}\wedge\pi_{2,3}^*\widehat{F_{A+a}}\wedge\widetilde{\pi}_{1,3}^*\holo{A+a}{\bullet}{\bullet}\right],
\end{align*}
where $\widetilde{\pi}_{j,3}:\lo M\times \triangle_3\to \lo M\times \triangle_2$ forgets the $j$-th point of the $3$-simplex.
We have preferred to adopt the notation
\[
\underset{t_2\geq s\geq t_1}\int\holo{A+a}{t_1}{s}\wedge \widehat{F_{A+a}}_s \wedge \holo{A+a}{s}{t_2}
\]
for the third term in the three above expressions, where $t_1\le t_2$ can be fixed or can be understood as variables, given the case in the specific context.
\end{Rem}

\end{document}